\documentclass[11pt,a4paper]{article}

\usepackage[left=2cm,right=2cm,    top=2cm,bottom=2cm,bindingoffset=0cm]{geometry}

\usepackage{subfig}
\usepackage{flushend}
\usepackage{indentfirst}
\usepackage{graphics}
\usepackage{amsmath}
\usepackage{epstopdf}
\usepackage{float}
\usepackage{amssymb}
\usepackage[font=footnotesize,labelfont=bf]{caption}
\usepackage[labelsep=period]{caption}

\begin{document}

\newtheorem{ozn}{Definition}[section]
\newtheorem{thm}{Theorem}[section]
\newtheorem{prop}{Proposition}[section]
\newtheorem{nas}{Corollary}[section]
\newtheorem{zau}{Remark}[section]
\newtheorem{lema}{Lemma}[section]
\newtheorem{pry}{Example}[section]

\newcommand{\eps}{\varepsilon}
\newcommand{\me}{\mathbf}
\newcommand{\mr}{\mathbb}
\newcommand{\mt}{\mathsf}
\newcommand{\md}{\mathcal}
\newcommand{\ld}{\left}
\newcommand{\rd}{\right}
\newcommand{\kd}{\bigtriangleup}
\newcommand{\ip}{\int_{-\pi}^{\pi}}
\newcommand{\iii}{\int_{-\infty}^{\infty}}

\newcommand{\be}{\begin{equation}}
\newcommand{\ee}{\end{equation}}

\newcommand{\bem}{\begin{multline}}
\newcommand{\eem}{\end{multline}}

\newcommand{\bml}{\begin{multline*}}
\newcommand{\eml}{\end{multline*}}

\newcommand{\beg}{\begin{gather}}
\newcommand{\eeg}{\end{gather}}

\title{\textbf{Robust Forecasting of Sequences with Periodically Stationary Long Memory Multiplicative Seasonal Increments Observed with Noise and Cointegrated Sequences}}

\date{}

\maketitle

\noindent  STATISTICS, OPTIMIZATION AND INFORMATION COMPUTING\\
Statistics Opt. Inform. Comput., Vol. 10, March 2022, pp 295-338.\\
doi:10.19139/soic-2310-5070-1408

\vspace{20pt}

\author{\textbf{Maksym Luz}$^{1}$, \textbf{Mykhailo Moklyachuk}$^{2,*}$,\\\\
 {$^{1}$BNP Paribas Cardif, Kyiv, Ukraine, \\
$^{2,*}$Taras Shevchenko National University of Kyiv, Kyiv 01601, Ukraine}\\
$^{*}$Corresponding Author: Moklyachuk@gmail.com}\\\\\\

\noindent \textbf{Abstract.} \hspace{2pt}
The problem of optimal estimation of linear functionals
constructed from unobserved values of stochastic sequence with periodically stationary increments based on observations of the sequence with a periodically stationary noise is considered.
For sequences with known spectral densities, we obtain formulas for calculating values of the mean square errors and the spectral characteristics of the optimal estimates of the functionals.
Formulas that determine the least favorable spectral densities and minimax (robust) spectral
characteristics of the optimal linear estimates of functionals
are proposed in the case where spectral densities of the sequence
are not exactly known while some sets of admissible spectral densities are given.
\\

\noindent \textbf{Keywords:} \hspace{2pt} Periodically stationary sequence, SARFIMA, fractional integration, optimal linear estimate, mean square error, least favourable spectral density
matrix, minimax spectral characteristic\\

\noindent \textbf{AMS Subject Classification:} \hspace{2pt}
 Primary: 60G10, 60G25, 60G35, Secondary: 62M20, 93E10

\maketitle

\vspace{10pt}

\noindent {\bf{DOI}:} 10.19139/soic-2310-5070-1408

\section{Introduction}
\label{intro}

Non-stationary and long memory time series models are wildly used in different fields of economics, finance, climatology, air pollution, signal processing etc.
(see, for example, papers by Dudek and Hurd \cite{Dudek-Hurd}, Johansen and Nielsen  \cite{Johansen}, Reisen et al.\cite{Reisen2014}).
A core example -- a general multiplicative model, or $SARIMA (p, d, q)\times(P, D, Q)_s$ -- was introduced in the book by Box and Jenkins \cite{Box_Jenkins}. It includes both  integrated and  seasonal factors:
\be
 \Psi (B ^ s) \psi (B) (1-B) ^ d (1-B ^ s) ^ Dx_t = \Theta (B ^ s) \theta (B) \eps_t,
 \label{seasonal_3_model}
 \ee
 where $ \eps_t $, $ t \in \mr Z $, is a sequence of zero mean i.i.d. random variables,
and where $ \Psi (z) $ and $ \Theta (z) $ are two polynomials of degrees of $ P $ and $ Q $ respectively which have roots outside the unit circle.
The parameters $d$ and $D$ are allowed to be fractional. In the case where $|d+D|<1/2$ and $|D|<1/2$, the process  $\eqref{seasonal_3_model}$ is stationary and invertible.
The paper  by Porter-Hudak \cite{Porter-Hudak} illustrates an application of a seasonal ARFIMA model to the analysis of the monetary aggregates used by U.S. Federal Reserve.
Another model of fractional integration is GARMA processes  described by the equation (see Gray,  Cheng and  Woodward \cite{Gray})
\be
 (1-2uB+B^2) ^ dx_t = \eps_t,\quad |u|\leq1. \label{GARMA_model}
 \ee
For the resent results dedicated to the statistical inference for seasonal long-memory sequences,
we refer  to the paper by Tsai,  Rachinger and  Lin \cite{Tsai}, who developed methods of estimation of parameters in case of measurement errors.
In their paper Baillie,   Kongcharoen and  Kapetanios \cite{Baillie} compared MLE and semiparametric estimation procedures for prediction problems based on ARFIMA models.
 Based on simulation study, they indicate better performance of  MLE predictor than the one  based on the two-step local Whittle estimation.
Hassler and  Pohle \cite{Hassler} (see also Hassler \cite{Hassler_book}) assess a predictive  performance of various methods of  forecasting  of inflation and return volatility time series and show strong evidences for models with a fractional integration component.
One of the fields of interests related to time series analysis is optimal filtering. It aims to remove the unobserved components, such as trends, seasonality or noise signal, from the observed data  \cite{Ansley,Eiurridge}.

Another type of non-stationary processes are periodically correlated, or cyclostationary, processes introduced by Gladyshev \cite{Gladyshev}, which belong to the class of processes with time-dependent spectrum and are widely used in signal processing and communications
(see Gardner \cite{Gardner1994, Gardner2006}, Hurd and Miamee \cite{HurdMiamee},  Napolitano \cite{Napolitano} for a review of the works on cyclostationarity and its applications).
Periodic time series are considered as an extension of a SARIMA model (see Lund \cite{Lund} for a test assessing if a PARMA model is preferable to a SARMA one) and are suitable for forecasting stream flows with quarterly, monthly or weekly cycles (see Osborn \cite{Osborn}).
Baek, Davis and Pipiras \cite{Baek} have introduced a periodic dynamic factor model (PDFM) with periodic vector autoregressive (PVAR) factors, in contrast to seasonal VARIMA factors.
Basawa, Lund and Shao \cite{Basawa} have investigated first-order seasonal autoregressive processes with periodically varying parameters.

Methods of parameters estimations and filtering usually do not take into account the issues arising from real data, namely, the presence of outliers, measurement errors, incomplete information about the spectral, or model, structure etc. From this point of view, we see an increasing interest to robust methods of estimation that are reasonable in such cases (see
Reisen,  et al. \cite{Reisen2018}, Solci at al. \cite{Solci} for the examples of robust estimates of SARIMA and PAR models).
The paper by Grenander \cite{Grenander} should be marked as the first one where the minimax (robust) extrapolation problem for stationary processes was
formulated as a game of two players and solved.
Hosoya \cite{Hosoya}, Kassam  \cite{Kas1982}, Kassam and Poor \cite{Kassam_Poor}, Franke \cite{Franke1985}, Vastola and  Poor \cite{VastPoor1984}, Moklyachuk \cite{Moklyachuk,Moklyachuk2015},  Liu et al. \cite{Liu}  studied minimax (robust) extrapolation (prediction), interpolation (missing values estimation) and filtering (smoothing) problems for the stationary sequences and processes. Recent results of minimax extrapolation problems for stationary vector processes and periodically correlated processes belong to
Moklyachuk and Masyutka \cite{Moklyachuk:2008b,Mokl_Mas_pred,2012}
and  Moklyachuk and  Golichenko (Dubovets'ka) \cite{Dubovetska_filt,Dubovetska_extr,Moklyachuk2016} respectively.
Stationary sequences associated with a periodically correlated sequence are investigated by
Makagon et al.  \cite{ Makagon1999,Makagon2011}.
Processes with stationary and periodically stationary  increments are investigated by Luz and Moklyachuk  \cite{Luz_Mokl_extra,Luz_Mokl_book,Luz_Mokl_extra_GMI,Luz_Mokl_Bull}. We also mention works  by
Moklyachuk and Sidei \cite{Sidei_extr}, Moklyachuk, Masyutka and Sidei \cite{Sidei_book}, who derive minimax estimates for stationary processes from observations with missed values. Kozak, Luz and Moklyachuk \cite{Kozak_Mokl,Kozak_Luz_Mokl} studied the estimation problems for stochastic sequences with periodically stationary increments.

This article is dedicated to the robust forecasting problem for stochastic sequences with periodically stationary long memory multiple seasonal increments, or sequences with periodically stationary general multiplicative (GM) increments, introduced by Luz and Moklyachuk \cite{Luz_Mokl_extra_GMI}. Estimates of the unknown values of the
 sequence with periodically stationary GM  increments are based on observations of the sequence with the stationary noise sequence.

The article is organized as follows.
In Section $\ref{spectral_ theory}$, we recall definitions of
 generalized multiple (GM)  increment sequence $\chi_{\overline{\mu},\overline{s}}^{(d)}(\vec{\xi}(m))$ and stochastic
sequences $\xi(m)$ with periodically stationary (periodically correlated, cyclostationary) GM increments.
The spectral theory of vector-valued GM increment sequences is discussed.
Section $\ref{classical_extrapolation}$ deals with the classical forecasting   problem for the linear functionals  $A\xi $ which are constructed from unobserved values of the sequence $\xi(m)$ when the spectral densities of the sequence $\xi(m)$ and a noise  sequence $\eta(m)$ are known.
Estimates are obtained by applying the Hilbert space projection technique to the vector-valued sequence $\vec \xi(m)+ \vec \eta(m)$ with stationary GM  increments under the stationary noise sequence $\vec \eta(m)$ uncorrelated with $\vec \xi(m)$.
The case of non-stationary fractional integration is discussed as well.
Section $\ref{minimax_extrapolation}$ is dedicated to the minimax (robust) estimates in cases, where spectral densities of sequences are not exactly known
while some sets of admissible spectral densities are specified. We illustrate the proposed technique on the particular types of the sets, which are generalizations of the sets of admissible spectral densities described in a survey article by
Kassam and Poor \cite{Kassam_Poor} for stationary stochastic processes.

\section{Stochastic sequences with periodically stationary generalized multiple increments}\label{spectral_ theory}

\subsection{Definition and spectral representation of a periodically stationary GM increment}

In this section, we present definition, justification and a brief review of the spectral theory of stochastic sequences with periodically stationary multiple seasonal increments.
This type of stochastic sequences will allow us to deal with a wide range of non-stationarity in time series analysis.

 Consider a   stochastic sequence $\xi(m)$, $m\in\mathbb Z$ defined on a probability space $(\Omega, \cal F, \mt P)$.
 Denote by $B_{\mu}$ a backward shift operator   with the step $\mu\in
\mathbb Z$, such that $B_{\mu}\xi(m)=\xi(m-\mu)$; $B:=B_1$. Then $B_{\mu}^s=B_{\mu}B_{\mu}\cdots B_{\mu}$.

Define the incremental operator
\[
\chi_{\overline{\mu},\overline{s}}^{(d)}(B)
:=(1-B_{\mu_1}^{s_1})^{d_1}(1-B_{\mu_2}^{s_2})^{d_2}\cdots(1-B_{\mu_r}^{s_r})^{d_r}
=\sum_{k=0}^{n(\gamma)}e_{\gamma}(k)B^k,
\]
where
$d:=d_1+d_2+\ldots+d_r$, $\overline{d}=(d_1,d_2,\ldots,d_r)\in (\mr N^*)^r$,
 $\overline{s}=(s_1,s_2,\ldots,s_r)\in (\mr N^*)^r$
and $\overline{\mu}=(\mu_1,\mu_2,\ldots,\mu_r)\in (\mr N^*)^r$ or $\in (\mr Z\setminus\mr N)^r$, $n(\gamma):=\sum_{i=1}^r\mu_is_id_i$. Here $\mr N^*=\mr N\setminus\{0\}$.   The explicit formula for the coefficients $e_{\gamma}(k)$ is given in \cite{Luz_Mokl_extra_GMI}.

\begin{ozn}\label{def_multiplicative_Pryrist}
For a stochastic sequence $\xi(m)$, $m\in\mathbb Z$, the
sequence
\begin{eqnarray}
\nonumber
& &\chi_{\overline{\mu},\overline{s}}^{(d)}(\xi(m)):=\chi_{\overline{\mu},\overline{s}}^{(d)}(B)\xi(m)
=(1-B_{\mu_1}^{s_1})^{d_1}(1-B_{\mu_2}^{s_2})^{d_2}\cdots(1-B_{\mu_r}^{s_r})^{d_r}\xi(m)=
\\& &=\sum_{l_1=0}^{d_1}\ldots \sum_{l_r=0}^{d_r}(-1)^{l_1+\ldots+ l_r}{d_1 \choose l_1}\cdots{d_r \choose l_r}\xi(m-\mu_1s_1l_1-\cdots-\mu_rs_rl_r)
\label{GM_Pryrist}
\end{eqnarray}
is called a \emph{stochastic  generalized multiple (GM)  increment sequence} of differentiation   order
$d$
with a fixed seasonal  vector $\overline{s}\in (\mr N^*)^r$
and a varying step $\overline{\mu}\in (\mr N^*)^r$ or $\in (\mr Z\setminus\mr N)^r$.
\end{ozn}

\begin{ozn}
\label{oznStPryrostu}
A stochastic GM increment sequence $\chi_{\overline{\mu},\overline{s}}^{(d)}(\xi(m))$  is called   a wide sense
stationary if the mathematical expectations
\begin{eqnarray*}
\mt E\chi_{\overline{\mu},\overline{s}}^{(d)}(\xi(m_0))& = &c^{(d)}_{\overline{s}}(\overline{\mu}),
\\
\mt E\chi_{\overline{\mu}_1,\overline{s}}^{(d)}(\xi(m_0+m))\chi_{\overline{\mu}_2,\overline{s}}^{(d)}(\xi(m_0))
& = & D^{(d)}_{\overline{s}}(m;\overline{\mu}_1,\overline{\mu}_2)
\end{eqnarray*}
exist for all $m_0,m,\overline{\mu},\overline{\mu}_1,\overline{\mu}_2$ and do not depend on $m_0$.
The function $c^{(d)}_{\overline{s}}(\overline{\mu})$ is called a mean value  and the function $D^{(d)}_{\overline{s}}(m;\overline{\mu}_1,\overline{\mu}_2)$ is
called a structural function of the stationary GM increment sequence (of a stochastic sequence with stationary GM increments).
\\
The stochastic sequence $\xi(m)$, $m\in\mathbb   Z$
determining the stationary GM increment sequence
$\chi_{\overline{\mu},\overline{s}}^{(d)}(\xi(m))$ by   \eqref{GM_Pryrist} is called a stochastic
sequence with stationary GM increments (or GM increment sequence of order $d$).
\end{ozn}

\begin{zau}
For spectral properties of one-pattern increment sequence $\chi_{\mu,1}^{(n)}(\xi(m)):=\xi^{(n)}(m,\mu)=(1-B_{\mu})^n\xi(m)$
see, e.g., \cite{Luz_Mokl_book}, p. 1-8;  \cite{Yaglom}, p. 390--430.
The corresponding results for continuous time increment process $\xi^{(n)}(t,\tau)=(1-B_{\tau})^n\xi(t)$ are described in \cite{Yaglom:1955}, \cite{Yaglom}.
\end{zau}

\subsection{Definition and spectral representation of stochastic sequences with periodically stationary GM increment}

In this subsection, we present definition, justification and a brief review of the spectral theory of stochastic sequences with periodically stationary GM increments.

\begin{ozn}
\label{OznPeriodProc}
A stochastic sequence $\xi(m)$, $m\in\mathbb Z$ is called a stochastic sequence with periodically stationary (periodically correlated) GM increments with period $T$ if the mathematical expectations
\begin{eqnarray*}
\mt E\chi_{\overline{\mu},T\overline{s}}^{(d)}(\xi(m+T)) & = & \mt E\chi_{\overline{\mu},T\overline{s}}^{(d)}(\xi(m))=c^{(d)}_{T\overline{s}}(m,\overline{\mu}),
\\
\mt E\chi_{\overline{\mu}_1,T\overline{s}}^{(d)}(\xi(m+T))\chi_{\overline{\mu}_2,T\overline{s}}^{(d)}(\xi(k+T))
& = & D^{(d)}_{T\overline{s}}(m+T,k+T;\overline{\mu}_1,\overline{\mu}_2)
= D^{(d)}_{T\overline{s}}(m,k;\overline{\mu}_1,\overline{\mu}_2)
\end{eqnarray*}
exist for every  $m,k,\overline{\mu}_1,\overline{\mu}_2$ and  $T>0$ is the least integer for which these equalities hold.
\end{ozn}

Using  Definition \ref{OznPeriodProc}, one can directly check that the sequence
\begin{equation}
\label{PerehidXi}
\xi_{p}(m)=\xi(mT+p-1), \quad p=1,2,\dots,T; \quad m\in\mathbb Z
\end{equation}
forms a vector-valued sequence
$\vec{\xi}(m)=\left\{\xi_{p}(m)\right\}_{p=1,2,\dots,T}, m\in\mathbb Z$
with stationary GM increments by the relation
\[
\chi_{\overline{\mu},\overline{s}}^{(d)}(\xi_p(m))=\chi_{\overline{\mu},T\overline{s}}^{(d)}(\xi(mT+p-1)),\quad p=1,2,\dots,T,
\]
where $\chi_{\overline{\mu},\overline{s}}^{(d)}(\xi_p(m))$ is the GM increment of the $p$-th component of the vector-valued sequence $\vec{\xi}(m)$.

The spectral structure of the GM increment is described in the following theorem \cite{Karhunen}, \cite{Luz_Mokl_extra_GMI}.

\begin{thm}\label{thm1}
1. The mean value and the structural function
 of the vector-valued stochastic stationary
GM increment sequence $\chi_{\overline{\mu},\overline{s}}^{(d)}(\vec{\xi}(m))$ can be represented in the form
\begin{eqnarray}
\label{serFnaR_vec}
c^{(d)}_{ \overline{s}}(\overline{\mu})& = &c\prod_{i=1}^r\mu_i^{d_i},
\\
\label{strFnaR_vec}
 D^{(d)}_{\overline{s}}(m;\overline{\mu}_1,\overline{\mu}_2)& = &\int_{-\pi}^{\pi}e^{i\lambda
m} \chi_{\overline{\mu}_1}^{(d)}(e^{-i\lambda})\chi_{\overline{\mu}_2}^{(d)}(e^{i\lambda})\frac{1}
{|\beta^{(d)}(i\lambda)|^2}dF(\lambda),
\end{eqnarray}
where
\[\chi_{\overline{\mu}}^{(d)}(e^{-i\lambda})=\prod_{j=1}^r(1-e^{-i\lambda\mu_js_j})^{d_j}, \quad \beta^{(d)}(i\lambda)= \prod_{j=1}^r\prod_{k_j=-[s_j/2]}^{[s_j/2]}(i\lambda-2\pi i k_j/s_j)^{d_j},
\]
 $c$ is a vector, $F(\lambda)$ is the matrix-valued spectral function of the stationary stochastic sequence $\chi_{\overline{\mu},\overline{s}}^{(d)}(\vec{\xi}(m))$. The vector $c$
and the matrix-valued function $F(\lambda)$ are determined uniquely by the GM
increment sequence $ \chi_{\overline{\mu},\overline{s}}^{(d)}(\vec \xi(m))$.

2. The stationary GM increment sequence $\chi_{\overline{\mu},\overline{s}}^{(d)}(\vec{\xi}(m))$ admits the spectral representation
\begin{equation}
\label{SpectrPred_vec}
\chi_{\overline{\mu},\overline{s}}^{(d)}(\vec{\xi}(m))
=\int_{-\pi}^{\pi}e^{im\lambda}\chi_{\overline{\mu}}^{(d)}(e^{-i\lambda})\frac{1}{\beta^{(d)}(i\lambda)}d\vec{Z}_{\xi^{(d)}}(\lambda),
\end{equation}
where $d\vec{Z}_{\xi^{(d)}}(\lambda)=\{Z_{ p}(\lambda)\}_{p=1}^{T}$ is a (vector-valued) stochastic process with uncorrelated increments on $[-\pi,\pi)$ connected with the spectral function $F(\lambda)$ by
the relation
\[
 \mt E(Z_{p}(\lambda_2)-Z_{p}(\lambda_1))(\overline{ Z_{q}(\lambda_2)-Z_{q}(\lambda_1)})
 =F_{pq}(\lambda_2)-F_{pq}(\lambda_1),\]
 \[  -\pi\leq \lambda_1<\lambda_2<\pi,\quad p,q=1,2,\dots,T.
\]
\end{thm}

Consider another vector-valued stochastic sequence with the stationary GM
increments $\vec \zeta (m)=\vec \xi(m)+\vec \eta(m)$, where $\vec\eta(m)$ is a vector-valued stationary stochastic sequence, uncorrelated with $\vec\xi(m)$, with a spectral representation
\[
 \vec\eta(m)=\int_{-\pi}^{\pi}e^{i\lambda m} d\vec Z_{\eta}(\lambda),\]
where  $Z_{\eta}(\lambda)=\{Z_{\eta,p}(\lambda)\}_{p=1}^T$, $\lambda\in [-\pi,\pi)$, is a stochastic process with uncorrelated increments, that corresponds to the spectral function $G(\lambda)$ \cite{Hannan}.
The stochastic stationary GM increment $\chi_{\overline{\mu},\overline{s}}^{(d)}(\vec{\zeta}(m))$ allows the spectral representation
\begin{eqnarray*}
 \chi_{\overline{\mu},\overline{s}}^{(d)}(\vec{\zeta}(m))&=&\int_{-\pi}^{\pi}e^{i\lambda m}\frac{\chi_{\overline{\mu}}^{(d)}(e^{-i\lambda})}{\beta^{(d)}(i\lambda)}
 d\vec Z_{\xi^{(n)}}(\lambda)
  +\int_{-\pi}^{\pi}e^{i\lambda m}\chi_{\overline{\mu}}^{(d)}(e^{-i\lambda}) d\vec Z_{\eta }(\lambda),
 \end{eqnarray*}
while $d\vec Z_{\eta }(\lambda)=(\beta^{(d)}(i\lambda))^{-1} d\vec Z_{\eta^{(n)}}(\lambda)$,
$\lambda\in[-\pi,\pi)$. Therefore, in the case where the spectral functions $F(\lambda)$ and $G(\lambda)$ have the spectral density matrices $f(\lambda)=\{f_{ij}(\lambda)\}_{i,j=1}^T$ and $g(\lambda)=\{g_{ij}(\lambda)\}_{i,j=1}^T$, the spectral density matrix $p(\lambda)=\{p_{ij}(\lambda)\}_{i,j=1}^{T}$ of the stochastic sequence $\vec \zeta(m)$ is determined by the formula
\[
 p(\lambda)=f(\lambda)+|\beta^{(d)}(i\lambda)|^2g(\lambda).\]

\subsection{Moving average representation of  periodically stationary GM increment}
\label{moving_average}

Denote by $ H=L_2(\Omega, \cal F, \mt P)$ the Hilbert space of random variables $\zeta$ with zero first moment, ${\mt E}{\zeta}=0$, finite second moment, ${\mt E}|{\zeta}|^2<\infty$, endowed with the inner product $\langle \zeta,\eta\rangle={\mt E}{\zeta\overline{\eta}}$.
Denote by $H(\vec{\xi}^{(d)})$ the closed linear subspace of the space $H$ generated by
components
$\{\chi_{\overline{\mu},\overline{s}}^{(d)}(\xi_p(m)), p=1,\dots,T;\, m\in \mathbb Z \}$
of the stationary stochastic GM increment sequence $\vec{\xi}^{(d)}=\{\chi_{\overline{\mu},\overline{s}}^{(d)}(\xi_p(l))\}_{p=1}^{T}$, $\overline{\mu}>\overline{0}$,
and denote by $H^{q}(\vec{\xi}^{(d)})$ the closed linear subspace generated by components  $\{\chi_{\overline{\mu},\overline{s}}^{(d)}(\xi_p(m)),\, p=1,\dots,T;\, m\leqslant q \}$, $q\in \mathbb Z$.
Define a subspace
\[ S(\vec{\xi}^{(d)})=\bigcap_{q\in \mathbb{Z}} H^{q}(\vec{\xi}^{(d)})
\]
  of the Hilbert space $H(\vec{\xi}^{(d)})$. Then the  space $H(\vec{\xi}^{(d)})$ admits a decomposition
$ H(\vec{\xi}^{(d)})=S(\vec{\xi}^{(d)})\oplus R(\vec{\xi}^{(d)}) $
where $R(\vec{\xi}^{(d)})$ is the orthogonal complement of the subspace $S(\vec{\xi}^{(d)})$ in the space $H(\vec{\xi}^{(d)})$.
\begin{ozn}
A stationary (wide sense) stochastic GM increment sequence $\chi_{\overline{\mu},\overline{s}}^{(d)}(\vec{\xi}(m))=\{\chi_{\overline{\mu},\overline{s}}^{(d)}(\xi_p(m))\}_{p=1}^{T}$ is called regular if $H(\vec{\xi}^{(d)})=R(\vec{\xi}^{(d)})$,
and it is called singular if
$H(\vec{\xi}^{(d)})=S(\vec{\xi}^{(d)})$.
\end{ozn}

\begin{thm}
A stationary stochastic GM increment sequence $\chi_{\overline{\mu},\overline{s}}^{(d)}(\vec{\xi}(m))
=\{\chi_{\overline{\mu},\overline{s}}^{(d)}(\xi_p(m))\}_{p=1}^{T}$ is uniquely represented in the form
\begin{equation} \label{rozklad}
\chi_{\overline{\mu},\overline{s}}^{(d)}(\xi_p(m))
=\chi_{\overline{\mu},\overline{s}}^{(d)}(\xi_{S,p}(m))+\chi_{\overline{\mu},\overline{s}}^{(d)}(\xi_{R,p}(m))
\end{equation}
where  $\chi_{\overline{\mu},\overline{s}}^{(d)}(\xi_{R,p}(m)), p=1,\dots,T$ is a regular stationary GM   increment sequence and
$\chi_{\overline{\mu},\overline{s}}^{(d)}(\xi_{S,p}(m)), p=1,\dots,T$ is a singular stationary GM increment sequence.
The GM increment sequences
$\chi_{\overline{\mu},\overline{s}}^{(d)}(\xi_{R,p}(m)), p=1,\dots,T$  and
$\chi_{\overline{\mu},\overline{s}}^{(d)}(\xi_{S,p}(m)), p=1,\dots,T$
are
orthogonal for all $ m,k\in\mathbb{Z} $. They are defined by the formulas
\begin{eqnarray*} \chi_{\overline{\mu},\overline{s}}^{(d)}(\xi_{S,p}(m))&=&\mt E[\chi_{\overline{\mu},\overline{s}}^{(d)}(\xi_{p}(m))|S(\vec{\xi}^{(d)})],
\\
\chi_{\overline{\mu},\overline{s}}^{(d)}(\xi_{R,p}(m))&=&\chi_{\overline{\mu},\overline{s}}^{(d)}(\xi_{p}(m))
-\chi_{\overline{\mu},\overline{s}}^{(d)}(\xi_{S,p}(m)), \quad p=1,\dots,T.
\end{eqnarray*}
\end{thm}

Consider an innovation sequence  ${\vec\varepsilon(u)=\{\varepsilon_k(u)\}_{k=1}^q, u \in\mathbb Z}$ for a regular stationary GM
increment, namely, the sequence of uncorrelated random
variables such that $\mathsf{E} \varepsilon_k(u)\overline{\varepsilon}_j(v)=\delta_{kj}\delta_{uv}$,   $\mathsf{E} |\varepsilon_k(u)|^2=1, k,j=1,\dots,q; u \in\mathbb Z$, and $H^{r}(\vec\xi^{(d)} )=H^{r}(\vec\varepsilon)$ holds true for all $r \in \mathbb Z$, where  $H^r(\vec\varepsilon)$ is
the Hilbert space generated by elements $ \{ \varepsilon_k(u):k=1,\dots,q; u\leq r\}$,
 $\delta_{kj}$ and $\delta_{uv}$ are Kronecker symbols.

\begin{thm}\label{thm 4}
A   stationary GM increment sequence
$\chi_{\overline{\mu},\overline{s}}^{(d)}(\vec{\xi}(m))$ is regular if and only if there exists an
innovation sequence ${\vec\varepsilon(u)=\{\varepsilon_k(u)\}_{k=1}^q, u \in\mathbb Z}$
and a sequence of matrix-valued
functions $\varphi^{(d)}(k,\overline{\mu}) =\{\varphi^{(d)}_{ij}(k,\overline{\mu}) \}_{i=\overline{1,T}}^{j=\overline{1,q}}$, $k\geq0$, such that
\begin{equation}\label{odnostRuhSer}
\sum_{k=0}^{\infty}
\sum_{i=1}^{T}
\sum_{j=1}^{q}
|\varphi^{(d)}_{ij}(k,\overline{\mu})|^2
<\infty,\quad
\chi_{\overline{\mu},\overline{s}}^{(d)}(\vec{\xi}(m))=
\sum_{k=0}^{\infty}\varphi^{(d)}(k,\overline{\mu})\vec\varepsilon (m-k).
\end{equation}
Representation (\ref{odnostRuhSer}) is called the canonical
moving average representation of the stochastic stationary GM increment
sequence $\chi_{\overline{\mu},\overline{s}}^{(d)}(\vec{\xi}(m))$.
\end{thm}

The spectral function $F(\lambda)$ of  a stationary GM increment sequence $\chi_{\overline{\mu},\overline{s}}^{(d)}(\vec{\xi}(m))$
which admits the canonical representation
 $(\ref{odnostRuhSer})$  has the spectral density  $f(\lambda)=\{f_{ij}(\lambda)\}_{i,j=1}^T$ admitting the canonical
factorization
\begin{equation}\label{SpectrRozclad_f}
f(\lambda)=
\varphi(e^{-i\lambda})\varphi^*(e^{-i\lambda}),
\end{equation}
where the function
$\varphi(z)=\sum_{k=0}^{\infty}\varphi(k)z^k$ has
analytic in the unit circle $\{z:|z|\leq1\}$
components
$\varphi_{ij}(z)=\sum_{k=0}^{\infty}\varphi_{ij}(k)z^k; i=1,\dots,T; j=1,\dots,q$. Based on moving average  representation
$(\ref{odnostRuhSer})$ define
\[\varphi_{\overline{\mu}}(z)=\sum_{k=0}^{\infty} \varphi^{(d)}(k,\overline{\mu})z^k=\sum_{k=0}^{\infty}\varphi_{\overline{\mu}}(k)z^k.\]
 Then the following relation holds true:
\begin{equation}
\varphi_{\overline{\mu}}(e^{-i\lambda})
\varphi^*_{\overline{\mu}}(e^{-i\lambda})=
 \frac{|\chi_{\overline{\mu}}^{(d)}(e^{-i\lambda})|^2}{|\beta^{(d)}(i\lambda)|^2}f(\lambda)=
 \prod_{j=1}^r\frac{\ld|1-e^{-i\lambda\mu_js_j}\rd|^{2d_j}}{\prod_{k_j=-[s_j/2]}^{[s_j/2]}|\lambda-2\pi  k_j/s_j|^{2d_j}}f(\lambda).
\label{dd}
\end{equation}

The one-sided moving average representation (\ref{odnostRuhSer}) and
relation  (\ref{dd}) are used for finding the mean square optimal
estimates of unobserved values of vector-valued  sequences with stationary GM increments.

\subsection{Stochastic sequences with GM fractional increments}\label{fractional_interpolation}

Now we extend the definition of  the  GM increment sequence $\chi_{\overline{\mu},\overline{s}}^{(d)}(\vec{\xi}(m))$ of  the positive integer orders $(d_1,\ldots,d_r)$ to the fractional ones.
Within the subsection, we put the step $\overline{\mu}=(1,1,\ldots,1)$. Following the results of \cite{Luz_Mokl_extra_GMI}, represent the increment operator $\chi_{\overline{s}}^{(d)}(B)$  in the form
\be\label{FM_increment}
\chi_{\overline{s}}^{(R+D)}(B)=(1-B)^{R_0+D_0}\prod_{j=1}^r(1-B^{s_j})^{R_j+D_j},
\ee
where $(1-B)^{R_0+D_0}$ is an integrating component, $R_j$, $j=0,1,\ldots, r$, are non-negative integer numbers, $1<s_1<\ldots<s_r$.  Below we describe representations $d_j=R_j+D_j$, $j=0,1,\ldots, r$, of the increment orders $d_j$ by stating  conditions on the fractional parts $D_j$, such that the increment sequence $$\vec y(m):=(1-B)^{R_0}\prod_{j=1}^r(1-B^{s_j})^{R_i}\vec{\xi}(m)$$ is  a stationary  fractionally integrated seasonal stochastic  sequence.
For example, in case of single  increment pattern $(1-B^{s^*})^{R^*+D^*}$, this condition is $|D^*|<1/2$.

\begin{ozn}\label{def_fract_Pryrist}
A  sequence $\chi_{\overline{s}}^{(R+D)}(\vec \xi(m))$ is called  \emph{a fractional multiple (FM) increment sequence}.
\end{ozn}

Consider the generating function of the Gegenbauer polynomial:
\[
(1-2 u B+B^2)^{-d}=\sum_{n=0}^{\infty}C_n^{(d)}(u)B^n,
\quad
C_n^{(d)}(u)=\sum_{k=0}^{[n/2]}\frac{(-1)^k(2u)^{n-2k}\Gamma(d-k+n)}{k!(n-2k)!\Gamma(d)}.
\]

The following lemma and  theorem hold true \cite{Luz_Mokl_extra_GMI}.

\begin{lema}\label{frac_incr_2}
Define the sets $\md M_j=\{\nu_{k_j}=2\pi k_j/s_j: k_j=0,1,\ldots, [s_j/2]\}$, $j=0,1,\ldots, r$, and the set $\md M=\bigcup_{j=0}^r \md M_j$. Then the multiple seasonal increment operator admits the following representation:
\begin{eqnarray*}
\chi_{\overline{s}}^{(D)}(B)& := &(1-B)^{D_0}\prod_{j=1}^r(1-B^{s_j})^{D_j}
=\prod_{\nu \in \md M}(1-2\cos \nu B+B^2)^{\widetilde{D}_{\nu}}
\\
& = &(1-B)^{D_0+D_1+\ldots+D_r}(1+B)^{D_{\pi}}\prod_{\nu \in \md M\setminus\{0,\pi\}}(1-2\cos \nu B+B^2)^{D_{\nu}}
\\
& = &\ld(\sum_{m=0}^{\infty}G^+_{k^*}(m)B^m\rd)^{-1}=\sum_{m=0}^{\infty}G^-_{k^*}(m)B^m,
\end{eqnarray*}
where
\begin{eqnarray}
\label{Gegenbauer_GI+}
G^+_{k^*}(m)& = &\sum_{0\leq n_1,\ldots,n_{k^*}\leq m, n_1+\ldots+n_{k^*}=m}\prod_{\nu \in \md M}C_{n_{\nu}}^{(\widetilde{D}_{\nu})}(\cos\nu),
\\
\label{Gegenbauer_GI-}
G^-_{k^*}(m)& = &\sum_{0\leq n_1,\ldots,n_{k^*}\leq m, n_1+\ldots+n_{k^*}=m}\prod_{\nu \in \md M}C_{n_{\nu}}^{(-\widetilde{D}_{\nu})}(\cos\nu).
\end{eqnarray}
 $k^*=|\md M|$,  $D_{\nu}=\sum_{j=0}^rD_j \mr I \{\nu\in \md M_j\}$, $\widetilde{D}_{\nu}=D_{\nu}$ for $\nu \in \md M\setminus\{0,\pi\}$, $\widetilde{D}_{\nu}=D_{\nu}/2$ for $\nu=0$ and $\nu=\pi$.
\end{lema}

\begin{thm}\label{thm_frac}
Assume that for a stochastic vector-valued sequence $\vec \xi(m)$ and fractional differencing orders $d_j=R_j+D_j$, $j=0,1,\ldots, r$, the FM increment sequence $\chi_{\overline{1},\overline{s}}^{(R+D)}(\vec \xi(m))$ generated by increment operator (\ref{FM_increment})  is a stationary sequence with a bounded from zero and infinity spectral density $\widetilde{f}_{\overline{1}}(\lambda)$. Then for the non-negative integer numbers $R_j$, $j=0,1,\ldots, r$, the GM increment sequence $\chi_{\overline{1},\overline{s}}^{(R)}(\vec \xi(m))$    is stationary if $-1/2< D_{\nu}<1/2$ for all $\nu\in \md M$, where $D_{\nu}$ are defined by real numbers $D_j$, $j=0,1,\ldots, r$
in Lemma \ref{frac_incr_2} and it is long memory if $0< D_{\nu}<1/2$ for at least one $\nu\in \md M$, and invertible if $-1/2< D_{\nu}<0$. The spectral density $f(\lambda)$ of the stationary GM increment sequence $\chi_{\overline{1},\overline{s}}^{(R)}(\vec \xi(m))$ admits a representation
\[
 f(\lambda)=|\beta^{(R)}(i\lambda)|^2 \ld|\chi_{\overline{1}}^{(R)}(e^{-i\lambda})\rd|^{-2}\ld|\chi_{\overline{1}}^{(D)}(e^{-i\lambda})\rd|^{-2} \widetilde{f}_{\overline{1}}(\lambda)=:\ld|\chi_{\overline{1}}^{(D)}(e^{-i\lambda})\rd|^{-2} \widetilde{f} (\lambda),
  \]
  where
  \[
  \ld|\chi_{\overline{1}}^{(D)}(e^{-i\lambda})\rd|^{-2}=\ld|\sum_{m=0}^{\infty}G^+_{k^*}(m) e^{-i\lambda m}\rd|^2=\ld|\sum_{m=0}^{\infty}G^-_{k^*}(m) e^{-i\lambda m}\rd|^{-2}.
\]
\end{thm}

The further properties of the spectral density $f(\lambda)$ and the structural function $D^{(R)}_{ \overline{s}}(m,\overline{1},\overline{1})$ of a stationary GM increment vector sequence $\chi_{\overline{1},\overline{s}}^{(R)}(\vec \xi(m))$ as well as for examples of an application of Theorem \ref{thm_frac}, can be found in the works by  Palma and Bondon \cite{Palma-Bondon}, Giraitis and Leipus \cite{Giraitis}, Luz and Moklyachuk \cite{Luz_Mokl_extra_GMI}.

\section{Hilbert space projection method of forecasting}\label{classical_extrapolation}

\subsection{Forecasting of vector-valued stochastic sequences with stationary GM increments}

Consider a vector-valued stochastic sequence $\vec{\xi}(m)$ with stationary GM increments  constructed from the sequence $\xi(m)=\{\xi_{p}(m)\}_{p=1}^{T}$ with the help of transformation \eqref{PerehidXi}.
Let the stationary GM increment sequence $\chi_{\overline{\mu},\overline{s}}^{(d)}(\vec{\xi}(m))=\{\chi_{\overline{\mu},\overline{s}}^{(d)}(\xi_p(m))\}_{p=1}^{T}$
has an absolutely continuous spectral function $F(\lambda)$
and the spectral density $f(\lambda)=\{f_{ij}(\lambda)\}_{i,j=1}^{T}$.

Let  $\vec{\eta}(m)=\{\eta_{p}(m)\}_{p=1}^{T}$ be an uncorrelated with
the sequence $\vec\xi(m)$ stationary stochastic sequence with absolutely
continuous spectral function $G(\lambda)$ and spectral density
$g(\lambda)=\{g_{ij}(\lambda)\}_{i,j=1}^{T}$.

Without loss of generality assume that
 the mean
values of the increment sequence $ \chi_{\overline{\mu},\overline{s}}^{(d)}(\vec{\xi}(m))$ and the stationary
sequence  $\vec{\eta}(m)$ equal to $\vec 0$. We will also consider the increment step $\overline{\mu}>\overline{0}$.

\textbf{Extrapolation (forecasting) problem.} Consider the problem of mean square optimal linear estimation of the functionals
\begin{equation}
A\vec{\xi}=\sum_{k=0}^{\infty}(\vec{a}(k))^{\top}\vec{\xi}(k), \quad
A_{N}\vec{\xi}=\sum_{k=0}^{N}(\vec{a}(k))^{\top}\vec{\xi}(k),
\end{equation}
which depend on unobserved values of the stochastic sequence $\vec{\xi}(k)$ with stationary GM
increments.
Estimates are based on observations of the sequence $\vec\zeta(m)=\vec\xi(m)+\vec\eta(m)$ at points $m=-1,-2,\ldots$.

First of all we indicate some conditions which are necessary for solving the considered problem.  Assume that coefficients $\vec{a}(k)=\{a_{p}(k)\}_{p=1}^{T}$, $k\geq0$, and the linear transformation $D^{\mu}$ which is defined in the following part of the section satisfy the conditions
 \be\label{umovana a_e_d}
\sum_{k=0}^{\infty}\|\vec{a}(k)\|<\infty,\quad
\sum_{k=0}^{\infty}(k+1)\|\vec{a}(k)\|^{2}<\infty,\ee
\be\label{umovana a_mu_e_d}
\sum_{k=0}^{\infty}\|(D^{\mu}\me a)_k\|<\infty,\quad
\sum_{k=0}^{\infty}(k+1)\|(D^{\mu}\me a)_k\|^2<\infty.\ee

Assume also that spectral densities $f(\lambda)$ and $g(\lambda)$ satisfy the minimality condition
\be
 \ip \text{Tr}\left[ \frac{|\beta^{(d)}(i\lambda)|^2}{|\chi_{\overline{\mu}}^{(d)}(e^{-i\lambda})|^2}\ld(f(\lambda)+|\beta^{(d)}(i\lambda)|^2 g(\lambda)\rd)^{-1}\right]
 d\lambda<\infty.
\label{umova11_e_st.n_d}
\ee
This is the necessary and sufficient condition under which the mean square errors of the optimal estimates of the functionals $A\vec\xi$ and $A_N\vec\xi$ are not equal to $0$.

We apply the  Hilbert space estimation technique proposed by Kolmogorov \cite{Kolmogorov} which can be described as a $3$-stage procedure:
(i) define a target element (to be estimated) of the space $H=L_2(\Omega, \mathcal{F},\mt P)$ of random variables $\gamma$ which have zero mean values and finite variances, $\mt E\gamma=0$, $\mt E|\gamma|^2<\infty$, endowed with the inner product $\langle \gamma_1;\gamma_2\rangle={\mt E}{\gamma_1\overline{\gamma_2}}$,
(ii) define a subspace of $H$ generated by observations,
(iii) find an estimate of the target element as an orthogonal projection on the defined subspace.

\emph{Stage i}. The  functional $A\vec{\xi}$ does not  belong to the space $H=L_2(\Omega, \mathcal{F},\mt P)$.
With the help of the following lemma we  describe representations of the functional   as a sum of a functional with finite second moments  belonging to $H$ and a functional depending on the observed values of the sequence $\vec\zeta(k)$ (``initial values'') (for more details see \cite{Luz_Mokl_book, Luz_Mokl_extra_GMI}).

\begin{lema}\label{lema predst A}
The functional $A\vec\xi$ admits the representation
\be \label{zobrazh A_N_i_st.n_d}
    A\vec\xi=A\vec\zeta-A\vec\eta=H\vec\xi-V\vec\zeta,
\ee
where
\[
    H\vec\xi:=B\chi\vec\zeta-A\vec\eta,\]
    \[
    A\vec{\zeta}=\sum_{k=0}^{\infty}(\vec{a}(k))^{\top}\vec{\zeta}(k),\quad
    A\vec{\eta}=\sum_{k=0}^{\infty}(\vec{a}(k))^{\top}\vec{\eta}(k),\]
\[
B\chi\vec{\zeta}=\sum_{k=0}^{\infty}(\vec{b}_{N}(k))^{\top}\chi_{\overline{\mu},\overline{s}}^{(d)}(\vec{\zeta}(k)),
\quad
V\vec{\zeta}=\sum_{k=-n(\gamma)}^{-1} (\vec{v}_{N}(k))^{\top}\vec{\zeta}(k),
\]
the coefficients
$\vec{b}(k)=\{b_{p}(k)\}_{p=1}^{T}, k=0,1,\dots$ and
 $\vec{v}(k)=\{ v_{p}(k)\}_{p=1}^{T}, k=-1,-2,\dots,-n(\gamma)$
are calculated by the formulas
\begin{eqnarray}\label{koefv_N_diskr}
 \vec{v}(k)& = &\sum_{l=0}^{ k+n(\gamma)} \mt{diag}_T(e_{\nu}(l-k))\vec{b}(l), \, k=-1,-2,\dots,-n(\gamma),
\\
 \label{determ_b_N}
\vec{b}(k)& = &\sum_{m=k}^{\infty}\mt{diag}_T(d_{\overline{\mu}}(m-k))\vec{a}(m)=(D^{\overline{\mu}}{\me a})_k, \, k=0,1,\dots,
\end{eqnarray}
  $D^{\overline{\mu}}$ is the  linear transformation  determined by a matrix with the entries $(D^{\overline{\mu}})(k,j)=\mt{diag}_T(d_{\overline{\mu}}(j-k))$ if
$0\leq k\leq j$, and $(D^{\overline{\mu}})(k,j)=0$ if $0\leq j<k$, $\mt{diag}_T(x)$ denotes a $T\times T$ diagonal matrix with the entry $x$ on its diagonal, $\me a=((\vec{a}(0))^{\top},(\vec{a}(1))^{\top}, \ldots)^{\top}$, coefficients $\{d_{\overline{\mu}}(k):k\geq0\}$ are determined by the relationship
\[
 \sum_{k=0}^{\infty}d_{\overline{\mu}}(k)x^k
=\prod_{i=1}^r\ld(\sum_{j_i=0}^{\infty}x^{\mu_is_ij_i}\rd)^{d_i}.\]

\end{lema}

\begin{nas}\label{nas predst A_N}
The functional $A_N\vec\xi$ admits the representation
\be \label{zobrazh A_N_i_st.n_d}
    A_N\vec\xi=A_N\vec\zeta-A_N\vec\eta=H_N\vec\xi-V_N\vec\zeta,
\ee
where
\[
    H_N\vec\xi:=B_N\chi\vec\zeta-A_N\vec\eta,\]
    \[
    A_{N}\vec{\zeta}=\sum_{k=0}^{N}(\vec{a}(k))^{\top}\vec{\zeta}(k),\quad
    A_{N}\vec{\eta}=\sum_{k=0}^{N}(\vec{a}(k))^{\top}\vec{\eta}(k),\]
\[
B_{N}\chi\vec{\zeta}=\sum_{k=0}^{N}(\vec{b}_{N}(k))^{\top}\chi_{\overline{\mu},\overline{s}}^{(d)}(\vec{\zeta}(k)),
\quad
V_{N}\vec{\zeta}=\sum_{k=-n(\gamma)}^{-1} (\vec{v}_{N}(k))^{\top}\vec{\zeta}(k),
\]
the coefficients
$\vec{b}_{N}(k)=\{b_{N,p}(k)\}_{p=1}^{T}, k=0,1,\dots,N$ and
 $\vec{v}_{N}(k)=\{ v_{N,p}(k)\}_{p=1}^{T}, k=-1,-2,\dots,-n(\gamma)$
are calculated by the formulas
\begin{eqnarray}\label{koefv_N_diskr}
 \vec{v}_N(k)& = &\sum_{l=0}^{N \wedge k+n(\gamma)} \mt{diag}_T(e_{\nu}(l-k))\vec{b}_N(l), \, k=-1,-2,\dots,-n(\gamma),
\\
 \label{determ_b_N}
\vec{b}_N(k)& = &\sum_{m=k}^N\mt{diag}_T(d_{\overline{\mu}}(m-k))\vec{a}(m)=(D^{\overline{\mu}}_{N}{\me a}_{N})_k, \, k=0,1,\dots,N,
\end{eqnarray}
 $D^{\overline{\mu}}_{N}$ is the  linear transformation  determined by an infinite matrix with the entries\\ $(D^{\overline{\mu}}_{N})(k,j)=\mt{diag}_T(d_{\overline{\mu}}(j-k))$ if
$0\leq k\leq j\leq N$, and $(D^{\overline{\mu}}_{N})(k,j)=0$ if $j<k$ or $j,k>N$; \\
$\me a_N=((\vec{a}(0))^{\top},(\vec{a}(1))^{\top}, \ldots,(\vec{a}(N))^{\top},\vec 0 \ldots)^{\top}.$

\end{nas}

The functional $H\vec\xi$ from representation (\ref{zobrazh A_N_i_st.n_d}) has finite variance and the functional  $V\vec\zeta$ depends on the known observations of the stochastic sequence $\vec\zeta(k)$ at points $k=-n(\gamma),-n(\gamma)+1,\ldots,-1$. Therefore, estimates $\widehat{A}\vec\xi$ and $\widehat{H}\vec\xi$ of the functionals $A\vec\xi$
and $H\vec\xi$ and the mean-square errors
$\Delta(f,g;\widehat{A}\vec\xi)=\mt E |A\vec\xi-\widehat{A}\vec\xi|^2$ and $\Delta(f,g;\widehat{H}\vec\xi)=\mt E
|H\vec\xi-\widehat{H}\vec\xi|^2$ of the estimates $\widehat{A}\vec\xi$ and $\widehat{H}\vec\xi$ satisfy the following relations
\be\label{mainformula}
    \widehat{A}\vec\xi=\widehat{H}\vec\xi-V\vec\zeta,\ee
\be\label{mainformula2}
    \Delta(f,g;\widehat{A}\vec\xi)
    =\mt E |A\vec\xi-\widehat{A}\vec\xi|^2=
  \mt E|H\vec\xi-\widehat{H}\vec\xi|^2=\Delta(f,g;\widehat{H}\vec\xi).
    \ee
Therefore, the estimation problem for the functional $A\vec\xi$ is equivalent to the one for the functional $H\vec\xi$.
This problem can be solved by applying the Hilbert space projection method proposed by Kolmogorov \cite{Kolmogorov}.

The  functional $H\vec\xi$ admits the spectral representation
\begin{equation*}
H\vec\xi=
\int_{-\pi}^{\pi}
\left(\vec{B}_{\overline{\mu}}(e^{i\lambda})\right)^{\top}
\frac{\chi_{\overline{\mu}}^{(d)}(e^{-i\lambda})}{\beta^{(d)}(i\lambda)}
d\vec{Z}_{\xi^{(d)}+\eta^{(d)}}(\lambda)-
\int_{-\pi}^{\pi}\left(\vec{A}(e^{i\lambda})\right)^{\top}d\vec{Z}_{\eta}(\lambda),
\end{equation*}
where
\begin{equation*}
\vec{B}_{\overline{\mu}}(e^{i\lambda})=\sum_{k=0}^{\infty}\vec{b}_{\overline{\mu}}(k)e^{i\lambda k}
=\sum_{k=0}^{\infty}(D^{\overline{\mu}}\me a)_ke^{i\lambda k},\quad
 \vec{A}(e^{i\lambda})=\sum_{k=0}^{\infty}\vec{a}(k)e^{i\lambda k}.
\end{equation*}

\emph{Stage (ii).} Introduce the following notations. Denote by
$H^{0-}(\xi^{(d)}_{\overline{\mu},\overline{s}}+\eta^{(d)}_{\overline{\mu},\overline{s}})$ the closed linear subspace generated by values
$\{\chi_{\overline{\mu},\overline{s}}^{(d)}(\vec{\xi}(k))+\chi_{\overline{\mu},\overline{s}}^{(d)}(\vec{\eta}(k)):k=-1,-2,-3,\dots\}$, $\overline{\mu}>\vec 0$
of the observed GM increments
in the Hilbert space $H=L_2(\Omega,\mathcal{F},\mt P)$ of random variables $\gamma$ with zero mean value, $\mt E\gamma=0$, finite variance, $\mt E|\gamma|^2<\infty$, and the inner product $(\gamma_1,\gamma_2)=\mt E\gamma_1\overline{\gamma_2}$.

Denote by $L_2^{0-}(f(\lambda)+|{\beta^{(d)}(i\lambda)}|^2 g(\lambda))$  the closed linear subspace of the Hilbert space
$L_2(f(\lambda)+|{\beta^{(d)}(i\lambda)}|^2 g(\lambda))$  of vector-valued functions with the inner product $\langle g_1,g_2\rangle=\ip (g_1(\lambda))^{\top}(f(\lambda)+|{\beta^{(d)}(i\lambda)}|^2 g(\lambda))\overline{g_2(\lambda)}d\lambda$ which is generated by the functions
\[
 e^{i\lambda k}\chi_{\overline{\mu}}^{(d)}(e^{-i\lambda})\frac{1}{\beta^{(d)}(i\lambda)}\vec\delta_l,\quad \vec\delta_l=\{\delta_{lp}\}_{p=1}^T,\,
 l=1,\dots,T;\,\,\, k \leq -1,
\]
where $\delta_{lp}$ are Kronecker symbols.

The representation
\[
 \chi_{\overline{\mu},\overline{s}}^{(d)}(\vec{\xi}(k))+\chi_{\overline{\mu},\overline{s}}^{(d)}(\vec{\eta}(k))=\ip e^{i\lambda
 k}\chi_{\overline{\mu}}^{(d)}(e^{-i\lambda})\dfrac{1}{\beta^{(d)}(i\lambda)}dZ_{\xi^{(d)}+\eta^{(d)}}(\lambda)\]
yields a one-to-one correspondence  between elements
\[
e^{i\lambda k}\chi_{\overline{\mu}}^{(d)}(e^{-i\lambda})\dfrac{1}{\beta^{(d)}(i\lambda)}\vec \delta_l\]
from the space
$L_2^{0-}(f(\lambda)+|{\beta^{(d)}(i\lambda)}|^2 g(\lambda))    $
and elements
$\chi_{\overline{\mu},\overline{s}}^{(d)}(\vec{\xi}(k))+\chi_{\overline{\mu},\overline{s}}^{(d)}(\vec{\eta}(k))$
from the space $
H^{0-}(\xi^{(d)}_{\overline{\mu},\overline{s}}+\eta^{(d)}_{\overline{\mu},\overline{s}}).
$

Relation (\ref{mainformula}) implies that every linear estimate $\widehat{A}\vec\xi$ of the functional $A\vec\xi$
can be represented in the form
\be \label{otsinka A_e_d}
 \widehat{A}\vec\xi=
 \ip
(\vec{h}_{\overline{\mu}}(\lambda))^{\top}d\vec{Z}_{\xi^{(d)}+\eta^{(d)}}(\lambda)-
\sum_{k=-n(\gamma)}^{-1}(\vec v_{\mu}(k))^{\top}(\vec\xi(k)+\vec\eta(k)),
\ee
 where
$\vec{h}_{\overline{\mu}}(\lambda)=\{h_{p}(\lambda)\}_{p=1}^{T}$ is the spectral characteristic of the optimal estimate $\widehat{H}\vec\xi$.

\emph{Stage (iii).}
At this stage we find the mean square optimal estimate
$\widehat{H}\vec\xi$ as a projection of the element $H\vec\xi$ on the
subspace
$H^{0-}(\xi^{(d)}_{\overline{\mu},\overline{s}}+\eta^{(d)}_{\overline{\mu},\overline{s}}) $.  This projection is
determined by two conditions:

1) $ \widehat{H}\vec\xi\in H^{0-}(\xi^{(d)}_{\overline{\mu},\overline{s}}+\eta^{(d)}_{\overline{\mu},\overline{s}}) $;

2) $(H\vec\xi-\widehat{H}\vec\xi)
\perp
H^{0-}(\xi^{(d)}_{\overline{\mu},\overline{s}}+\eta^{(d)}_{\overline{\mu},\overline{s}}) $.

The second condition implies the following relation which holds true for all $k\leq-1$
\begin{multline*}
\int_{-\pi}^{\pi}
\bigg[\bigg(
\vec{B}_{\overline{\mu}}(e^{i\lambda})^{\top}
\frac{\chi_{\overline{\mu}}^{(d)}(e^{-i\lambda})}{\beta^{(d)}(i\lambda)}
-
\vec{h}_{\overline{\mu}}(\lambda)\bigg)^{\top}
\ld(f(\lambda)+|{\beta^{(d)}(i\lambda)}|^2 g(\lambda)\rd)
-
\\
-
(\vec{A}(e^{i\lambda}))^{\top}g(\lambda){\overline{\beta^{(d)}(i\lambda)}}\bigg]
\frac{\overline{
\chi_{\overline{\mu}}^{(d)}(e^{-i\lambda})}}
{\overline{\beta^{(d)}(i\lambda)}}
e^{-i\lambda k}d\lambda=\vec 0.
 \end{multline*}

\noindent This relation allows us to derive the spectral characteristic
$\vec{h}_{\overline{\mu}}(\lambda)$ of the estimate $\widehat{H}\vec\xi$ which can be represented in the form
\begin{multline} \label{spectr A}
(\vec{h}_{\overline{\mu}}(\lambda))^{\top}=
(\vec{B}_{\overline{\mu}}(e^{i\lambda}))^{\top}
\frac{\chi_{\overline{\mu}}^{(d)}(e^{-i\lambda})}{\beta^{(d)}(i\lambda)}
-
(\vec{A}_{\overline{\mu}}(e^{i\lambda}))^{\top}g(\lambda)\frac{\overline{\beta^{(d)}(i\lambda)}}
{\overline{
\chi_{\overline{\mu}}^{(d)}(e^{-i\lambda})}}
\ld(f(\lambda)+|{\beta^{(d)}(i\lambda)}|^2 g(\lambda)\rd)^{-1}
-
\\
-(\vec{C}_{\overline{\mu}}(e^{i\lambda}))^{\top}
\frac{\overline{\beta^{(d)}(i\lambda)}}
{\overline{
\chi_{\overline{\mu}}^{(d)}(e^{-i\lambda})}}
\ld(f(\lambda)+|{\beta^{(d)}(i\lambda)}|^2 g(\lambda)\rd)^{-1},
\end{multline}
where
\[
\vec{A}_{\overline{\mu}}(e^{i\lambda})=\vec{A}(e^{i\lambda})\overline{
\chi_{\overline{\mu}}^{(d)}(e^{-i\lambda})}=\sum_{k=0}^{\infty}\vec{a}_{\overline{\mu}}(k)e^{ik\lambda},
\]
\begin{equation} \label{coeff a_mu}
 \vec a_{\overline{\mu}}(m)=\sum_{l=\max\ld\{m-n(\gamma),0\rd\}}^{m}e_{\gamma}(m-l)\vec a(l),\quad  m\geq 0,
 \end{equation}
\[
\vec{C}_{\overline{\mu}}(e^{i \lambda})=\sum_{k=0}^{\infty}\vec{c}_{\overline{\mu}}(k)e^{ik\lambda},\]
coefficients $\vec{c}_{\overline{\mu}}(k)=\{c_{\overline{\mu},p}(k)\}_{p=1}^T, k=0,1,\dots,$  are unknown and have to be found.

It follows from condition 1) that the following equations should be satisfied  for $j\geq 0$
\begin{multline} \label{eq_C}
\int_{-\pi}^{\pi} \biggl[(\vec{B}_{\overline{\mu}}(e^{i\lambda}))^{\top}-
(\vec{A}_{\overline{\mu}}(e^{i\lambda}))^{\top}
\frac{|{\beta^{(d)}(i\lambda)}|^2}{|\chi_{\overline{\mu}}^{(d)}(e^{-i\lambda})|^2}g(\lambda)
\ld(f(\lambda)+|{\beta^{(d)}(i\lambda)}|^2 g(\lambda)\rd)^{-1}
-
\\
-
(\vec{C}_{\overline{\mu}}(e^{i\lambda}))^{\top}
\frac{|{\beta^{(d)}(i\lambda)}|^2}{|\chi_{\overline{\mu}}^{(d)}(e^{-i\lambda})|^2}
\ld(f(\lambda)+|{\beta^{(d)}(i\lambda)}|^2 g(\lambda)\rd)^{-1}
\biggr]e^{-ij\lambda}d\lambda=0.
\end{multline}

 Define for $ k,j\geq 0$ the Fourier coefficients of the corresponding functions
\[
T^{\overline{\mu}}_{k,j}=\frac{1}{2\pi}\int_{-\pi}^{\pi}
e^{-i\lambda(k-j)}
 \frac{|\beta^{(d)}(i\lambda)|^2}{|\chi_{\overline{\mu}}^{(d)}(e^{-i\lambda})|^2}\ld[g(\lambda)
\ld(f(\lambda)+|{\beta^{(d)}(i\lambda)}|^2 g(\lambda)\rd)^{-1}\rd]^{\top}
d\lambda;
\]
\[
P_{k,j}^{\overline{\mu}}=\frac{1}{2\pi}\int_{-\pi}^{\pi} e^{-i\lambda(k-j)}
 \frac{|\beta^{(d)}(i\lambda)|^2}{|\chi_{\overline{\mu}}^{(d)}(e^{-i\lambda})|^2}
\ld[\ld(f(\lambda)+|{\beta^{(d)}(i\lambda)}|^2 g(\lambda)\rd)^{-1}\rd]^{\top}
d\lambda;
\]
 \[
 Q_{k,j}=\frac{1}{2\pi}\int_{-\pi}^{\pi}
e^{-i\lambda(k-j)}\ld[f(\lambda)\ld(f(\lambda)+|{\beta^{(d)}(i\lambda)}|^2 g(\lambda)\rd)^{-1}g(\lambda)\rd]^{\top}
d\lambda.
\]

  \noindent Making use of the defined Fourier coefficients, relation \eqref{eq_C} can be presented as a system  linear equations
  \begin{equation} \label{linear equations1}
    \vec{b}_{\overline{\mu}}(j)-\sum_{ m=0}^{\infty}T^{\overline{\mu}}_{j,m}\vec{a}_{\overline{\mu}}(m)
    =\sum_{k=0}^{\infty}P_{j,k}^{\overline{\mu}}\vec{c}_{\overline{\mu}}(k),\,\, j \geq 0,
    \end{equation}
determining the unknown coefficients ${\vec c}_{\overline{\mu}}(k)$, $k\geq0$.
This  system of equations can be written in the form
\[D^{\overline{\mu}}\me a-\me T_{\overline{\mu}}\me a_{\overline{\mu}}=\me P_{\overline{\mu}}\me c_{\overline{\mu}},\]
 where $$\me a_{\overline{\mu} }=((\vec{a}_{\overline{\mu}}(0))^{\top},(\vec{a}_{\overline{\mu}}(1))^{\top},(\vec{a}_{\overline{\mu}}(2))^{\top}, \ldots)^{\top},$$
$$\me c_{\overline{\mu} }=((\vec{c}_{\overline{\mu}}(0))^{\top},(\vec{c}_{\overline{\mu}}(1))^{\top},(\vec{c}_{\overline{\mu}}(2))^{\top}, \ldots)^{\top},$$
$$\me a=((\vec{a}(0))^{\top},(\vec{a}(1))^{\top},(\vec{a}(2))^{\top}, \ldots)^{\top},$$ $\me P_{\overline{\mu}}$ and $\me T_{\overline{\mu}}$ are
linear operators in the space $\ell_2$ defined by matrices with the $T\times T$ matrix entries
 $(\me P_{\overline{\mu}})_{l,k}=P_{l,k}^{\overline{\mu}}$, $l,k\geq0$ and $(\me T_{\overline{\mu}})_{l, k} =T^{\overline{\mu}}_{l,k}$, $l,k\geq0$; the linear transformation $D^{\overline{\mu}}$ is defined in Lemma \ref{lema predst A}.

Consequently, the unknown coefficients $\vec{c}_{\overline{\mu}}(k)$, $k\geq0$, which determine
the spectral characteristic $\vec{h}_{\overline{\mu}}(\lambda)$ are calculated by the formula
\begin{equation}\label{meq-e}
\vec{c}_{\overline{\mu} }(k)=(\me P_{\overline{\mu}}^{-1}D^{\overline{\mu}}\me a-\me P_{\overline{\mu}}^{-1}\me T_{\overline{\mu}}\me a_{\overline{\mu}})_k,\quad k\geq 0,
 \end{equation}
where $(\me P_{\overline{\mu}}^{-1}D^{\overline{\mu}}\me a-\me
P_{\overline{\mu}}^{-1}\me T_{\overline{\mu}}\me a_{\overline{\mu}})_k$, $k\geq 0$, is the
$k$th $T$-dimension vector element of the vector  $\me P_{\overline{\mu}}^{-1}D^{\overline{\mu}}\me a-\me P_{\overline{\mu}}^{-1}\me T_{\overline{\mu}}\me a_{\overline{\mu}}$ and
the function $\vec{C}_{\overline{\mu}}(e^{i \lambda})$
is of the form
\[
\vec{C}_{\overline{\mu}}(e^{i \lambda})=\sum_{k=0}^{\infty}
(\me P_{\overline{\mu}}^{-1}D^{\overline{\mu}}\me a-\me P_{\overline{\mu}}^{-1}\me T_{\overline{\mu}}\me a_{\overline{\mu}})_k
e^{ik\lambda}.
\]

\begin{zau}
The problem of projection of the element ${H}\vec\xi$ of the Hilbert space $H$ on the closed convex set $H^{0-}(\xi^{(n)}_{\overline{\mu}}+\eta^{(n)}_{\overline{\mu}} )$ has a unique solution for each non-zero coefficients $\{\vec{a}(0),\vec{a}(1)),\vec{a}(2), \ldots\}$, satisfying conditions
(\ref{umovana a_e_d}) -- (\ref{umovana a_mu_e_d}). Therefore, equation (\ref{meq-e}) has a unique solution for each vector $D^{\overline{\mu}} \me a$, which implies an existence of the inverse operator $\me P^{-1}_{\overline{\mu}}$.
\end{zau}

The spectral characteristic $\vec{h}_{\overline{\mu}}(\lambda)$ of the optimal estimate $\widehat{H}\vec\xi$ of the functional $H\vec\xi$ can be calculated by the formula
\[
(\vec{h}_{\overline{\mu}}(\lambda))^{\top}=(\vec{B}_{\overline{\mu}}(e^{i\lambda}))^{\top}
\frac{\chi_{\overline{\mu}}^{(d)}(e^{-i\lambda})}{\beta^{(d)}(i\lambda)}
-\frac{\overline{\beta^{(d)}(i\lambda)}}
{\overline{\chi_{\overline{\mu}}^{(d)}(e^{-i\lambda})}}(\vec{A}_{\overline{\mu}}(e^{i\lambda}))^{\top}g(\lambda)
(f(\lambda)+|\beta^{(d)}(i\lambda)|^2g(\lambda))^{-1}
-
\]
\begin{equation}\label{spectr A_e_d}
-
\frac{\overline{\beta^{(d)}(i\lambda)}}
{\overline{\chi_{\overline{\mu}}^{(d)}(e^{-i\lambda})}}
\left(
\sum_{k=0}^{\infty}(\me P_{\overline{\mu}}^{-1}D^{\overline{\mu}}\me a-\me P_{\overline{\mu}}^{-1}\me T_{\overline{\mu}}\me a_{\overline{\mu}})_k e^{ik\lambda}
\right)^{\top}(f(\lambda)+|\beta^{(d)}(i\lambda)|^2g(\lambda))^{-1},
\end{equation}

The value of the mean square error of the estimate $\widehat{A}\vec\xi$ is calculated by the formula
\[
\Delta(f,g;\widehat{A}\vec\xi)=\Delta(f,g;\widehat{H}\vec\xi)= \mt E|H\vec\xi-\widehat{H}\vec\xi|^2=
\]
\[
=
\frac{1}{2\pi}\int_{-\pi}^{\pi}
\left[(\vec{A}_{\overline{\mu}}(e^{i\lambda}))^{\top}g(\lambda) +
\left(\sum_{k=0}^{\infty}(\me P_{\overline{\mu}}^{-1}D^{\overline{\mu}}\me a-\me P_{\overline{\mu}}^{-1}\me T_{\overline{\mu}}\me a_{\overline{\mu}})_k e^{ik\lambda}
\right)^{\top}
\right]
\times
\]
\[
\times
\frac{|\beta^{(d)}(i\lambda)|^2}{|\chi_{\overline{\mu}}^{(d)}(e^{-i\lambda})|^2}(f(\lambda)+|\beta^{(d)}(i\lambda)|^2g(\lambda))^{-1}\, f(\lambda)\, (f(\lambda)+|\beta^{(d)}(i\lambda)|^2g(\lambda))^{-1}
\times
\]
\[
\times
\left[g(\lambda)\overline{\vec{A}_{\overline{\mu}}(e^{i\lambda})} +
\overline{\sum_{k=0}^{\infty}(\me P_{\overline{\mu}}^{-1}D^{\overline{\mu}}\me a-\me P_{\overline{\mu}}^{-1}\me T_{\overline{\mu}}\me a_{\overline{\mu}})_k e^{ik\lambda}}
\right]
d\lambda+
\]
\[
+\frac{1}{2\pi}\int_{-\pi}^{\pi}
\left[{\overline{\chi_{\overline{\mu}}^{(d)}(e^{-i\lambda})}}(\vec{A}(e^{i\lambda}))^{\top}f(\lambda) -|\beta^{(d)}(i\lambda)|^2
\left(\sum_{k=0}^{\infty}(\me P_{\overline{\mu}}^{-1}D^{\overline{\mu}}\me a-\me P_{\overline{\mu}}^{-1}\me T_{\overline{\mu}}\me a_{\overline{\mu}})_k e^{ik\lambda}
\right)^{\top}
\right]
\times
\]
\[
\times
\frac{1}{|\chi_{\overline{\mu}}^{(d)}(e^{-i\lambda})|^2}(f(\lambda)+|\beta^{(d)}(i\lambda)|^2g(\lambda))^{-1}\,g(\lambda)\,(f(\lambda)+|\beta^{(d)}(i\lambda)|^2g(\lambda))^{-1}
\times
\]
\[
\times
\left[f(\lambda)\overline{\vec{A}_{\overline{\mu}}(e^{i\lambda})} -|\beta^{(d)}(i\lambda)|^2
\overline{\sum_{k=0}^{\infty}(\me P_{\overline{\mu}}^{-1}D^{\overline{\mu}}\me a-\me P_{\overline{\mu}}^{-1}\me T_{\overline{\mu}}\me a_{\overline{\mu}})_k e^{ik\lambda}}
\right]
d\lambda=
\]
\be\label{pohybkaA}
=\ld\langle D^{\overline{\mu}}\me a- \me T_{\overline{\mu}}\me
 a_{\overline{\mu}},\me P_{\overline{\mu}}^{-1}D^{\overline{\mu}}\me a-\me P_{\overline{\mu}}^{-1}\me T_{\overline{\mu}}\me a_{\overline{\mu}}\rd\rangle+\ld\langle\me Q\me a,\me
 a\rd\rangle,
 \ee

\noindent where $\me Q$ is a linear operator in the space $\ell_2$ defined by the matrix with the $T\times T$ matrix  elements $(\me Q)_{l,k}=Q_{l,k}$, $l,k\geq0$;  $\langle \vec x, \vec y\rangle=\sum_{k=0}^{\infty}(\vec x(k))^{\top}\overline{\vec y}(k)$ for vectors $\vec x=((\vec x(0))^{\top},(\vec x(1))^{\top},(\vec x(2))^{\top},\ldots)^{\top}$, $\vec y=((\vec y(0))^{\top},(\vec y(1))^{\top},(\vec y(2))^{\top},\ldots)^{\top}$.

\begin{thm}\label{thm1_e_st.n_d}
Let $\vec\xi(m)$, $m\in\mr Z$, be a stochastic sequence which defines
stationary $n$th increment sequence $\chi_{\overline{\mu},\overline{s}}^{(d)}(\vec{\xi}(m))$ with
absolutely continuous spectral function $F(\lambda)$ which has
spectral density $f(\lambda)$. Let $\vec\eta(m)$, $m\in\mr Z$, be an
uncorrelated with the sequence $\vec\xi(m)$ stationary stochastic
sequence with  absolutely continuous spectral function $G(\lambda)$
which has spectral density $g(\lambda)$. Let the minimality condition
(\ref{umova11_e_st.n_d}) be satisfied. Let coefficients
$\vec{a}(k)$, $k\geq0$,
satisfy  conditions (\ref{umovana a_e_d}) -- (\ref{umovana a_mu_e_d}). The optimal linear
estimate $\widehat{A}\vec\xi$ of the functional $A\vec\xi$ which depend on the unknown values of
elements $\vec\xi(m)$, $m\geq0$, based on observations of the sequence
$\vec\xi (m )+\vec\eta (m )$ at points $m=-1,-2,\ldots$ is calculated by
formula (\ref{otsinka A_e_d}). The spectral characteristic
$\vec{h}_{\overline{\mu}}(\lambda)$ of the optimal estimate $\widehat{A}\vec\xi$ is
calculated by formula (\ref{spectr A_e_d}). The value of the
mean-square error  $\Delta(f,g;\widehat{A}\vec\xi)$ is calculated by
formula (\ref{pohybkaA}).
\end{thm}

\begin{nas}
The spectral characteristic $\vec{h}_{\overline{\mu}}(\lambda)$ admits the representation $\vec{h}_{\overline{\mu}}(\lambda)=\vec{h}_{\overline{\mu}}^{1}(\lambda)-\vec{h}_{\overline{\mu}}^2(\lambda)$, where
\begin{multline}
\label{poh A1_e_st.n_d}
(\vec{h}_{\overline{\mu}}^1(\lambda))^{\top}=(\vec{B}_{\overline{\mu}}(e^{i\lambda}))^{\top}
\frac{\chi_{\overline{\mu}}^{(d)}(e^{-i\lambda})}{\beta^{(d)}(i\lambda)}-
\\
-
\frac{\overline{\beta^{(d)}(i\lambda)}}
{\overline{\chi_{\overline{\mu}}^{(d)}(e^{-i\lambda})}}
\left(
\sum_{k=0}^{\infty}(\me (P_{\overline{\mu}}^{-1}D^{\overline{\mu}}\me a)_k e^{ik\lambda}
\right)^{\top}(f(\lambda)+|\beta^{(d)}(i\lambda)|^2g(\lambda))^{-1},
\end{multline}
\begin{multline}\label{poh A2_e_st.n_d}
(\vec{h}_{\overline{\mu}}^2(\lambda))^{\top}
=
\frac{\overline{\beta^{(d)}(i\lambda)}}
{\overline{\chi_{\overline{\mu}}^{(d)}(e^{-i\lambda})}}
(\vec{A}_{\overline{\mu}}(e^{i\lambda}))^{\top}
 g(\lambda)(f(\lambda)+|\beta^{(d)}(i\lambda)|^2g(\lambda))^{-1}-
 \\
 -
\frac{\overline{\beta^{(d)}(i\lambda)}
}
{\overline{\chi_{\overline{\mu}}^{(d)}(e^{-i\lambda})}}\left(
\sum_{k=0}^{\infty}(\me P_{\overline{\mu}}^{-1}\me T_{\overline{\mu}}\me a_{\overline{\mu}})_k e^{ik\lambda}
\right)^{\top}(f(\lambda)+|\beta^{(d)}(i\lambda)|^2g(\lambda))^{-1}.
\end{multline}
Here $\vec{h}_{\overline{\mu}}^1(\lambda)$ and $\vec{h}_{\overline{\mu}}^2(\lambda)$ are spectral characteristics of the optimal estimates $\widehat{B}\vec\zeta$ and $\widehat{A}\vec\eta$ of the functionals $B\vec\zeta$ and $A\vec\eta$ respectively based on observations $\vec\xi(m)+\vec\eta(m)$ at points $m=-1,-2,\ldots$.
\end{nas}

From Theorem $\ref{thm1_e_st.n_d}$  obtain the optimal estimate
$\widehat{A}_N\vec\xi$ of the functional $A_N\vec\xi$ of the unknown values of elements
$\vec\xi (m )$, $m=0,1,2,\ldots,N$, based on observations of the sequence
$\vec\xi (m )+\vec\eta (m )$ at  points $m=-1,-2,\ldots$. Let $\vec{a}(k)=0$,
$k>N$. Then the spectral characteristic $\vec{h}_{\overline{\mu}, N}(\lambda)$ of the
linear estimate
\be \label{otsinka A_N_e_st.n_d}
 \widehat{A}_N\vec\xi=\ip
(\vec{h}_{\overline{\mu},N}(\lambda))^{\top}d\vec{Z}_{\xi^{(n)}+\eta^{(n)}}(\lambda)-\sum_{k=-n(\gamma)}^{-1}(\vec v_{\overline{\mu},N}(k))^{\top}(\vec\xi(k)+\vec\eta(k)),
\ee
 is calculated by the formula
\begin{multline}\label{spectr A_e_dN}
(\vec{h}_{\overline{\mu},N}(\lambda))^{\top}
=(\vec{B}_{\overline{\mu},N}(e^{i\lambda}))^{\top}
\frac{\chi_{\overline{\mu}}^{(d)}(e^{-i\lambda})}{\beta^{(d)}(i\lambda)}
-
\frac{\overline{\beta^{(d)}(i\lambda)}}
{\overline{\chi_{\overline{\mu}}^{(d)}(e^{-i\lambda})}}(\vec{A}_{\overline{\mu,N}}(e^{i\lambda}))^{\top}
g(\lambda)
(f(\lambda)+|\beta^{(d)}(i\lambda)|^2g(\lambda))^{-1}
-
\\
-
\frac{\overline{\beta^{(d)}(i\lambda)}}
{\overline{\chi_{\overline{\mu}}^{(d)}(e^{-i\lambda})}}
\left(
\sum_{k=0}^{\infty}(\me P_{\overline{\mu}}^{-1}D_N^{\overline{\mu}}\me a_N-\me P_{\overline{\mu}}^{-1}\me T_{\overline{\mu},N}\me a_{\overline{\mu},N})_k e^{ik\lambda}
\right)^{\top}(f(\lambda)+|\beta^{(d)}(i\lambda)|^2g(\lambda))^{-1},
\end{multline}
where
\[
\vec{B}_{\overline{\mu},N}(e^{i\lambda})=\sum_{k=0}^{N}(D^{\overline{\mu}}_N\me a_N)_ke^{i\lambda k},
\quad \vec{A}_N(e^{i\lambda })=\sum_{k=0}^{N}\vec{a}(k)e^{i\lambda k},\quad
\vec{A}_{\overline{\mu},N}(e^{i\lambda })=\sum_{k=0}^{N}\vec{a}_{\overline{\mu},N}(k)e^{i\lambda k},
\]
 \[\me a_N=((\vec{a}(0))^{\top},(\vec{a}(1))^{\top},\ldots,(\vec{a}(N))^{\top},0, \ldots)^{\top},\]
 \[\me a_{\overline{\mu},N}=((\vec{a}_{\overline{\mu},N}(0))^{\top},(\vec{a}_{\overline{\mu},N}(1))^{\top},\ldots, (\vec{a}_{\overline{\mu},N}(N+n(\gamma)))^{\top},0,\ldots )^{\top},\]
\be\label{coeff a_N_mu}
\vec{a}_{\overline{\mu},N}(m)=\sum_{l=\max\ld\{m-n(\gamma),0\rd\}}^{\min\ld\{m,N\rd\}}e_{\gamma}(m-l)\vec{a}(l),\quad 0\leq m\leq N+n(\gamma),
\ee
operator $\me T_{\overline{\mu},N}$ is a linear operator in the space $\ell_2$ defined by the matrix with the $T\times T$ matrix  entries
  $(\me T_{\overline{\mu},N})_{l, m} =T^{\overline{\mu}}_{l,m}$, $l\geq0$, $0\leq m\leq N+n(\gamma)$, and $(\me T_{\overline{\mu},N})_{l, m} =0$, $l\geq0$, $m>N+n(\gamma)$.

The value of the mean square error of the estimate $\widehat{A}_N\vec\xi$ is calculated by the formula
\begin{multline*}
\Delta(f,g;\widehat{A}_N\vec\xi)=\Delta(f,g;\widehat{H}_N\vec\xi)= \mt E|H_N\vec\xi-\widehat{H}_N\vec\xi|^2=
\\
=
\frac{1}{2\pi}\int_{-\pi}^{\pi}
\frac{|\beta^{(d)}(i\lambda)|^2}{|\chi_{\overline{\mu}}^{(d)}(e^{-i\lambda})|^2}
\bigg[(\vec{A}_{\overline{\mu},N}(e^{i\lambda}))^{\top}g(\lambda) +
\\
+\left(\sum_{k=0}^{\infty}(\me P_{\overline{\mu}}^{-1}D_N^{\overline{\mu}}\me a_N-\me P_{\overline{\mu}}^{-1}\me T_{\overline{\mu},N}\me a_{\overline{\mu},N})_k e^{ik\lambda}
\right)^{\top}
\bigg]
\times
\\
\times
(f(\lambda)+|\beta^{(d)}(i\lambda)|^2g(\lambda))^{-1}\, f(\lambda)\, (f(\lambda)+|\beta^{(d)}(i\lambda)|^2g(\lambda))^{-1}
\times
\\
\times
\left[g(\lambda)\overline{\vec{A}_{\overline{\mu},N}(e^{i\lambda})} +
\overline{\sum_{k=0}^{\infty}(\me P_{\overline{\mu}}^{-1}D_N^{\overline{\mu}}\me a_N-\me P_{\overline{\mu}}^{-1}\me T_{\overline{\mu},N}\me a_{\overline{\mu},N})_k e^{ik\lambda}}
\right]
d\lambda
\end{multline*}
\begin{multline}
\label{pohybka33}
+\frac{1}{2\pi}\int_{-\pi}^{\pi}
\frac{1}{|\chi_{\overline{\mu}}^{(d)}(e^{-i\lambda})|^2}
\bigg[(\vec{A}_{\overline{\mu},N}(e^{i\lambda}))^{\top}f(\lambda) -
\\
-|\beta^{(d)}(i\lambda)|^2
\left(\sum_{k=0}^{\infty}(\me P_{\overline{\mu}}^{-1}D_N^{\overline{\mu}}\me a_N-\me P_{\overline{\mu}}^{-1}\me T_{\overline{\mu},N}\me a_{\overline{\mu},N})_k e^{ik\lambda}
\right)^{\top}
\bigg]
\times
\\
\times
(f(\lambda)+|\beta^{(d)}(i\lambda)|^2g(\lambda))^{-1}\,g(\lambda)\,(f(\lambda)+|\beta^{(d)}(i\lambda)|^2g(\lambda))^{-1}
\times
\\
\times
\bigg[f(\lambda)\overline{\vec{A}_N(e^{i\lambda})} -
\\
-|\beta^{(d)}(i\lambda)|^2
\overline{\sum_{k=0}^{\infty}(\me P_{\overline{\mu}}^{-1}D_N^{\overline{\mu}}\me a_N-\me P_{\overline{\mu}}^{-1}\me T_{\overline{\mu},N}\me a_{\overline{\mu},N})_k e^{ik\lambda}}
\bigg]
d\lambda=
\\
=\langle D_N^{\overline{\mu}}\me a_N- \me T_{\overline{\mu},N}\me
 a_{\overline{\mu},N},\me P_{\overline{\mu}}^{-1}D_N^{\overline{\mu}}\me a_N-\me P_{\overline{\mu}}^{-1}\me T_{\overline{\mu},N}\me a_{\overline{\mu},N}\rangle+\langle\me Q_N\me a_N,\me  a_N\rangle,
\end{multline}
 where
 $\me Q_{N}$ is a linear operator in the space $\ell_2$ defined by the matrix with the $T\times T$ matrix  elements $(\me Q_{N})_{l,k}=Q_{l,k}$, $0\leq l,k\leq N$, and $(\me Q_{N})_{l,k}=0$ otherwise.

The following theorem holds true.

\begin{thm}\label{thm2AN}
Let $\vec\xi(m)$, $m\in\mr Z\}$, be a stochastic sequence which defines
stationary $n$th increment sequence $\chi_{\overline{\mu},\overline{s}}^{(d)}(\vec{\xi}(m))$ with an
absolutely continuous spectral function $F(\lambda)$ which has
spectral density $f(\lambda)$. Let $\vec\eta(m)$, $m\in\mr Z$, be an
uncorrelated with the sequence $\vec\xi(m)$ stationary stochastic
sequence with an absolutely continuous spectral function
$G(\lambda)$ which has spectral density $g(\lambda)$. Let the minimality condition
(\ref{umova11_e_st.n_d}) be satisfied. The
optimal linear estimate $\widehat{A}_N\vec\xi$ of the functional
$A_N\vec\xi$ which depend on the unknown values of elements $\vec\xi(k)$, $k=0,1,2,\ldots,N$,   from
observations of the sequence $\vec\xi (m )+\vec\eta (m )$ at points
$m=-1,-2,\ldots$ is calculated by  formula (\ref{otsinka A_N_e_st.n_d}).
The spectral characteristic $\vec{h}_{\overline{\mu},N}(\lambda)$ of the optimal
estimate $\widehat{A}_N\vec\xi$ is calculated by  formula (\ref{spectr A_e_dN}). The value of the mean-square error
$\Delta(f,g;\widehat{A}_N\vec\xi)$ is calculated by  formula (\ref{pohybka33}).
\end{thm}

As a corollary from the proposed theorem, one can obtain the mean square optimal estimate of the unobserved value
$A_{N,p}\vec\xi=\xi_p(N)=(\vec\xi(N))^{\top}\boldsymbol{\delta}_p$, $p=1,2,\dots,T$, $N\geq0$
of the stochastic sequence with $n$th stationary increments based on observations of the sequence $\vec\xi (m )+\vec\eta (m )$ at points $m=-1,-2,\ldots$

\begin{nas}\label{nas xi_e_d}
The optimal linear estimate $\widehat{\xi}_p(N)$ of the unobserved value
$\xi_p(N)$, $p=1,2,\dots,T$, $N\geq0$, of the stochastic sequence with  stationary GM increments from observations of the sequence $\vec\xi (m )+\vec\eta (m )$ at points $m=-1,-2,\ldots$ is calculated by formula
\begin{multline} \label{est_xi_N}
 \widehat{\xi}_p(N)=\ip
(\vec h_{\overline{\mu},N,p}(\lambda))^{\top}d\vec Z_{\xi^{(n)}+\eta^{(n)}}(\lambda)
- \sum_{k=-n(\gamma)}^{-1}(\vec v_{\overline{\mu},N,p}(k))^{\top}(\vec\xi(k)+\vec\eta(k)).
\end{multline}

The spectral characteristic $\vec h_{\overline{\mu},N,p}(\lambda)$ of the estimate is calculated by the formula
\begin{multline}\label{sph_est_xi_N}
(\vec{h}_{\overline{\mu},N,p}(\lambda))^{\top}=\frac{\chi_{\overline{\mu}}^{(d)}(e^{-i\lambda})}{\beta^{(d)}(i\lambda)}\left(\boldsymbol{\delta}_p\sum_{k=0}^Nd_{\overline{\mu}}(N-k)e^{i\lambda k}\right)^{\top}
-
\\
-
\left(e^{i\lambda N}\boldsymbol{\delta}_p \right)^{\top}
g(\lambda)\overline{\beta^{(d)}(i\lambda)}
(f(\lambda)+|\beta^{(d)}(i\lambda)|^2g(\lambda))^{-1}
\\
-
\frac{\overline{\beta^{(d)}(i\lambda)}
\left(\sum_{k=0}^{\infty} (\me P_{\overline{\mu}}^{-1}\me d_{\overline{\mu},N,p}-\me P_{\overline{\mu}}^{-1}\me T_{\overline{\mu},N}\vec {\me a}_{\overline{\mu},N,p})_k e^{i\lambda k}\right)^{\top}}
{\overline{\chi_{\overline{\mu}}^{(d)}(e^{-i\lambda})}}
(f(\lambda)+|\beta^{(d)}(i\lambda)|^2g(\lambda))^{-1}.
\end{multline}
where
\begin{eqnarray*}
\me d_{\overline{\mu},N,p}&=&((d_{\overline{\mu}}(N)\boldsymbol{\delta}_p)^{\top},(d_{\overline{\mu}}(N-1)\boldsymbol{\delta}_p)^{\top},(d_{\overline{\mu}}(N-2)\boldsymbol{\delta}_p)^{\top},\ldots,(d_{\overline{\mu}}(0)\boldsymbol{\delta}_p)^{\top},0,\ldots)^{\top},
 \\
     \vec {\me a}_{\overline{\mu},N,p}&=&(0,\ldots,0,(\vec  a_{\overline{\mu},N,p}(N))^{\top},(\vec  a_{\overline{\mu},N,p}(N+1))^{\top},\ldots, (\vec  a_{\overline{\mu},N,p}(N+n(\gamma)))^{\top},0,\ldots )^{\top},
\\
   \vec   a_{\overline{\mu},N,p}(m)&=&e_{\gamma}(m-N)\vec \delta_p,\quad N\leq m\leq N+n(\gamma).
\end{eqnarray*}

The value of the mean square error of the optimal estimate is calculated by the formula
\[
\Delta(f,g;\widehat{\xi}_p(N))= \mt E|\xi_p(N)-\widehat{\xi}_p(N)|^2=
\]
\begin{multline*}
=
\frac{1}{2\pi}\int_{-\pi}^{\pi}
\frac{|\beta^{(d)}(i\lambda)|^2}{|\chi_{\overline{\mu}}^{(d)}(e^{-i\lambda})|^2}
\bigg[{\overline{\chi_{\overline{\mu}}^{(d)}(e^{-i\lambda})}}\left(e^{i\lambda N}\boldsymbol{\delta}_p \right)^{\top}g(\lambda) +
\\
+\left(\sum_{k=0}^{\infty} (\me P_{\overline{\mu}}^{-1}\me d_{\overline{\mu},N,p}-\me P_{\overline{\mu}}^{-1}\me T_{\overline{\mu},N}\vec {\me a}_{\overline{\mu},N,p})_k e^{i\lambda k}\right)^{\top}
\bigg]
\times
\\
\times
(f(\lambda)+|\beta^{(d)}(i\lambda)|^2g(\lambda))^{-1}\, f(\lambda)\, (f(\lambda)+|\beta^{(d)}(i\lambda)|^2g(\lambda))^{-1}
\times
\\
\times\bigg[g(\lambda){\chi_{\overline{\mu}}^{(d)}(e^{-i\lambda})}
{\left(e^{-i\lambda N}\boldsymbol{\delta}_p \right)} +
\overline{\sum_{k=0}^{\infty} (\me P_{\overline{\mu}}^{-1}\me d_{\overline{\mu},N,p}-\me P_{\overline{\mu}}^{-1}\me T_{\overline{\mu},N}\vec {\me a}_{\overline{\mu},N,p})_k e^{i\lambda k}}
\bigg]
d\lambda+
\end{multline*}
\begin{multline*}
+\frac{1}{2\pi}\int_{-\pi}^{\pi}
\frac{1}{|\chi_{\overline{\mu}}^{(d)}(e^{-i\lambda})|^2}
\bigg[{\overline{\chi_{\overline{\mu}}^{(d)}(e^{-i\lambda})}}\left(e^{i\lambda N}\boldsymbol{\delta}_p \right)^{\top}f(\lambda) -
\\
-
|\beta^{(d)}(i\lambda)|^2
\left(
\sum_{k=0}^{\infty} (\me P_{\overline{\mu}}^{-1}\me d_{\overline{\mu},N,p}-\me P_{\overline{\mu}}^{-1}\me T_{\overline{\mu},N}\vec {\me a}_{\overline{\mu},N,p})_k e^{i\lambda k}
\right)^{\top}
\bigg]
\times
\\
\times(f(\lambda)+|\beta^{(d)}(i\lambda)|^2g(\lambda))^{-1}\,g(\lambda)\,(f(\lambda)+|\beta^{(d)}(i\lambda)|^2g(\lambda))^{-1}
\times
\\
\times \bigg[f(\lambda){\chi_{\overline{\mu}}^{(d)}(e^{-i\lambda})}{\left(e^{-i\lambda N}\boldsymbol{\delta}_p \right)} -
|\beta^{(d)}(i\lambda)|^2
\overline{\sum_{k=0}^{\infty} (\me P_{\overline{\mu}}^{-1}\me d_{\overline{\mu},N,p}-\me P_{\overline{\mu}}^{-1}\me T_{\overline{\mu},N}\vec {\me a}_{\overline{\mu},N,p})_k e^{i\lambda k}}
\bigg]
d\lambda=
\end{multline*}
    \be\label{poh xi_p_e_st.n_d}
     = \langle \me d_{\overline{\mu},N,p}-\me T_{\overline{\mu},N}\vec {\me a}_{\overline{\mu},N,p},\me P_{\overline{\mu}}^{-1}\me d_{\overline{\mu},N,p}-\me
P_{\overline{\mu}}^{-1}\me T_{\overline{\mu},N}\vec {\me a}_{\overline{\mu},N,p}\rangle+ \langle \me Q_{0,0}\boldsymbol{\delta}_p,\boldsymbol{\delta}_p \rangle.
\ee
\end{nas}

\begin{zau}
The filtering problem in the presence of fractional integration can be solved using Theorem \ref{thm1_e_st.n_d}, Theorem \ref{thm2AN} and Corollary \ref{nas xi_e_d} under conditions of Theorem \ref{thm_frac} on the increment  orders $d_i$.
\end{zau}

\subsection{Forecasting based on factorizations of the spectral densities}

In  Theorem \ref{thm1_e_st.n_d}, Theorem \ref{thm2AN} and Corollary \ref{nas xi_e_d}, formulas for finding forecasts  of the linear functionals $A\xi$, $A_N\xi$
and the value $\xi(p)$, $p\geq0$, are derived  using the Fourier coefficients of
the functions
\[
\frac{|\beta^{(d)}(i\lambda)|^2}{|\chi_{\overline{\mu}}^{(d)}(e^{-i\lambda})|^2}
\ld[g(\lambda)(f(\lambda)+|\beta^{(d)}(i\lambda)|^2g(\lambda))^{-1}\rd]^{\top},\]
\[
\dfrac{|\beta^{(d)}(i\lambda)|^2}{|\chi_{\overline{\mu}}^{(d)}(e^{-i\lambda})|^2}
\ld[(f(\lambda)+|\beta^{(d)}(i\lambda)|^2g(\lambda))^{-1}\rd]^{\top}.\]
Assume that the following canonical factorizations take place
\begin{multline} \label{fakt1}
 \dfrac{|\chi_{\overline{\mu}}^{(d)}(e^{-i\lambda})|^2}{|\beta^{(d)}(i\lambda)|^2}
 (f(\lambda)+|\beta^{(d)}(i\lambda)|^2g(\lambda))
 =\Theta_{\overline{\mu}}(e^{-i\lambda})\Theta_{\overline{\mu}}^*(e^{-i\lambda}),\quad \Theta_{\overline{\mu}}(e^{-i\lambda})=\sum_{k=0}^{\infty}\theta_{\overline{\mu}}(k)e^{-i\lambda k},
\end{multline}
\be \label{fakt3}
 g(\lambda)=\sum_{k=-\infty}^{\infty}g(k)e^{i\lambda k}=\Phi(e^{-i\lambda})\Phi^*(e^{-i\lambda}), \quad
 \Phi(e^{-i\lambda})=\sum_{k=0}^{\infty}\phi(k)e^{-i\lambda k}.
 \ee
Define the matrix-valued function $\Psi_{\overline{\mu}}(e^{-i\lambda})= \{\Psi_{ij}(e^{-i\lambda})\}_{i=\overline{1,q}}^{j=\overline{1,T}}$ by the equation
\[
\Psi_{\overline{\mu}}(e^{-i\lambda})\Theta_{\overline{\mu}}(e^{-i\lambda})=E_q,
\]
where $E_q$ is an identity $q\times q$ matrix.
One can check that the following factorization takes place
\begin{multline} \label{fakt2}
\dfrac{|\beta^{(d)}(i\lambda)|^2}{|\chi_{\overline{\mu}}^{(d)}(e^{-i\lambda})|^2}
 (f(\lambda)+|\beta^{(d)}(i\lambda)|^2g(\lambda))^{-1} =
 \Psi_{\overline{\mu}}^*(e^{-i\lambda})\Psi_{\overline{\mu}}(e^{-i\lambda}), \quad
 \Psi_{\overline{\mu}}(e^{-i\lambda})=\sum_{k=0}^{\infty}\psi_{\overline{\mu}}(k)e^{-i\lambda k},
\end{multline}

\begin{zau}\label{remark_density_adjoint}
Any spectral density matrix $f(\lambda)$ is self-adjoint: $f(\lambda)=f^*(\lambda)$. Thus, $(f(\lambda))^{\top}=\overline{f(\lambda)}$. One can check that an inverse spectral density $f^{-1}(\lambda)$ is also self-adjoint $f^{-1}(\lambda)=(f^{-1}(\lambda))^*$ and $(f^{-1}(\lambda))^{\top}=\overline{f^{-1}(\lambda)}$.
\end{zau}

The following Lemmas provide factorizations of the operators $\me P_{\overline{\mu}}$ and $\me T_{\overline{\mu}}$, which contain  coefficients of factorizations (\ref{fakt1}) -- (\ref{fakt2}).

\begin{lema}\label{lema_fact_2}
Let factorization (\ref{fakt1})  takes place and let $q\times T$ matrix function $\Psi_{\overline{\mu}}(e^{-i\lambda})$ satisfy equation $\Psi_{\overline{\mu}}(e^{-i\lambda})\Theta_{\overline{\mu}}(e^{-i\lambda})=E_q$.
Define the linear operators
 $ \Psi_{\overline{\mu}}$ and  $ \Theta_{\overline{\mu}}$ in the space $\ell_2$ by the matrices with the matrix entries
 $( \Psi_{\overline{\mu}})_{k,j}=\psi_{\overline{\mu}}(k-j)$, $( \Theta_{\overline{\mu}})_{k,j}=\theta_{\overline{\mu}}(k-j)$ for $0\leq j\leq k$, $(
\Psi_{\overline{\mu}})_{k,j}=0$, $(
\Theta_{\overline{\mu}})_{k,j}=0$ for $0\leq k<j$.
Then:
\\
a) the linear operator $\me P_{\overline{\mu}}$  admits the factorization \[\me
P_{\overline{\mu}}=(\Psi_{\overline{\mu}})^{\top} \overline{\Psi}_{\overline{\mu}};\]
\\
b) the inverse operator $(\me
P_{\overline{\mu}})^{-1}$ admits the factorization
\[
 (\me
P_{\overline{\mu}})^{-1}= \overline{\Theta}_{\overline{\mu}}(\Theta_{\overline{\mu}})^{\top}.\]
\end{lema}

$\mathrm{\mbox{Proof.}}$ Making use of  factorization (\ref{fakt2}), obtain the relation
\begin{eqnarray*}
 && \dfrac{|\beta^{(d)}(i\lambda)|^2}{|\chi_{\overline{\mu}}^{(d)}(e^{-i\lambda})|^2}
 \ld[(f(\lambda)+|\beta^{(d)}(i\lambda)|^2g(\lambda))^{-1}\rd]^{\top}
 \\
 &=&
 \sum_{m=-\infty}^{\infty}P_{\overline{\mu}}(m)e^{i\lambda m}
 =
 (\Psi_{\overline{\mu}}(e^{-i\lambda}))^{\top}\overline{\Psi_{\overline{\mu}}(e^{-i\lambda})}
\\
 &=&\sum_{m=-\infty}^{-1}\sum_{k=-m}^{\infty}\psi^{\top}_{\overline{\mu}}(k)
 \overline{\psi}_{\overline{\mu}}(k+m)e^{i\lambda m}
\sum_{m=0}^{\infty}
 \sum_{k=0}^{\infty}\psi^{\top}_{\overline{\mu}}(k)
 \overline{\psi}_{\overline{\mu}}(k+m)e^{i\lambda m}.\end{eqnarray*}
Thus, $P_{\overline{\mu}}(m)=\sum_{k=0}^{\infty}\psi^{\top}_{\overline{\mu}}(k)
 \overline{\psi}_{\overline{\mu}}(k+m)$, $m\geq0$, and $P_{\overline{\mu}}(-m)=(P_{\overline{\mu}}(m))^*$, $m\geq0$. In the case $i\geq j$, we have \[
 P^{\mu}_{i,j}=P_{\overline{\mu}}(i-j)
 =\sum_{l=i}^{\infty}\psi^{\top}_{\overline{\mu}}(l-i)\overline{\psi}_{\overline{\mu}}(l-j)
 =((\Psi_{\overline{\mu}})^{\top} \overline{\Psi}_{\overline{\mu}})_{i,j},\]
 and, in the case $i<j$, we have
 \[
 P^{\mu}_{i,j}=P_{\overline{\mu}}(i-j)=\overline{P_{\overline{\mu}}(j-i)}
 =\sum_{l=j}^{\infty}\psi^{\top}_{\overline{\mu}}(l-i)\overline{\psi}_{\overline{\mu}}(l-j)
 =((\Psi_{\overline{\mu}})^{\top} \overline{\Psi}_{\overline{\mu}})_{i,j},\]
which proves statement a).

 From factorizations (\ref{fakt1}) and (\ref{fakt2}),  obtain
\begin{eqnarray} \label{spivvidn1}
E_q=\Psi_{\overline{\mu}}(e^{-i\lambda})\Theta_{\overline{\mu}}(e^{-i\lambda})
 =\sum_{j=0}^{\infty}\ld(\sum_{k=0}^{j}\psi_{\overline{\mu}}(k)\theta_{\overline{\mu}}(j-k)\rd)e^{-i\lambda j}.\end{eqnarray}
Thus, we  derive the relations
\begin{eqnarray*}
 \mt{diag_q}(\delta_{i,j})&=&\sum_{k=0}^{i-j}\psi_{\overline{\mu}}(k)\theta_{\overline{\mu}}(i-j-k)
 =\sum_{p=i}^j\psi_{\overline{\mu}}(i-p)\theta_{\overline{\mu}}(p-j)
=(\Psi_{\overline{\mu}}\Theta_{\overline{\mu}})_{i,j}, \quad \Box\end{eqnarray*}
which imply  $\Psi_{\overline{\mu}}\Theta_{\overline{\mu}}=E_q$. Using $\me
P_{\overline{\mu}}=(\Psi_{\overline{\mu}})^{\top} \overline{\Psi}_{\overline{\mu}}$ we get $\me
(\Theta_{\overline{\mu}})^{\top} \me P_{\overline{\mu}}=\overline{\Psi}_{\overline{\mu}}$ and $\me
P_{\overline{\mu}}\overline{\Theta}_{\overline{\mu}}=(\Psi_{\overline{\mu}})^{\top}$, or $
(\Theta_{\overline{\mu}})^{\top} =\overline{\Psi}_{\overline{\mu}}(\me P_{\overline{\mu}})^{-1}$ and $\overline{\Theta}_{\overline{\mu}}=(\me P_{\overline{\mu}})^{-1}(\Psi_{\overline{\mu}})^{\top}$. The last two relations imply
\[
\overline{\Theta}_{\overline{\mu}}(\Theta_{\overline{\mu}})^{\top} =(\me P_{\overline{\mu}})^{-1}(\Psi_{\overline{\mu}})^{\top}\overline{\Psi}_{\overline{\mu}}(\me P_{\overline{\mu}})^{-1}=(\me P_{\overline{\mu}})^{-1}\me P_{\overline{\mu}}(\me P_{\overline{\mu}})^{-1}=(\me P_{\overline{\mu}})^{-1},\]
which proves statement b).  $\Box$

\begin{lema}\label{lema_fact_1}
 Let  factorizations (\ref{fakt1}) and (\ref{fakt3}) take place. Then the operator $\me T_{\overline{\mu}}$ admits the representation
\[
\me T_{\overline{\mu}}=(\Psi_{\overline{\mu}})^{\top} \me Z_{\overline{\mu}},
\]
 where $\me Z_{\overline{\mu}}$ is a linear operator in the space $\ell_2$ defined by a matrix with the entries
 \[
 (\me Z_{\overline{\mu}})_{k,j}=\sum_{l=j}^{\infty}\overline{\psi}_{\overline{\mu}}(l-j)\overline{g}(l-k),\quad g(k)=\sum_{m=\max\{0,-k\}}^{\infty}\phi(m)\phi^*(k+m),\quad k,j\geq0.
 \]
\end{lema}

$\mathrm{\mbox{Proof.}}$ Factorizations (\ref{fakt3}), (\ref{fakt2}) and Remark \ref{remark_density_adjoint} imply
\[
\frac{|\beta^{(d)}(i\lambda)|^2}{|\chi_{\overline{\mu}}^{(d)}(e^{-i\lambda})|^2}
\ld[g(\lambda)(f(\lambda)+|\beta^{(d)}(i\lambda)|^2g(\lambda))^{-1}\rd]^{\top}
=(\Psi_{\overline{\mu}}(e^{-i\lambda}))^{\top}\overline{\Psi}_{\overline{\mu}}(e^{-i\lambda})\overline{g}(\lambda)
\]
\[
=\sum_{l=0}^{\infty}\psi_{\overline{\mu}}^{\top}(l)e^{-i\lambda l} \sum_{j=0}^{\infty}Z_{\overline{\mu}}(j)e^{i\lambda j}
=\sum_{k\in\mr Z} \sum_{j=0}^{\infty}\psi_{\overline{\mu}}^{\top}(l)Z_{\overline{\mu}}(l+k)e^{i\lambda k}.
\]
Then
\[
(\me T_{\overline{\mu}})_{k,j}=\me T_{\overline{\mu}}(k-j)=\sum_{m=k}^{\infty}\Psi_{\overline{\mu}}^{\top}(m-k)Z_{\overline{\mu}}(m-j)
=(\Psi_{\overline{\mu}}^{\top}\me Z_{\overline{\mu}})_{k,j}
\]
The representation for the entries $(\me Z_{\overline{\mu}})_{k,j}=Z_{\overline{\mu}}(k-j)$ follows from
\[
\sum_{j=0}^{\infty}Z_{\overline{\mu}}(j)e^{i\lambda j}
=\overline{\Psi}_{\overline{\mu}}(e^{-i\lambda})\overline{g}(\lambda)
=\sum_{k\in\mr Z} \sum_{l=0}^{\infty}\overline{\psi}_{\overline{\mu}}(l)\overline{g}_{\overline{\mu}}(l-k)e^{i\lambda k}.\quad \Box
\]

\begin{zau}\label{remark_operator_PTa}
Lemma $\ref{lema_fact_2}$ and Lemma $\ref{lema_fact_1}$ imply the factorization
\[
(\me P_{\overline{\mu}})^{-1}\me T_{\overline{\mu}} \me a_{\overline{\mu}}= \overline{\Theta}_{\overline{\mu}}(\Theta_{\overline{\mu}})^{\top}(\Psi_{\overline{\mu}})^{\top} \me Z_{\overline{\mu}}\me a_{\overline{\mu}}
=\overline{\Theta}_{\overline{\mu}} \me Z_{\overline{\mu}}\me a_{\overline{\mu}}=\overline{\Theta}_{\overline{\mu}} \me e_{\overline{\mu}},
\]
where $\me e_{\overline{\mu}}:=\me Z_{\overline{\mu}}\me a_{\overline{\mu}}$.
\end{zau}

Assuming that factorizations (\ref{fakt1}), (\ref{fakt3}), (\ref{fakt2}) take place and making use of   Remark \ref{remark_operator_PTa},   spectral characteristic  (\ref{spectr A_e_d}) and mean-squarer error  (\ref{pohybkaA}) can be presented in terms of the coefficients of the mentioned factorizations.
Make the following transformations:
\begin{eqnarray}
&&\notag\frac{|\beta^{(d)}(i\lambda)|^2}
{|\chi_{\overline{\mu}}^{(d)}(e^{-i\lambda})|^2}
\ld[(f(\lambda)+|\beta^{(d)}(i\lambda)|^2g(\lambda))^{-1}\rd]^{\top}\left(
\sum_{k=0}^{\infty}(\me P_{\overline{\mu}}^{-1}\me T_{\overline{\mu}}\me a_{\overline{\mu}})_k e^{ik\lambda}
\right)
\\\notag
&=&\ld(\sum_{k=0}^{\infty}\psi_{\overline{\mu}}^{\top}(k)e^{-i\lambda k}\rd)\sum_{j=0}^{\infty}\sum_{k=0}^{\infty}\overline{\psi}_{\overline{\mu}}(j)(\overline{\Theta}_{\overline{\mu}}\me e_{\overline{\mu}})_ke^{i\lambda(k+j)}
\\\notag&=&\ld(\sum_{k=0}^{\infty}\psi^{\top}_{\overline{\mu}}(k)e^{-i\lambda k}\rd)
\sum_{m=0}^{\infty}\sum_{p=0}^{m}\sum_{k=p}^m\overline{\psi}_{\overline{\mu}}(m-k)\overline{\theta}_{\overline{\mu}}(k-p)e_{\overline{\mu}}(p)e^{i\lambda m}
\\\notag&=&\ld(\sum_{k=0}^{\infty}\psi^{\top}_{\overline{\mu}}(k)e^{-i\lambda k}\rd)\sum_{m=0}^{\infty}\sum_{p=0}^{\infty}\mt{diag}(\delta_{m,p})e_{\overline{\mu}}(p)e^{i\lambda m}
\\&=&\ld(\sum_{k=0}^{\infty}\psi^{\top}_{\overline{\mu}}(k)e^{-i\lambda k}\rd)\sum_{m=0}^{\infty}e_{\overline{\mu}}(m)e^{i\lambda m},\label{simple_sp_char_part1}
\end{eqnarray}
where $ e_{\overline{\mu}}(m)=(\me Z_{\overline{\mu}}\me a_{\overline{\mu}})_m$, $m\geq0$, is the $m$-th vector entry of the vector
$\me e_{\overline{\mu}}=\me Z_{\overline{\mu}}\me a_{\overline{\mu}}$.

Using factorizations (\ref{fakt3}) and (\ref{fakt2}) conclude the following transformations:
\begin{align}
\notag & \frac{|\beta^{(d)}(i\lambda)|^2}
{\chi_{\overline{\mu}}^{(d)}(e^{-i\lambda})|^2}
\ld[(f(\lambda)+|\beta^{(d)}(i\lambda)|^2g(\lambda))^{-1}\rd]^{\top}(g(\lambda))^{\top}
\vec{A}_{\overline{\mu}}(e^{i\lambda})
 \\
 \notag&=\ld(\sum_{k=0}^{\infty}\psi^{\top}_{\overline{\mu}}(k)e^{-i\lambda k}\rd)\ld(\sum_{k=-\infty}^{\infty}Z_{\overline{\mu}}(k)e^{i\lambda k}\rd)\sum_{j=0}^{\infty}a_{\overline{\mu}}(j)e^{i\lambda j}
\\
\notag&=\ld(\sum_{k=0}^{\infty}\psi^{\top}_{\overline{\mu}}(k)e^{-i\lambda k}\rd)
\sum_{m=-\infty}^{\infty}\sum_{j=0}^{\infty}Z_{\overline{\mu}}(m-j)  a_{\overline{\mu}}(j)e^{i\lambda m}
\\
&=\ld(\sum_{k=0}^{\infty}\psi^{\top}_{\overline{\mu}}(k)e^{-i\lambda k}\rd)
\sum_{m=-\infty}^{\infty} e_{\overline{\mu}}(m)e^{i\lambda m}.
\label{simple_sp_char_part2}
\end{align}
 Making use of  (\ref{simple_sp_char_part1}) and (\ref{simple_sp_char_part2}),  formula (\ref{poh A2_e_st.n_d}) for the spectral characteristic $\overrightarrow{h}_{\overline{\mu}}^2(\lambda)$ of the optimal estimate $\widehat{\overrightarrow{A}}\eta$ can be presented as
\begin{eqnarray*}
 \vec h_{\overline{\mu}}^2(\lambda)&=&\frac
{\chi_{\overline{\mu}}^{(d)}(e^{-i\lambda})}{\beta^{(d)}(i\lambda)}
\ld(\sum_{k=0}^{\infty}\psi^{\top}_{\overline{\mu}}(k)e^{-i\lambda k}\rd)
\sum_{m=1}^{\infty} e_{\overline{\mu}}(-m)e^{-i\lambda m}
\\&=&\frac
{\chi_{\overline{\mu}}^{(d)}(e^{-i\lambda})}{\beta^{(d)}(i\lambda)}
\ld(\sum_{k=0}^{\infty}\psi^{\top}_{\overline{\mu}}(k)e^{-i\lambda k}\rd)
\sum_{m=1}^{\infty}\sum_{j=0}^{\infty}\sum_{l=0}^{\infty}\overline{\psi}_{\overline{\mu}}(l)\overline{g}(m+j+l)a_{\overline{\mu}}(j)e^{-i\lambda m}
\\&=&\frac
{\chi_{\overline{\mu}}^{(d)}(e^{-i\lambda})}{\beta^{(d)}(i\lambda)}\ld(\sum_{k=0}^{\infty}\psi^{\top}_{\overline{\mu}}(k)e^{-i\lambda k}\rd)\sum_{m=1}^{\infty}(\overline{\psi}_{\overline{\mu}} \me C_{\mu,g} )_m e^{-i\lambda m}
\\&=&\frac
{\chi_{\overline{\mu}}^{(d)}(e^{-i\lambda})}{\beta^{(d)}(i\lambda)}\Psi^{\top}_{\overline{\mu}}(e^{-i\lambda })C_{\mu,g} (e^{-i\lambda}),
\end{eqnarray*}
where $\overline{\psi}_{\overline{\mu}}=(\overline{\psi}_{\overline{\mu}}(0),\overline{\psi}_{\overline{\mu}}(1),\overline{\psi}_{\overline{\mu}}(2),\ldots)$,
\[
(\overline{\psi}_{\overline{\mu}} \me C_{\mu,g} )_m=\sum_{k=0}^{\infty}\overline{\psi}_{\overline{\mu}}(k)\me c_{\mu,g}(k+m),
\]
\[
\me c_{\mu,g}(m)=\sum_{k=0}^{\infty}\overline{g}(m+k)a_{\overline{\mu}}(k)
=\sum_{l=0}^{\infty}\overline{\phi}(l)\sum_{k=0}^{\infty}\phi^{\top}(l+m+k)a_{\overline{\mu}}(k)
=\sum_{l=0}^{\infty}\overline{\phi}(l)(\widetilde{\Phi}\me a_{\overline{\mu}})_{l+m},
\]
\[
(\widetilde{\Phi}\me a_{\overline{\mu}})_m=\sum_{k=0}^{\infty}\phi^{\top}(m+k)a_{\overline{\mu}}(k).
\]

Making use of the inverse operator $(\me
P_{\overline{\mu}})^{-1}$ factorization and following the transformation steps of \eqref{simple_sp_char_part1} obtain  the formula   for  the spectral characteristic $\vec h_{\overline{\mu}}^1(\lambda)$ of the optimal estimate $\widehat{\vec B}\xi$:

\begin{eqnarray*}
\vec{h}_{\overline{\mu}}^1(\lambda)
&=&\frac{\chi_{\overline{\mu}}^{(d)}(e^{-i\lambda})}{\beta^{(d)}(i\lambda)}
\Biggl(\vec{B}_{\overline{\mu}}(e^{i\lambda})
-
\frac{|\beta^{(d)}(i\lambda)|^2}
{|\chi_{\overline{\mu}}^{(d)}(e^{-i\lambda})|^2}
\ld[(f(\lambda)+|\beta^{(d)}(i\lambda)|^2g(\lambda))^{-1}\rd]^{\top}
\\
& &\times
\left(
\sum_{k=0}^{\infty}((\me P_{\overline{\mu}}^{-1}D^{\overline{\mu}}\me a)_k e^{ik\lambda}
\right)\Biggr)
\\
&=&\frac{\chi_{\overline{\mu}}^{(d)}(e^{-i\lambda})}{\beta^{(d)}(i\lambda)}
  \ld(\vec B_{\overline{\mu}}(e^{i\lambda})
 -\ld(\sum_{k=0}^{\infty}\psi_{\overline{\mu}}^{\top}(k)e^{-i\lambda k}\rd)
 \sum_{m=0}^{\infty}\sum_{p=0}^{\infty}\theta^{\top}_{\overline{\mu}}(p) \vec b_{\overline{\mu}}(b+m) e^{i\lambda m}\rd)
\\
&=&\frac{\chi_{\overline{\mu}}^{(d)}(e^{-i\lambda})}{\beta^{(d)}(i\lambda)}
  \ld(\vec B_{\overline{\mu}}(e^{i\lambda})
 -\Psi^{\top}_{\overline{\mu}}(e^{-i\lambda })\vec r_{\overline{\mu}} (e^{i\lambda})\rd)
   \\&=&\frac
{\chi_{\overline{\mu}}^{(d)}(e^{-i\lambda})}{\beta^{(d)}(i\lambda)}\Psi^{\top}_{\overline{\mu}}(e^{-i\lambda })\vec C_{\mu,1} (e^{-i\lambda}),
\end{eqnarray*}
where
\begin{eqnarray*}
\vec r_{\overline{\mu}} (e^{i\lambda})&=&
\sum_{m=0}^{\infty}( \theta^{\top}_{\overline{\mu}}D^{\mu}\me A)_m e^{i\lambda m}
 =
\sum_{m=0}^{\infty}\sum_{p=0}^{\infty}\theta^{\top}_{\overline{\mu}}(p) \vec b_{\overline{\mu}}(b+m) e^{i\lambda m},
\\
\vec C_{\mu,1} (e^{-i\lambda})&=&
\sum_{m=1}^{\infty}( \theta^{\top}_{\overline{\mu}}\widetilde{\me B}_{\overline{\mu}} )_m e^{-i\lambda m}
 =
\sum_{m=1}^{\infty}\sum_{p=m}^{\infty}\theta^{\top}_{\overline{\mu}}(p)\vec b_{\overline{\mu}}(p-m) e^{-i\lambda m}
 \\ &=&
\sum_{m=1}^{\infty}\sum_{p=0}^{\infty}\theta^{\top}_{\overline{\mu}}(m+p)\vec b_{\overline{\mu}}(p) e^{-i\lambda m},
\end{eqnarray*}
vector $\theta^{\top}_{\overline{\mu}}=((\theta_{\overline{\mu}}(0))^{\top},(\theta_{\overline{\mu}}(1))^{\top},(\theta_{\overline{\mu}}(2))^{\top},\ldots)$; $\me A$ is a linear symmetric operator determined by the matrix with the vector entries $(\me A)_{k,j}=\vec a(k+j)$, $k,j\geq0$; $\widetilde{\me B}_{\overline{\mu}}$ is a linear operator, which is determined by a matrix with the vector entries
 $( \widetilde{\me B}_{\overline{\mu}})_{k,j}=\vec b_{\overline{\mu}}(k-j)$ for $0\leq j\leq k$, $(
\widetilde{\me B}_{\overline{\mu}})_{k,j}=0$ for $0\leq k<j$.

Then the spectral characteristic $\vec h_{\overline{\mu}}(\lambda)$ of the estimate $\widehat{A}\xi$ can be calculated by the formula
\begin{eqnarray}
 \notag \vec h_{\overline{\mu}}(\lambda)&=&\frac{\chi_{\overline{\mu}}^{(d)}(e^{-i\lambda})}{\beta^{(d)}(i\lambda)}
\ld(\sum_{k=0}^{\infty}\psi^{\top}_{\overline{\mu}}(k) e^{-i\lambda k}\rd)
\sum_{m=1}^{\infty}\ld( \theta^{\top}_{\overline{\mu}}\widetilde{\me B}_{\overline{\mu}} -\overline{\psi}_{\overline{\mu}} \me C_{\mu,g}\rd)_m e^{-i\lambda m}
\\
\notag&=&\frac
{\chi_{\overline{\mu}}^{(d)}(e^{-i\lambda})}{\beta^{(d)}(i\lambda)}\Psi^{\top}_{\overline{\mu}}(e^{-i\lambda })\ld(\vec C_{\mu,1} (e^{-i\lambda})-\vec C_{\mu,g} (e^{-i\lambda})\rd)
 \\&=&\vec B_{\overline{\mu}}(e^{i\lambda})\frac{\chi_{\overline{\mu}}^{(d)}(e^{-i\lambda})}{\beta^{(d)}(i\lambda)}-\widetilde{h}_{\overline{\mu}}(\lambda),
 \label{simple_spectr A_e_st.n_d}\end{eqnarray}
\begin{eqnarray*}
\widetilde{\vec h}_{\overline{\mu}}(\lambda)&=&\frac{\chi_{\overline{\mu}}^{(d)}(e^{-i\lambda})}{\beta^{(d)}(i\lambda)}
\ld(\sum_{k=0}^{\infty}\psi^{\top}_{\overline{\mu}}(k) e^{-i\lambda k}\rd)\\
& &\times\ld(\sum_{m=0}^{\infty}( \theta^{\top}_{\overline{\mu}}D^{\mu}\me A)_m e^{i\lambda m}+\sum_{m=1}^{\infty}(\overline{\psi}_{\overline{\mu}} \me C_{\mu,g} )_m e^{-i\lambda m}\rd)
\\&=&
\frac{\chi_{\overline{\mu}}^{(d)}(e^{-i\lambda})}{\beta^{(d)}(i\lambda)}
\Psi^{\top}_{\overline{\mu}}(e^{-i\lambda })\ld(\vec r_{\overline{\mu}} (e^{i\lambda})+\vec C_{\mu,g} (e^{-i\lambda})\rd),
\end{eqnarray*}

The value of the mean square error of the estimate $\widehat{A}\xi$ is calculated by the formula
\begin{eqnarray}
 \notag \Delta\ld(f,g;\widehat{A}\xi\rd)&=&\Delta\ld(f,g;\widehat{H}\xi\rd)= \mt E\ld|H\xi-\widehat{H}\xi\rd|^2
 \\\notag&=&\frac{1}{2\pi}\int_{-\pi}^{\pi}(\vec A(e^{i\lambda}))^{\top}g(\lambda)\overline{\vec A(e^{i\lambda})}d\lambda
\\ \notag &&+
 \frac{1}{2\pi}\int_{-\pi}^{\pi}(\widetilde{\vec h}_{\overline{\mu}}(e^{i\lambda}))^{\top}(f(\lambda)
 +|\beta^{(d)}(i\lambda)|^2g(\lambda))\overline{\widetilde{\vec h}_{\overline{\mu}}(e^{i\lambda})}d\lambda
\\ \notag &&-\frac{1}{2\pi}\int_{-\pi}^{\pi}\frac{\beta^{(d)}(i\lambda)}{\chi_{\overline{\mu}}^{(d)}(e^{-i\lambda})}(\widetilde{h}_{\overline{\mu}}(e^{i\lambda}))^{\top}
g(\lambda)\overline{A_{\overline{\mu}}(e^{i\lambda})}d\lambda
\\\notag && -
 \frac{1}{2\pi}\int_{-\pi}^{\pi}\frac{\overline{\beta^{(d)}(i\lambda)}}{\overline{\chi_{\overline{\mu}}^{(d)}(e^{-i\lambda})}}(A_{\overline{\mu}}(e^{i\lambda}))^{\top}
g(\lambda)
\overline{\widetilde{\vec h}_{\overline{\mu}}(e^{i\lambda})}d\lambda
  \\\notag
 &=&\|\Phi^{\top}\me a_{\overline{\mu}}\|^2+\|\widetilde{\Phi}\me a_{\overline{\mu}}\|_1^2
 +\ld\langle \theta^{\top}_{\overline{\mu}}D^{\mu}\me A-\overline{\psi}_{\overline{\mu}} \me C_{\mu,g},\theta^{\top}_{\overline{\mu}}D^{\mu}\me A\rd\rangle
 \\ && -\ld\langle\theta^{\top}_{\overline{\mu}}D^{\mu}\me A,\me Z_{\overline{\mu}}\me a_{\overline{\mu}}\rd\rangle
 -\ld\langle \me Z_{\overline{\mu}}\me a_{\overline{\mu}},\overline{\psi}_{\overline{\mu}} \me C_{\mu,g}\rd\rangle_1,
\label{simple_poh A_e_st.n_d}\end{eqnarray}
where $\|\vec x\|_1^2=\langle \vec x, \vec x\rangle_1$, $\langle \vec x, \vec y\rangle_1=\sum_{k=1}^{\infty}(\vec x(k))^{\top}\overline{\vec y}(k)$ for the vectors $\vec x=((\vec x(0))^{\top},(\vec x(1))^{\top},(\vec x(2))^{\top},\ldots)^{\top}$, $\vec y=((\vec y(0))^{\top},(\vec y(1))^{\top},(\vec y(2))^{\top},\ldots)^{\top}$.

The obtained results are summarized in the form of the following theorem.

\begin{thm}\label{thm3_e_st.n_d}
Let $\vec \xi(m)$, $m\in\mr Z$, be a vector stochastic sequence which determines the stationary GM increment sequence $\chi_{\overline{\mu},\overline{s}}^{(d)}(\vec{\xi}(m))$ and let $\vec \eta(m)$, $m\in\mr Z$, be vector stationary stochastic  sequence, uncorrelated with $\vec \xi(m)$. Suppose that the coefficients
$\vec a(k)$, $k\geq0$,
satisfy conditions (\ref{umovana a_e_d}) -- (\ref{umovana a_mu_e_d}) and the spectral densities $f(\lambda)$ and $g(\lambda)$ of the stochastic sequences $\vec \xi(m)$ and $\vec \eta(m)$ admit the canonical factorizations (\ref{fakt1}) -- (\ref{fakt2}). Then the spectral characteristic
$\vec h_{\overline{\mu}}(\lambda)$ and the value of the mean square error $\Delta(f,g;\widehat{A}\vec \xi)$ of the optimal estimate $\widehat{A}\vec \xi$ of the functional $A\vec \xi$ based on observations of the sequence
$\vec \xi (m)+\vec \eta (m)$ at points $m=-1,-2,\ldots$ can be calculated by formulas (\ref{simple_spectr A_e_st.n_d}) and (\ref{simple_poh A_e_st.n_d}) respectively.
\end{thm}

\subsection{Forecasting of stochastic sequences with periodically stationary increment}

Consider the problem of mean square optimal linear estimation of the functionals
\begin{equation}
A{\vartheta}=\sum_{k=0}^{\infty}{a}^{(\vartheta)}(k)\vartheta(k), \quad
A_{M}{\vartheta}=\sum_{k=0}^{N}{a}^{(\vartheta)}(k)\vartheta(k)
\end{equation}
which depend on unobserved values of the stochastic sequence ${\vartheta}(m)$ with periodically stationary
increments. Estimates are based on observations of the sequence $\zeta(m)=\vartheta(m)+\eta(m)$ at points $m=-1,-2,\ldots$.

The functional $A{\vartheta}$ can be represented in the form
\begin{eqnarray}
\nonumber
A{\vartheta}& = &\sum_{k=0}^{\infty}{a}^{(\vartheta)}(k)\vartheta(k)=\sum_{m=0}^{\infty}\sum_{p=1}^{T}
{a}^{(\vartheta)}(mT+p-1)\vartheta(mT+p-1)
\\\nonumber
& = & \sum_{m=0}^{\infty}\sum_{p=1}^{T}a_p(m)\xi_p(m)=\sum_{m=0}^{\infty}(\vec{a}(m))^{\top}\vec{\xi}(m)=A\vec{\xi},
\end{eqnarray}
where
\begin{multline} \label{zeta}
\vec{\xi}(m)=({\xi}_1(m),{\xi}_2(m),\dots,{\xi}_T(m))^{\top},\hfill
 {\xi}_p(m)=\vartheta(mT+p-1);\,p=1,2,\dots,T;
\end{multline}
\begin{multline} \label{azeta}
 \vec{a}(m) =({a}_1(m),{a}_2(m),\dots,{a}_T(m))^{\top},\hfill
 {a}_p(m)=a^{(\vartheta)}(mT+p-1);\,p=1,2,\dots,T.
\end{multline}

\begin{thm}
\label{thm_est_Azeta}
Let a stochastic sequence ${\vartheta}(k)$ with periodically stationary increments generate by formula \eqref{zeta}
 a vector-valued stochastic sequence $\vec{\xi}(m) $ which determine a
stationary stochastic GM increment sequence
$\chi_{\overline{\mu},\overline{s}}^{(d)}(\vec{\xi}(m))$ with the spectral density matrix $f(\lambda)$.
Let $\vec\eta(m)$, $m\in\mr Z$, $\vec{\eta}(m)=({\eta}_1(m),{\eta}_2(m),\dots,{\eta}_T(m))^{\top},\,
 {\eta}_p(m)=\eta(mT+p-1);\,p=1,2,\dots,T, $
be uncorrelated with the sequence $\vec\xi(m)$ stationary stochastic
sequence with an absolutely continuous spectral function
$G(\lambda)$ which has spectral density matrix $g(\lambda)$. Let the minimality condition
(\ref{umova11_e_st.n_d}) be satisfied.
Let coefficients $\vec {a}(k), k\geqslant 0$ determined by formula \eqref{azeta}  satisfy conditions  (\ref{umovana a_e_d}) -- (\ref{umovana a_mu_e_d}).
Then the optimal linear estimate $\widehat{A}\vartheta$ of the functional $A\vartheta$ based on observations of the sequence
$\zeta(m)=\vartheta(m)+\eta(m)$ at points $m=-1,-2,\ldots$ is calculated by formula (\ref{otsinka A_e_d}).
The spectral characteristic
$\vec h_{\overline{\mu}}(\lambda)=\{h_{p}(\lambda)\}_{p=1}^{T}$ and the value of the mean square error $\Delta(f;\widehat{A}\zeta)$ of the optimal estimate $\widehat{A}\zeta$ are calculated by formulas
(\ref{spectr A_e_d}) and (\ref{pohybkaA}) respectively.
\end{thm}

The functional $A_M{\vartheta}$ can be represented in the form
\begin{eqnarray}
\nonumber
A_M{\vartheta}& = &\sum_{k=0}^{M}{a}^{(\vartheta)}(k)\vartheta(k)=\sum_{m=0}^{N}\sum_{p=1}^{T}
{a}^{(\vartheta)}(mT+p-1)\vartheta(mT+p-1)
\\\nonumber
& = &\sum_{m=0}^{N}\sum_{p=1}^{T}a_p(m)\xi_p(m)=\sum_{m=0}^{N}(\vec{a}(m))^{\top}\vec{\xi}(m)=A_N\vec{\xi},
\end{eqnarray}
where $N=[\frac{M}{T}]$, the sequence $\vec{\xi}(m) $ is determined by formula \eqref{zeta},
\begin{eqnarray}
\nonumber
(\vec{a}(m))^{\top}& = &({a}_1(m),{a}_2(m),\dots,{a}_T(m))^{\top},
\\\nonumber
 {a}_p(m)& = &a^{\vartheta}(mT+p-1);\,0\leq m\leq N; 1\leq p\leq T;\,mT+p-1\leq M;
\\  {a}_p(N)& = &0;\quad
M+1\leq NT+p-1\leq (N+1)T-1;1\leq p\leq T. \label{aNzeta}
\end{eqnarray}

Making use of the introduced notations and statements of Theorem \ref{thm2AN} we can claim that the following theorem holds true.

\begin{thm}
\label{thm_est_A_Nzeta}
Let a stochastic sequence ${\vartheta}(k)$ with periodically stationary increments generate by formula \eqref{zeta}
 a vector-valued stochastic sequence $\vec{\xi}(m) $ which determine a
stationary stochastic GM increment sequence
$\chi_{\overline{\mu},\overline{s}}^{(d)}(\vec{\xi}(m))$ with the spectral density matrix $f(\lambda)$.
Let $\{\vec\eta(m),m\in\mr Z\}$, $\vec{\eta}(m)=({\eta}_1(m),{\eta}_2(m),\dots,{\eta}_T(m))^{\top},\,
 {\eta}_p(m)=\eta(mT+p-1);\,p=1,2,\dots,T, $
be uncorrelated with the sequence $\vec\xi(m)$ stationary stochastic
sequence with an absolutely continuous spectral function
$G(\lambda)$ which has spectral density matrix $g(\lambda)$. Let the minimality condition
(\ref{umova11_e_st.n_d}) be satisfied.
Let coefficients $\vec {a}(k), k\geqslant 0$ be determined by formula \eqref{aNzeta}.
The optimal linear estimate $\widehat{A}_M\zeta$ of the functional $A_M\zeta=A_N\vec{\xi}$ based on observations of the sequence
$\zeta(m)=\vartheta(m)+\eta(m)$ at points $m=-1,-2,\ldots$ is calculated by formula \eqref{otsinka A_N_e_st.n_d}.
The spectral characteristic $\vec{h}_{\overline{\mu},N}(\lambda)=\{h_{\overline{\mu},N,p}(\lambda)\}_{p=1}^{T}$ and the value of the mean square error $\Delta(f;\widehat{A}_M\zeta)$
are calculated by formulas  \eqref{spectr A_e_dN} and \eqref{pohybka33} respectively.
\end{thm}

As a corollary from the proposed theorem, one can obtain the mean square optimal estimate of the unobserved value
$\vartheta(M)$, $M\geq0$ of a stochastic sequence ${\vartheta}(m)$ with periodically stationary increments
based on observations of the sequence $\zeta(m)=\vartheta(m)+\eta(m)$ at points $m=-1,-2,\ldots$
Making use of the notations
$\vartheta(M)=\vartheta_p(N)=(\vec\xi(N))^{\top}\boldsymbol{\delta}_p$, $N=[\frac{M}{T}]$, $p=M+1-NT$,
and the obtained results we can conclude that the following corollary holds true.

\begin{nas}\label{nas zeta_e_d}
Let a stochastic sequence ${\vartheta}(k)$ with periodically stationary increments generate by formula \eqref{zeta}
 a vector-valued stochastic sequence $\vec{\xi}(m) $ which determine the
stationary stochastic GM increment sequence
$\chi_{\overline{\mu},\overline{s}}^{(d)}(\vec{\xi}(m))$ with the spectral density matrix $f(\lambda)$.
Let $\vec\eta(m)$,  $m\in\mr Z$, $\vec{\eta}(m)=({\eta}_1(m),{\eta}_2(m),\dots,{\eta}_T(m))^{\top},\,
 {\eta}_p(m)=\eta(mT+p-1);\,p=1,2,\dots,T, $
be uncorrelated with the sequence $\vec\xi(m)$ stationary stochastic
sequence with an absolutely continuous spectral function
$G(\lambda)$ which has spectral density matrix $g(\lambda)$. Let the minimality condition
(\ref{umova11_e_st.n_d}) be satisfied.
The optimal linear estimate $\widehat{\vartheta}(M)$ of the unobserved value
$\vartheta(M)$, $M\geq0$ of a stochastic sequence ${\vartheta}(m)$ with periodically stationary increments
based on observations of the sequence $\zeta(m)=\vartheta(m)+\eta(m)$ at points $m=-1,-2,\ldots$
is calculated by formula \eqref{est_xi_N}.
The spectral characteristic $\vec h_{\overline{\mu},N,p}(\lambda)$ of the estimate is calculated by the formula
 \eqref{sph_est_xi_N}.
The value of the mean square error of the optimal estimate is calculated by the formula
 \eqref{poh xi_p_e_st.n_d}.
\end{nas}

\subsection{Forecasting of one class of cointegrated vector stochastic sequences}\label{classical_extrap_coint}

Consider two seasonal  vector stochastic sequences $\{\vec \xi(m),m\in\mr Z\}$ and $\{\vec \zeta(m),m\in\mr Z\}$ with absolutely continuous spectral functions
$F(\lambda)$ and $P(\lambda)$ and spectral densities
$f(\lambda)$ and $p(\lambda)$ respectively. Assume that both of them have the same order
$d$ and  seasonal  vector $\overline{s}$.

\begin{ozn}\label{ozn_coint}
Within this subsection, a pair of  seasonal  vector stochastic sequences $\{\vec\xi(m),m\in\mr Z\}$ and
$\{\vec\zeta(m),m\in\mr Z\}$ are called seasonally cointegrated if there exists a
constant $\alpha\neq0$ such that the linear combination  sequence
$\vec\zeta(m)-\alpha\vec \xi(m)$ is a stationary vector stochastic sequence.
\end{ozn}

Under the forecast of two seasonally cointegrated stochastic sequences
we understand the mean-square optimal linear estimates
of the functionals
\[A\vec{\xi}=\sum_{k=0}^{\infty}(\vec{a}(k))^{\top}\vec{\xi}(k), \quad
A_{N}\vec{\xi}=\sum_{k=0}^{N}(\vec{a}(k))^{\top}\vec{\xi}(k),\]
which depend on the unobserved values of the stochastic sequence $\vec \xi(m)$ based on observations of the stochastic sequence $\vec \zeta(m)$ at points $m=-1,-2,\ldots$.
Applying the results of Subsection \ref{classical_extrapolation}, the   forecasts can be found  under an assumption  that the vector sequences $\vec\xi(m)$ and $\vec\zeta(m)-\alpha\vec\xi(m)$
are uncorrelated.

Let the minimality condition holds true:
\be
 \ip \text{Tr}\left[ \frac{|\beta^{(d)}(i\lambda)|^2}{|\chi_{\overline{\mu}}^{(d)}(e^{-i\lambda})|^2}p(\lambda)^{-1}\right]
 d\lambda<\infty.\label{umova111_e_st.n_d}
\ee
Determine operators $\me P_{\overline{\mu}}^{\alpha}$, $\me T_{\overline{\mu}}^{\alpha}$,
$\me Q^{\alpha}$ with the help of the Fourier coefficients
\[
P_{k,j}^{\overline{\mu},\alpha}=\frac{1}{2\pi}\int_{-\pi}^{\pi} e^{-i\lambda (k-j)}
 \frac{|\beta^{(d)}(i\lambda)|^2}{|\chi_{\overline{\mu}}^{(d)}(e^{-i\lambda})|^2}
\ld[\ld(p(\lambda)\rd)^{-1}\rd]^{\top}
d\lambda;
\]
\[
T^{\overline{\mu},\beta}_{k,j}=\frac{1}{2\pi}\int_{-\pi}^{\pi}
e^{-i\lambda (k-j)}
 \frac{1}{|\chi_{\overline{\mu}}^{(d)}(e^{-i\lambda})|^2}\ld[(p(\lambda)-\alpha^2f(\lambda))
\ld(p(\lambda)\rd)^{-1}\rd]^{\top}
d\lambda;
\]
 \[
 Q_{k,j}^{,\alpha}=\frac{1}{2\pi}\int_{-\pi}^{\pi}
e^{-i\lambda (k-j)}\frac{1}{|\beta^{(d)}(i\lambda)|^2}\ld[f(\lambda)\ld(p(\lambda)\rd)^{-1}(p(\lambda)-\alpha^2f(\lambda))\rd]^{\top}
d\lambda.
\]
of the functions
\be\label{functions for
beta-operators}
 \frac{|\beta^{(d)}(i\lambda)|^2}{|\chi_{\overline{\mu}}^{(d)}(e^{-i\lambda})|^2}
\ld[\ld(p(\lambda)\rd)^{-1}\rd]^{\top},
 \,
 \frac{1}{|\chi_{\overline{\mu}}^{(d)}(e^{-i\lambda})|^2}\ld[(p(\lambda)-\alpha^2f(\lambda))
\ld(p(\lambda)\rd)^{-1}\rd]^{\top},
\ee
\[
 \frac{1}{|\chi_{\overline{\mu}}^{(d)}(e^{-i\lambda})|^2}\ld[f(\lambda)
\ld(p(\lambda)\rd)^{-1}(p(\lambda)-\alpha^2f(\lambda))\rd]^{\top}
\]
 in the same way as we defined operators $\me P_{\overline{\mu}}$, $\me T_{\overline{\mu}}$, $\me Q$. Theorem $\ref{thm1_e_st.n_d}$ implies that the spectral characteristic
$h^{\alpha}_{\overline{\mu}}(\lambda)$ of the optimal estimate
 \be \label{otsinka A_co_e_st.n_d} \widehat{A}\vec\xi=\ip
 (h^{\alpha}_{\overline{\mu}}(\lambda))^{\top}d\vec Z_{\zeta^{(d)}}(\lambda)-\sum_{k=-n(\nu)}^{-1}(\vec v_{\overline{\mu}}(k))^{\top}\zeta(k), \ee
of the functional $A\vec\xi$ is calculated by the formula

\begin{multline}\label{spectr A_co_e_st.n_d}
(\vec{h}_{\overline{\mu}}^{\alpha}(\lambda))^{\top}=(\vec{B}_{\overline{\mu}}(e^{i\lambda}))^{\top}
\frac{\chi_{\overline{\mu}}^{(d)}(e^{-i\lambda})}{\beta^{(d)}(i\lambda)}-
\\
-
\frac{1}{\overline{\chi_{\overline{\mu}}^{(d)}(e^{-i\lambda})}\beta^{(d)}(i\lambda)}
(\vec{A}_{\overline{\mu}}(e^{i\lambda}))^{\top}(p(\lambda)-\alpha^2f(\lambda))
(p(\lambda))^{-1}
-
\\
-
\frac{\overline{\beta^{(d)}(i\lambda)}
}
{\overline{\chi_{\overline{\mu}}^{(d)}(e^{-i\lambda})}}\left(
\sum_{k=0}^{\infty}((\me P^{\alpha}_{\overline{\mu}})^{-1}D^{\overline{\mu}}\me a-(\me P^{\alpha}_{\overline{\mu}})^{-1}\me T^{\alpha}_{\overline{\mu}}\me a_{\overline{\mu}})_k e^{ik\lambda}
\right)^{\top}(p(\lambda))^{-1}.
\end{multline}

The value of the mean square error of the estimate $\widehat{A}\vec\xi$ is calculated by the formula
\begin{multline*}
\Delta(f,g;\widehat{A}\vec\xi)=
\frac{1}{2\pi}\int_{-\pi}^{\pi}
\Biggl[\frac{1}{|\beta^{(d)}(i\lambda)|^2}(\vec{A}_{\overline{\mu}}(e^{i\lambda}))^{\top}(p(\lambda)-\alpha^2f(\lambda)) +
\\
+
\left(\sum_{k=0}^{\infty}((\me P^{\alpha}_{\overline{\mu}})^{-1}D^{\overline{\mu}}\me a-(\me P^{\alpha}_{\overline{\mu}})^{-1}\me T^{\alpha}_{\overline{\mu}}\me a_{\overline{\mu}})_k e^{ik\lambda}
\right)^{\top}
\Biggr]
\\
\times
\frac{|\beta^{(d)}(i\lambda)|^2}{|\chi_{\overline{\mu}}^{(d)}(e^{-i\lambda})|^2}(p(\lambda))^{-1}\, f(\lambda)\, (p(\lambda))^{-1}
\times
\end{multline*}
\[
\times
\left[\frac{1}{|\beta^{(d)}(i\lambda)|^2}
(p(\lambda)-\alpha^2f(\lambda))\overline{\vec{A}(e^{i\lambda})} +
\overline{\sum_{k=0}^{\infty}((\me P^{\alpha}_{\overline{\mu}})^{-1}D^{\overline{\mu}}\me a-(\me P^{\alpha}_{\overline{\mu}})^{-1}\me T^{\alpha}_{\overline{\mu}}\me a_{\overline{\mu}})_k e^{ik\lambda}}
\right]
d\lambda+
\]
\[
+\frac{1}{2\pi}\int_{-\pi}^{\pi}
\left[(\vec{A}_{\overline{\mu}}(e^{i\lambda}))^{\top}f(\lambda) -|\beta^{(d)}(i\lambda)|^2
\left(\sum_{k=0}^{\infty}((\me P^{\alpha}_{\overline{\mu}})^{-1}D^{\overline{\mu}}\me a-(\me P^{\alpha}_{\overline{\mu}})^{-1}\me T^{\alpha}_{\overline{\mu}}\me a_{\overline{\mu}})_k e^{ik\lambda}
\right)^{\top}
\right]
\]
\[
\times
\frac{1}{|\chi_{\overline{\mu}}^{(d)}(e^{-i\lambda})|^2}\frac{1}{|\beta^{(d)}(i\lambda)|^2}(p(\lambda))^{-1}\,(p(\lambda)-\alpha^2f(\lambda))\,(p(\lambda))^{-1}
\times
\]
\[
\times
\left[f(\lambda)\overline{\vec{A}_{\overline{\mu}}(e^{i\lambda})} -|\beta^{(d)}(i\lambda)|^2
\overline{\sum_{k=0}^{\infty}((\me P^{\alpha}_{\overline{\mu}})^{-1}D^{\overline{\mu}}\me a-(\me P^{\alpha}_{\overline{\mu}})^{-1}\me T^{\alpha}_{\overline{\mu}}\me a_{\overline{\mu}})_k e^{ik\lambda}}
\right]
d\lambda=
\]
\be\label{poh A_co_e_st.n_d}
=\ld\langle D^{\overline{\mu}}\me a- \me T^{\alpha}_{\overline{\mu}}\me
 a_{\overline{\mu}},(\me P^{\alpha}_{\overline{\mu}})^{-1}D^{\overline{\mu}}\me a-(\me P^{\alpha}_{\overline{\mu}})^{-1}\me T^{\alpha}_{\overline{\mu}}\me a_{\overline{\mu}}\rd\rangle+\ld\langle\me Q^{\alpha}\me a,\me
 a\rd\rangle.
 \ee

\begin{thm}\label{thm4_e_st.n_d}
Let $\vec \xi(m)$, $m\in\mr Z$, and $\vec \zeta(m)$, $m\in\mr Z$ be seasonally cointegrated stochastic sequences with the spectral densities
 $f(\lambda)$ and $p(\lambda)$ respectively. Suppose that the spectral density
 $p(\lambda)$ satisfy the minimality condition
(\ref{umova111_e_st.n_d}) and the coefficients $\vec a(k)$, $k\geq0$,
satisfy conditions
(\ref{umovana a_e_d}) -- (\ref{umovana a_mu_e_d}). If the stochastic sequences $\xi(m)$ and
$\vec \zeta(m)-\alpha\vec \xi(m)$ are uncorrelated, then the spectral characteristic $\vec h^{\alpha}_{\overline{\mu}}(\lambda)$ and the value of the mean square error $\Delta(f,g;\widehat{A}\xi)$ of the optimal estimate $\widehat{A}\vec \xi$ (\ref{otsinka A_co_e_st.n_d}) of the functional $A\vec \xi$ based on observations of the sequence
$\vec \zeta(m)$ at points $m=-1,-2,\ldots$ are calculated by formulas
 (\ref{spectr A_co_e_st.n_d}) and
 (\ref{poh A_co_e_st.n_d}) respectively.
\end{thm}

Let operators $\me P_{\overline{\mu}}^{\alpha}$, $\me T_{\overline{\mu},N}^{\alpha}$, $\me
Q_{N}^{\alpha}$ be defined by the Fourier coefficients of the functions (\ref{functions for beta-operators}) in the same way as we defined
 operators $\me P_{\overline{\mu}}$, $\me T_{\overline{\mu},N}$, $\me Q_{N}$.
Theorem \ref{thm2AN} implies that the spectral characteristic
$\vec h^{\alpha}_{\overline{\mu},N}(\lambda)$ of the optimal estimate
 \be \label{otsinka A_N_co_e_st.n_d} \widehat{A}_N\vec\xi=\ip
 \vec h^{\alpha}_{\overline{\mu},N}(\lambda)d\vec Z_{\zeta^{(n)}}(\lambda)-\sum_{k=-n(\gamma)}^{-1}(\vec v_{\overline{\mu},N}(k))^{\top}\vec \zeta(k) \ee
of the functional $A_N\vec \xi$ is calculated by the formula
\begin{multline}\label{spectr A_N_co_e_st.n_d}
(\vec{h}^{\alpha}_{\overline{\mu},N}(\lambda))^{\top}
=(\vec{B}_{\overline{\mu},N}(e^{i\lambda}))^{\top}
\frac{\chi_{\overline{\mu}}^{(d)}(e^{-i\lambda})}{\beta^{(d)}(i\lambda)}
-
\\-
\frac{1}{\overline{\chi_{\overline{\mu}}^{(d)}(e^{-i\lambda})}\beta^{(d)}(i\lambda)}
(\vec{A}_{\overline{\mu},N}(e^{i\lambda}))^{\top}(p(\lambda)-\alpha^2f(\lambda))
(p(\lambda))^{-1}
-
\\
-
\frac{\overline{\beta^{(d)}(i\lambda)}
}
{\overline{\chi_{\overline{\mu}}^{(d)}(e^{-i\lambda})}}\left(\sum_{k=0}^{\infty}((\me P^{\alpha}_{\overline{\mu}})^{-1}D^{\overline{\mu}}_N\me a_N-(\me P^{\alpha}_{\overline{\mu}})^{-1}\me T^{\alpha}_{\overline{\mu},N}\me a_{\overline{\mu},N})_k e^{ik\lambda}
\right)^{\top}(p(\lambda))^{-1}.
\end{multline}

The value of the mean square error of the estimate $\widehat{A}_N\xi$ is calculated by the formula
\begin{multline*}
\Delta(f,g;\widehat{A}_N\vec\xi)=
\frac{1}{2\pi}\int_{-\pi}^{\pi}
\Biggl[\frac{1}{|\beta^{(d)}(i\lambda)|^2}(\vec{A}_{\overline{\mu},N}(e^{i\lambda}))^{\top}(p(\lambda)-\alpha^2f(\lambda)) +
\\
+\left(\sum_{k=0}^{\infty}((\me P^{\alpha}_{\overline{\mu}})^{-1}D^{\overline{\mu}}_N\me a_N-(\me P^{\alpha}_{\overline{\mu}})^{-1}\me T^{\alpha}_{\overline{\mu},N}\me a_{\overline{\mu},N})_k e^{ik\lambda}
\right)^{\top}
\Biggr]
\times
\\
\times
\frac{|\beta^{(d)}(i\lambda)|^2}{|\chi_{\overline{\mu}}^{(d)}(e^{-i\lambda})|^2}(p(\lambda))^{-1}\, f(\lambda)\, (p(\lambda))^{-1}
\times
\\
\times
\Biggl[\frac{1}{|\beta^{(d)}(i\lambda)|^2}(p(\lambda)-\alpha^2f(\lambda))\overline{\vec{A}_{\overline{\mu},N}(e^{i\lambda})}
 +
 \\
 +
\overline{\sum_{k=0}^{\infty}((\me P^{\alpha}_{\overline{\mu}})^{-1}D^{\overline{\mu}}_N\me a_N- (\me P^{\alpha}_{\overline{\mu}})^{-1}\me T^{\alpha}_{\overline{\mu},N}\me a_{\overline{\mu},N})_k e^{ik\lambda}}
\Biggr]
d\lambda+
\end{multline*}
\begin{multline*}
+\frac{1}{2\pi}\int_{-\pi}^{\pi}
\Biggl[\frac{1}{|\beta^{(d)}(i\lambda)|^2}(\vec{A}_{\overline{\mu},N}(e^{i\lambda}))^{\top}f(\lambda) -
\\
-
\left(\sum_{k=0}^{\infty}((\me P^{\alpha}_{\overline{\mu}})^{-1}D^{\overline{\mu}}_N\me a_N-(\me P^{\alpha}_{\overline{\mu}})^{-1}\me T^{\alpha}_{\overline{\mu},N}\me a_{\overline{\mu},N})_k e^{ik\lambda}
\right)^{\top}
\Biggr]
\times
\\
\times
\frac{|\beta^{(d)}(i\lambda)|^2}{|\chi_{\overline{\mu}}^{(d)}(e^{-i\lambda})|^2}
(p(\lambda))^{-1}\,(p(\lambda)-\alpha^2f(\lambda))\,(p(\lambda))^{-1}
\times
\\
\times
\left[\frac{1}{|\beta^{(d)}(i\lambda)|^2}f(\lambda)\overline{\vec{A}_{\overline{\mu},N}(e^{i\lambda})}
 -
\overline{\sum_{k=0}^{\infty}((\me P^{\alpha}_{\overline{\mu}})^{-1}D^{\overline{\mu}}_N\me a_N-(\me P^{\alpha}_{\overline{\mu}})^{-1}\me T^{\alpha}_{\overline{\mu},N}\me a_{\overline{\mu},N})_k e^{ik\lambda}}
\right]
d\lambda=
\end{multline*}
\be\label{poh A_N_co_e_st.n_d}
=\ld\langle D^{\overline{\mu}}_N\me a_N- \me T^{\alpha}_{\overline{\mu},N}\me
 a_{\overline{\mu},N},(\me P^{\alpha}_{\overline{\mu}})^{-1}D^{\overline{\mu}}_N\me a_N-(\me P^{\alpha}_{\overline{\mu}})^{-1}\me T^{\alpha}_{\overline{\mu},N}\me a_{\overline{\mu},N}\rd\rangle+\ld\langle\me Q^{\alpha}_N\me a_N,\me
 a_N\rd\rangle.
 \ee

\begin{thm}\label{thm5_e_st.n_d}
Let $\vec \xi(m)$, $m\in\mr Z$ and $\vec \zeta(m)$, $m\in\mr Z$ be seasonally cointegrated stochastic sequences with the spectral densities
 $f(\lambda)$ and $p(\lambda)$ respectively. Suppose that the spectral density
 $p(\lambda)$ satisfy the minimality condition
(\ref{umova111_e_st.n_d}). If the stochastic sequences $\vec \xi(m)$ and
$\vec \zeta(m)-\vec \alpha\xi(m)$ are uncorrelated, then the spectral characteristic $\vec h^{\alpha}_{\overline{\mu},N}(\lambda)$ and the value of the mean square error $\Delta(f,g;\widehat{A}_N\vec \xi)$ of the optimal estimate $\widehat{A}_N\vec \xi$ (\ref{otsinka A_N_co_e_st.n_d}) of the functional $A_N\vec \xi$ based on observations of the sequence
$\vec \zeta(m)$ at points $m=-1,-2,\ldots$ are calculated by formulas (\ref{spectr A_N_co_e_st.n_d}) and
(\ref{poh A_N_co_e_st.n_d}) respectively.
\end{thm}

Let the spectral densities $f(\lambda)$ and $p(\lambda)$ admit the factorizations
\be \label{fakt1_coi}
 \dfrac{|\chi_{\overline{\mu}}^{(d)}(e^{-i\lambda})|^2}{|\beta^{(d)}(i\lambda)|^2}
 p(\lambda)
 =\Theta^{\alpha}_{\overline{\mu}}(e^{-i\lambda})(\Theta^{\alpha}_{\overline{\mu}}(e^{-i\lambda}))^*,\quad \Theta^{\alpha}_{\overline{\mu}}(e^{-i\lambda})=\sum_{k=0}^{\infty}\theta^{\alpha}_{\overline{\mu}}(k)e^{-i\lambda k},\ee
\be \label{fakt2_coi}
\dfrac{|\beta^{(d)}(i\lambda)|^2}{|\chi_{\overline{\mu}}^{(d)}(e^{-i\lambda})|^2}
 (p(\lambda))^{-1} =
 (\Psi_{\overline{\mu}}^{\alpha}(e^{-i\lambda}))^*\Psi_{\overline{\mu}}^{\alpha}(e^{-i\lambda}), \quad
 \Psi_{\overline{\mu}}^{\alpha}(e^{-i\lambda})=\sum_{k=0}^{\infty}\psi_{\overline{\mu}}^{\alpha}(k)e^{-i\lambda k},\ee
 \be \label{fakt3_coi}
|\beta^{(d)}(i\lambda)|^{-2}( p(\lambda)-\alpha^2f(\lambda))=\Phi^{\alpha}(e^{-i\lambda})(\Phi^{\alpha}(e^{-i\lambda}))^*, \quad
 \Phi^{\alpha}(e^{-i\lambda})=\sum_{k=0}^{\infty}\phi^{\alpha}(k)e^{-i\lambda k}.\ee
Let operators and vectors $\me C^{\alpha}_{\mu,g}$, $\widetilde{\Phi}^{\alpha}$,  $\Phi^{\alpha}$, $\me Z^{\alpha}_{\overline{\mu}}$, $\theta_{\overline{\mu}}^{\alpha}$, $\overline{\psi}^{\alpha}_{\overline{\mu}}$ be defined by coefficients of the canonical factorizations (\ref{fakt1_coi}) -- (\ref{fakt3_coi}) in the same way as were defined
  operators and vectors $\me C_{\mu,g}$, $\widetilde{\Phi}$,  $\Phi$, $\me Z_{\overline{\mu}}$, $\theta_{\overline{\mu}}$, $\overline{\psi}_{\overline{\mu}}$.
From Theorem $\ref{thm3_e_st.n_d}$ we obtain that the spectral characteristic
$\vec h^{\alpha}_{\overline{\mu}}(\lambda)$ of the optimal estimate
 $\widehat{A}\vec \xi$
of the functional $A\vec \xi$ can be calculated by the formula
\be\label{simple_spectr A_co_e_st.n_d}
\vec h^{\alpha}_{\overline{\mu}}(\lambda)=\frac{\chi_{\overline{\mu}}^{(d)}(e^{-i\lambda})}{\beta^{(d)}(i\lambda)}
\ld(\sum_{k=0}^{\infty}(\psi^{\alpha}_{\overline{\mu}}(k))^{\top} e^{-i\lambda k}\rd)
\sum_{m=1}^{\infty}\ld( (\theta_{\overline{\mu}}^{\alpha})^{\top}\widetilde{\me B}_{\overline{\mu}} -\overline{\psi}^{\alpha}_{\overline{\mu}} \me C^{\alpha}_{\mu,g}\rd)_m e^{-i\lambda m}.
\ee
The value of the mean square error of the estimate $\widehat{A}\vec \xi$ is calculated by the formula
\begin{eqnarray}
 \notag \Delta(f,g;\widehat{A}\vec\xi)&=&
  \|(\Phi^{\alpha})^{\top}\me a_{\overline{\mu}}\|^2+\|\widetilde{\Phi}^{\alpha}\me a_{\overline{\mu}}\|_1^2
 +\ld\langle (\theta^{\alpha}_{\overline{\mu}})^{\top}D^{\mu}\me A-\overline{\psi}^{\alpha}_{\overline{\mu}} \me C^{\alpha}_{\mu,g},(\theta_{\overline{\mu}}^{\alpha})^{\top}D^{\mu}\me A\rd\rangle
 \\ && -\ld\langle(\theta_{\overline{\mu}}^{\alpha})^{\top}D^{\mu}\me A,\me Z^{\alpha}_{\overline{\mu}}\me a_{\overline{\mu}}\rd\rangle
 -\ld\langle \me Z^{\alpha}_{\overline{\mu}}\me a_{\overline{\mu}},\overline{\psi}^{\alpha}_{\overline{\mu}} \me C^{\alpha}_{\mu,g}\rd\rangle_1.\label{simple_poh A_co_e_st.n_d}
\end{eqnarray}

\begin{thm}\label{thm6_e_st.n_d}
Let $\vec \xi(m)$, $m\in\mr Z\}$, and $\vec \zeta(m)$, $m\in\mr Z$, be seasonally cointegrated stochastic vector sequences with the spectral densities
 $f(\lambda)$ and $p(\lambda)$
 which admit the canonical factorizations (\ref{fakt1_coi}) -- (\ref{fakt3_coi}). Suppose that the coefficients $\vec a(k)$, $k\geq0$,
satisfy conditions
(\ref{umovana a_e_d}) -- (\ref{umovana a_mu_e_d}). Then the spectral characteristic $h^{\beta}_{\overline{\mu}}(\lambda)$ and the value of the mean square error $\Delta(f,g;\widehat{A}\vec \xi)$ of the optimal estimate $\widehat{A}\vec \xi$ of the functional $A\vec \xi$ based on observations of the sequence
$\vec \zeta(m)$ at points $m=-1,-2,\ldots$ are calculated by formulas
 (\ref{simple_spectr A_co_e_st.n_d}) and (\ref{simple_poh A_co_e_st.n_d}) respectively.
\end{thm}

\section{Minimax (robust) method of forecasting}\label{minimax_extrapolation}

Values of the mean square errors and the spectral characteristics of the optimal estimates
of the functionals ${A}\vec\xi$ and ${A}_N\vec\xi$
depending on the unobserved values of a stochastic sequence $\vec{\xi}(m)$ which determine a stationary stochastic GM increment sequence
$\chi_{\overline{\mu},\overline{s}}^{(d)}(\vec{\xi}(m))$ with the spectral density matrix $f(\lambda)$
based on observations of the sequence
$\vec\xi(m)+\vec\eta(m)$ at points $ m=-1,-2,\dots$ can be calculated by formulas
(\ref{spectr A_e_d}), (\ref{pohybkaA}) and
(\ref{spectr A_e_dN}), \eqref{pohybka33}
respectively, under the condition that
spectral densities
$f(\lambda)$ and $g(\lambda)$ of stochastic sequences $\vec\xi(m)$ and
$\vec\eta(m)$ are exactly known.

In practical cases, however, spectral densities of sequences usually are not exactly known.
If in such cases a set $\md D=\md D_f\times\md D_g$ of admissible spectral densities is defined,
the minimax (robust) approach to
estimation of linear functionals depending on unobserved values of stochastic sequences with stationary GM increments may be applied.
This method consists in finding an estimate that minimizes
the maximal values of the mean square errors for all spectral densities
from a given class $\md D=\md D_f\times\md D_g$ of admissible spectral densities
simultaneously.

To formalize this approach we present the following definitions.

\begin{ozn}
For a given class of spectral densities $\mathcal{D}=\md
D_f\times\md D_g$ the spectral densities
$f^0(\lambda)\in\mathcal{D}_f$, $g^0(\lambda)\in\md D_g$ are called
least favorable in the class $\mathcal{D}$ for the optimal linear
forecasting  of the functional $A\vec \xi$  if the following relation holds
true:
\[\Delta(f^0,g^0)=\Delta(h(f^0,g^0);f^0,g^0)=
\max_{(f,g)\in\mathcal{D}_f\times\md
D_g}\Delta(h(f,g);f,g).\]
\end{ozn}

\begin{ozn}
For a given class of spectral densities $\mathcal{D}=\md
D_f\times\md D_g$ the spectral characteristic $h^0(\lambda)$ of
the optimal linear estimate of the functional $A\vec \xi$ is called
minimax-robust if there are satisfied the conditions
\[h^0(\lambda)\in H_{\mathcal{D}}=\bigcap_{(f,g)\in\mathcal{D}_f\times\md D_g}L_2^{0-}(f(\lambda)+|\beta^{(d)}(i\lambda)|^2g(\lambda)),\]
\[\min_{h\in H_{\mathcal{D}}}\max_{(f,g)\in \mathcal{D}_f\times\md D_g}\Delta(h;f,g)=\max_{(f,g)\in\mathcal{D}_f\times\md
D_g}\Delta(h^0;f,g).\]
\end{ozn}

Taking into account the introduced definitions and the derived relations we can verify that the following lemmas hold true.

\begin{lema}
Spectral densities $f^0\in\mathcal{D}_f$,
$g^0\in\mathcal{D}_g$ which satisfy condition (\ref{umova11_e_st.n_d})
are least favorable in the class $\md D=\md D_f\times\md D_g$ for
the optimal linear forecasting of the functional $A\vec\xi$ if operators
$\me P_{\overline{\mu}}^0$, $\me T_{\overline{\mu}}^0$, $\me Q^0$  determined by the
Fourier coefficients of the functions
\[
\frac{|\beta^{(d)}(i\lambda)|^2}{|\chi_{\overline{\mu}}^{(d)}(e^{-i\lambda})|^2}
\ld[g^0(\lambda)(f^0(\lambda)+|\beta^{(d)}(i\lambda)|^2g^0(\lambda))^{-1}\rd]^{\top},
\]
\[
\dfrac{|\beta^{(d)}(i\lambda)|^2}{|\chi_{\overline{\mu}}^{(d)}(e^{-i\lambda})|^2}
\ld[f^0(\lambda)+|\beta^{(d)}(i\lambda)|^2g^0(\lambda))^{-1}\rd]^{\top},\]
\[
\ld[f^0(\lambda)(f^0(\lambda)+|\beta^{(d)}(i\lambda)|^2g^0(\lambda))^{-1}g^0(\lambda)\rd]^{\top}
\]
determine a solution of the constrained optimisation problem
\[
\max_{(f,g)\in \mathcal{D}_f\times\md D_g}(\langle D^{\overline{\mu}}\me a- \me T_{\overline{\mu}}\me
 a_{\overline{\mu}},\me P_{\overline{\mu}}^{-1}D^{\overline{\mu}}\me a-\me P_{\overline{\mu}}^{-1}\me T_{\overline{\mu}}\me a_{\overline{\mu}}\rangle+\langle\me Q\me a,\me a\rangle)
\]
\be
    = \langle D^{\overline{\mu}}\me a- \me T^0_{\overline{\mu}}\me
    a_{\overline{\mu}},(\me P^0_{\overline{\mu}})^{-1}D^{\overline{\mu}}\me a-(\me P^0_{\overline{\mu}})^{-1}\me T^0_{\overline{\mu}}\me a_{\overline{\mu}}\rangle+\langle\me Q^0\me a,\me a\rangle.
\label{minimax1}
\ee
The minimax spectral characteristic $h^0=h_{\overline{\mu}}(f^0,g^0)$ is calculated by formula (\ref{spectr A_e_d}) if
$h_{\overline{\mu}}(f^0,g^0)\in H_{\mathcal{D}}$.
\end{lema}

\begin{lema} The spectral densities $f^0\in\mathcal{D}_f$,
$g^0\in\mathcal{D}_g$ which admit  canonical factorizations (\ref{dd}), (\ref{fakt1}) and (\ref{fakt3})
are least favourable densities in the class $\mathcal{D}$ for the optimal linear forecasting
of the functional $A\vec \xi$ based on observations of the sequence $\vec \xi(m)+\vec \eta(m)$ at points $m=-1,-2,\ldots$ if the matrix coefficients
of  canonical factorizations
(\ref{fakt1}) and (\ref{fakt3})
determine a solution to the constrained optimization problem
\begin{eqnarray}
\notag
\|\Phi^{\top}\me a_{\overline{\mu}}\|^2&+&\|\widetilde{\Phi}\me a_{\overline{\mu}}\|_1^2
 +\ld\langle \theta^{\top}_{\overline{\mu}}D^{\mu}\me A-\overline{\psi}_{\overline{\mu}} \me C_{\mu,g},\theta^{\top}_{\overline{\mu}}D^{\mu}\me A\rd\rangle
 \\&-&  \ld\langle\theta^{\top}_{\overline{\mu}}D^{\mu}\me A,\me Z_{\overline{\mu}}\me a_{\overline{\mu}}\rd\rangle
 -\ld\langle \me Z_{\overline{\mu}}\me a_{\overline{\mu}},\overline{\psi}_{\overline{\mu}} \me C_{\mu,g}\rd\rangle_1\rightarrow\sup, \label{simple_minimax1_e_st.n_d}
\end{eqnarray}
\begin{eqnarray*}
 f(\lambda)&=&\frac{|\beta^{(d)}(i\lambda)|^2}{|\chi_{\overline{\mu}}^{(d)}(e^{-i\lambda})|^2}
\Theta_{\overline{\mu}}(e^{-i\lambda})\Theta_{\overline{\mu}}^*(e^{-i\lambda})
-|\beta^{(d)}(i\lambda)|^2\Phi(e^{-i\lambda})\Phi^*(e^{-i\lambda})\in \mathcal{D}_f,
\end{eqnarray*}
\[
 g(\lambda)=\Phi(e^{-i\lambda})\Phi^*(e^{-i\lambda})\in \mathcal{D}_g.\]
The minimax spectral characteristic $\vec h^0=\vec h_{\overline{\mu}}(f^0,g^0)$ is calculated by formula (\ref{simple_spectr A_e_st.n_d}) if
$\vec h_{\overline{\mu}}(f^0,g^0)\in H_{\mathcal{D}}$.
\end{lema}

\begin{lema} The spectral density $g^0\in\mathcal{D}_g$ which admits  canonical factorizations (\ref{fakt1}), (\ref{fakt3}) with the known spectral density $f(\lambda)$ is the least favourable in the class $\mathcal{D}_g$ for the optimal linear forecasting
of the functional $A\xi$ based on observations of the sequence $\vec \xi(m)+\vec \eta(m)$ at points $m=-1,-2,\ldots$ if the matrix coefficients
of the canonical factorizations
\be \label{fakt24_1_lf}
 f(\lambda)+|\beta^{(d)}(i\lambda)|^2g^0(\lambda)=\frac{|\beta^{(d)}(i\lambda)|^2}{|\chi_{\overline{\mu}}^{(d)}(e^{-i\lambda})|^2}
\ld(\sum_{k=0}^{\infty}\theta^0_{\overline{\mu}}(k)e^{-i\lambda k}\rd)\ld(\sum_{k=0}^{\infty}\theta^0_{\overline{\mu}}(k)e^{-i\lambda k}\rd)^*,\ee
\be \label{fakt24_2_lf}
 g^0(\lambda)=\ld(\sum_{k=0}^{\infty}\phi^0(k)e^{-i\lambda k}\rd)\ld(\sum_{k=0}^{\infty}\phi^0(k)e^{-i\lambda k}\rd)^*
\ee
determine a solution to the constrained optimization problem
\begin{eqnarray}
\notag
\|\Phi^{\top}\me a_{\overline{\mu}}\|^2&+&\|\widetilde{\Phi}\me a_{\overline{\mu}}\|_1^2
 +\ld\langle \theta^{\top}_{\overline{\mu}}D^{\mu}\me A-\overline{\psi}_{\overline{\mu}} \me C_{\mu,g},\theta^{\top}_{\overline{\mu}}D^{\mu}\me A\rd\rangle
 \\&-&  \ld\langle\theta^{\top}_{\overline{\mu}}D^{\mu}\me A,\me Z_{\overline{\mu}}\me a_{\overline{\mu}}\rd\rangle
 -\ld\langle \me Z_{\overline{\mu}}\me a_{\overline{\mu}},\overline{\psi}_{\overline{\mu}} \me C_{\mu,g}\rd\rangle_1\rightarrow\sup,
 \label{simple_minimax2_e_st.n_d}
 \end{eqnarray}
\[
 g(\lambda)=\Phi(e^{-i\lambda})\Phi^*(e^{-i\lambda})\in \mathcal{D}_g.\]
The minimax spectral characteristic $\vec h^0=\vec h_{\overline{\mu}}(f,g^0)$ is calculated by formula (\ref{simple_spectr A_e_st.n_d}) if
$\vec h_{\overline{\mu}}(f,g^0)\in H_{\mathcal{D}}$.
\end{lema}

\begin{lema} The spectral density $f^0\in\mathcal{D}_f$ which admits  canonical factorizations (\ref{dd}), (\ref{fakt1}) with the known spectral density $g(\lambda)$ is the least favourable spectral density in the class
 $\md D_f$ for the optimal linear forecasting
of the functional $A\vec \xi$ based on observations of the sequence $\vec \xi(m)+\vec \eta(m)$ at points $m=-1,-2,\ldots$ if matrix coefficients
of the canonical factorization
\be \label{fakt2_lf}
f^0(\lambda)+|\beta^{(d)}(i\lambda)|^2g(\lambda)=\frac{|\beta^{(d)}(i\lambda)|^2}{|\chi_{\overline{\mu}}^{(d)}(e^{-i\lambda})|^2}
\ld(\sum_{k=0}^{\infty}\theta^0_{\overline{\mu}}(k)e^{-i\lambda k}\rd)\ld(\sum_{k=0}^{\infty}\theta^0_{\overline{\mu}}(k)e^{-i\lambda k}\rd)^*,\ee
and the equastion $\Psi^0_{\overline{\mu}}(e^{-i\lambda})\Theta^0_{\overline{\mu}}(e^{-i\lambda})=E_q$ determine a solution to the constrained optimization problem
\begin{eqnarray}
\ld\langle \theta^{\top}_{\overline{\mu}}D^{\mu}\me A-\overline{\psi}_{\overline{\mu}} \me C_{\mu,g},\theta^{\top}_{\overline{\mu}}D^{\mu}\me A\rd\rangle
- \ld\langle\theta^{\top}_{\overline{\mu}}D^{\mu}\me A,\me Z_{\overline{\mu}}\me a_{\overline{\mu}}\rd\rangle
 -\ld\langle \me Z_{\overline{\mu}}\me a_{\overline{\mu}},\overline{\psi}_{\overline{\mu}} \me C_{\mu,g}\rd\rangle_1\rightarrow\sup,
 \label{simple_minimax3_e_st.n_d}\end{eqnarray}
\begin{eqnarray*}
 f(\lambda)&=&\frac{|\beta^{(d)}(i\lambda)|^2}{|\chi_{\overline{\mu}}^{(d)}(e^{-i\lambda})|^2}
\Theta_{\overline{\mu}}(e^{-i\lambda})\Theta_{\overline{\mu}}^*(e^{-i\lambda})
-|\beta^{(d)}(i\lambda)|^2\Phi(e^{-i\lambda})\Phi^*(e^{-i\lambda})\in \mathcal{D}_f\end{eqnarray*}
for the fixed matrix coefficients $\{\phi(k):k\geq0\}$. The minimax spectral characteristic $\vec h^0=\vec h_{\overline{\mu}}(f^0,g)$ is calculated by formula (\ref{simple_spectr A_e_st.n_d}) if
$\vec h_{\overline{\mu}}(f^0,g)\in H_{\mathcal{D}}$.
\end{lema}

For more detailed analysis of properties of the least favorable spectral densities and minimax-robust spectral characteristics we observe that
the minimax spectral characteristic $h^0$ and the least favourable spectral densities $(f^0,g^0)$
form a saddle
point of the function $\Delta(h;f,g)$ on the set
$H_{\mathcal{D}}\times\mathcal{D}$.

The saddle point inequalities
\[\Delta(h;f^0,g^0)\geq\Delta(h^0;f^0,g^0)\geq\Delta(h^0;f,g)
\quad\forall f\in \mathcal{D}_f,\forall g\in \mathcal{D}_g,\forall
h\in H_{\mathcal{D}}\] hold true if $h^0=h_{\overline{\mu}}(f^0,g^0)$ and
$h_{\overline{\mu}}(f^0,g^0)\in H_{\mathcal{D}}$, where $(f^0,g^0)$  is a
solution of the  constrained optimization problem
\be  \label{cond-extr1}
\widetilde{\Delta}(f,g)=-\Delta(h_{\overline{\mu}}(f^0,g^0);f,g)\to
\inf,\quad (f,g)\in \mathcal{D},\ee
where the functional $\Delta(h_{\overline{\mu}}(f^0,g^0);f,g)$ is calculated by the formula
\begin{multline*}
\Delta(h_{\overline{\mu}}(f^0,g^0);f,g)=
\\
=\frac{1}{2\pi}\int_{-\pi}^{\pi}
\frac{|\beta^{(d)}(i\lambda)|^2}{|\chi_{\overline{\mu}}^{(d)}(e^{-i\lambda})|^2}
\left[(\vec{A}_{\overline{\mu}}(e^{i\lambda}))^{\top}g^0(\lambda) +
(\vec{C}^0_{\overline{\mu}}(e^{i \lambda}))^{\top}
\right]
\times
\\
\times
(f^0(\lambda)+|\beta^{(d)}(i\lambda)|^2g^0(\lambda))^{-1}\, f(\lambda)\, (f^0(\lambda)+|\beta^{(d)}(i\lambda)|^2g^0(\lambda))^{-1}
\times
\\
\times
\left[g^0(\lambda)\overline{\vec{A}_{\overline{\mu}}(e^{i\lambda})} +
\overline{\vec{C}^0_{\overline{\mu}}(e^{i \lambda})}
\right]
d\lambda+
\\
+\frac{1}{2\pi}\int_{-\pi}^{\pi}
\left(\frac{|\beta^{(d)}(i\lambda)|^2}{|\chi_{\overline{\mu}}^{(d)}(e^{-i\lambda})|^2}\right)^2
\left[\frac{\chi_{\overline{\mu}}^{(d)}(e^{-i\lambda})}{|\beta^{(d)}(i\lambda)|^2}
(\vec{A}_{\overline{\mu}}(e^{i\lambda}))^{\top}f^0(\lambda) -\chi_{\overline{\mu}}^{(d)}(e^{-i\lambda})
(\vec{C}^0_{\overline{\mu}}(e^{i \lambda}))^{\top}
\right]
\times
\\
\times
(f^0(\lambda)+|\beta^{(d)}(i\lambda)|^2g^0(\lambda))^{-1}\, g(\lambda)\, (f^0(\lambda)+|\beta^{(d)}(i\lambda)|^2g^0(\lambda))^{-1}
\times
\\
\times
\left[\frac{\overline{\chi_{\overline{\mu}}^{(d)}(e^{-i\lambda})}}{|\beta^{(d)}(i\lambda)|^2}
f^0(\lambda)\overline{\vec{A}_{\overline{\mu}}(e^{i\lambda})} -\overline{\chi_{\overline{\mu}}^{(d)}(e^{-i\lambda})}
\overline{\vec{C}^0_{\overline{\mu}}(e^{i \lambda})}
\right]
d\lambda,
\end{multline*}
where
\[
\vec{C}^0_{\overline{\mu}}(e^{i \lambda})=\sum_{k=0}^{\infty}(\me P^0_{\overline{\mu}})^{-1}D^{\overline{\mu}}\me a-(\me P^0_{\overline{\mu}})^{-1}\me T^0_{\overline{\mu}}\me a_{\overline{\mu}})_k e^{ik\lambda},
\]
or it is calculated by the formula
\begin{multline*}
\Delta\ld(h_{\mu}(f^0,g^0);f,g\rd)=
\\
=\frac{1}{2\pi}\ip\frac{|\chi_{\overline{\mu}}^{(d)}(e^{-i\lambda})|^2}{|\beta^{(d)}(i\lambda)|^2}
(\me r^0_{\overline{\mu},f}(e^{-i\lambda}))^{\top}\Psi^0_{\overline{\mu}}(e^{-i\lambda }))f(\lambda)
(\Psi^0_{\overline{\mu}}(e^{-i\lambda }))^*\overline{\me r^0_{\overline{\mu},f}(e^{-i\lambda})}d\lambda
+
\\ +
\frac{1}{2\pi}\ip
(\me r^0_{\overline{\mu},g}(e^{-i\lambda}))^{\top}\Psi^0_{\overline{\mu}}(e^{-i\lambda })g(\lambda)
(\Psi^0_{\overline{\mu}}(e^{-i\lambda }))^*\overline{\me r^0_{\overline{\mu},g}(e^{-i\lambda})}d\lambda,
\end{multline*}
where
\[
 \me r^0_{\overline{\mu},f}(e^{-i\lambda})=\sum_{m=0}^{\infty}( (\theta^0)^{\top}_{\overline{\mu}}D^{\mu}\me A)_m e^{i\lambda m}+\sum_{m=1}^{\infty}(\overline{\psi}^0_{\overline{\mu}} \me C^0_{\mu,g} )_m e^{-i\lambda m},
\]
\[
 \me r^0_{\overline{\mu},g}(e^{-i\lambda})=\chi_{\overline{\mu}}^{(d)}(e^{-i\lambda})\ld(\sum_{m=0}^{\infty}( (\theta^0)^{\top}_{\overline{\mu}}D^{\mu}\me A)_m e^{i\lambda m}+\sum_{m=1}^{\infty}(\overline{\psi}_{\overline{\mu}}^0 \me C^0_{\mu,g} )_m e^{-i\lambda m}\rd)-(\Theta^0_{\overline{\mu}}(e^{-i\lambda }))^{\top}A(e^{i\lambda}).
\]

The constrained optimization problem (\ref{cond-extr1}) is equivalent to the unconstrained optimisation problem
\be  \label{uncond-extr}
\Delta_{\mathcal{D}}(f,g)=\widetilde{\Delta}(f,g)+ \delta(f,g|\mathcal{D}_f\times
\mathcal{D}_g)\to\inf,\ee
 where $\delta(f,g|\mathcal{D}_f\times
\mathcal{D}_g)$ is the indicator function of the set
$\mathcal{D}=\mathcal{D}_f\times\mathcal{D}_g$.
 Solution $(f^0,g^0)$ to this unconstrained optimization problem is characterized by the condition $0\in
\partial\Delta_{\mathcal{D}}(f^0,g^0)$, where
$\partial\Delta_{\mathcal{D}}(f^0,g^0)$ is the subdifferential of the functional $\Delta_{\mathcal{D}}(f,g)$ at point $(f^0,g^0)\in \mathcal{D}=\mathcal{D}_f\times\mathcal{D}_g$, that is the set of all continuous linear functionals $\Lambda$ on $L_1\times L_1$ which satisfy the inequality
$\Delta_{\mathcal{D}}(f,g)-\Delta_{\mathcal{D}}(f^0,g^0)\geq \Lambda ((f,g)-(f^0,g^0)), (f,g)\in \mathcal{D}$ (see  \cite{Moklyachuk2015,Rockafellar} for more details).
 This condition makes it possible to find the least favourable spectral densities in some special classes of spectral densities $\mathcal{D}=\mathcal{D}_f\times\mathcal{D}_g$.

In the case of cointegrated vector sequences (in terms of Subsection \ref{classical_extrap_coint}) we have the following optimization problem for determining the least favourable spectral densities:
\be\label{zad_bezum_extr_e_coi}
 \Delta_{\mathcal{D}}(f,p)
 \\=\widetilde{\Delta}(f,p)+\delta(f,p|\mathcal{D}_f\times
 \mathcal{D}_p)\to\inf,\ee
\[
 \widetilde{\Delta}(f,p)=-\Delta\ld(h^{\alpha}_{\overline{\mu}}(f^0,p^0);f,p\rd),\]
\[
\Delta\ld(h^{\alpha}_{\overline{\mu}}(f^0,p^0);f,p\rd)=
\]
\[
=
\frac{1}{2\pi}\int_{-\pi}^{\pi}
\left[\frac{1}{|\beta^{(d)}(i\lambda)|^2}
(\vec{A}_{\overline{\mu}}(e^{i\lambda}))^{\top}(p^0(\lambda)-\alpha^2f^0(\lambda)) +
\left(C^{\alpha,0}_{\overline{\mu}}(e^{i\lambda})
\right)^{\top}
\right]
\times
\]
\[
\times
\frac{|\beta^{(d)}(i\lambda)|^2}{|\chi_{\overline{\mu}}^{(d)}(e^{-i\lambda})|^2}(p^0(\lambda))^{-1}\, f(\lambda)\, (p^0(\lambda))^{-1}
\times
\]
\[
\times
\left[\frac{1}{|\beta^{(d)}(i\lambda)|^2}
(p^0(\lambda)-\alpha^2f^0(\lambda))\overline{\vec{A}_{\overline{\mu}}(e^{i\lambda})} +
\overline{C^{\alpha,0}_{\overline{\mu}}(e^{i\lambda})}
\right]
d\lambda+
\]
\[
+\frac{1}{2\pi}\int_{-\pi}^{\pi}
\left[\frac{1}{|\beta^{(d)}(i\lambda)|^2}(\vec{A}_{\overline{\mu}}(e^{i\lambda}))^{\top}f^0(\lambda) -
\left(C^{\alpha,0}_{\overline{\mu}}(e^{i\lambda})
\right)^{\top}
\right]
\times
\]
\[
\times
\frac{|\beta^{(d)}(i\lambda)|^2}{|\chi_{\overline{\mu}}^{(d)}(e^{-i\lambda})|^2}
(p^0(\lambda))^{-1}\,(p(\lambda)-\alpha^2f(\lambda))\,(p^0(\lambda))^{-1}
\times
\]
\[
\times
\left[\frac{1}{|\beta^{(d)}(i\lambda)|^2}f^0(\lambda)\overline{\vec{A}_{\overline{\mu}}(e^{i\lambda})} -
\frac{|\beta^{(d)}(i\lambda)|^2}{|\chi_{\overline{\mu}}^{(d)}(e^{-i\lambda})|^2}
\overline{C^{\alpha,0}_{\overline{\mu}}(e^{i\lambda})}
\right]
d\lambda,
\]
\[
C^{\alpha,0}_{\overline{\mu}}(e^{i\lambda})=\sum_{k=0}^{\infty}
 \ld(\ld((\me P_{\overline{\mu}}^{\alpha})^0\rd)^{-1}\ld(D^{\overline{\mu}}\me a
 -(\me T_{\overline{\mu}}^{\alpha})^0 \me a_{\overline{\mu}}\rd)\rd)_ke^{i\lambda k}.\]

A solution $(f^0,p^0)$ of this optimization problem is characterized by the condition
 $0\in\partial\Delta_{\mathcal{D}}(f^0,p^0)$.

The form of the functionals $\Delta(h_{\overline{\mu}}(f^0,g^0);f,g)$, $\Delta(h_{\overline{\mu}}(f^0,p^0);f,p)$ is convenient for application the Lagrange method of indefinite multipliers for
finding solution to the problem (\ref{uncond-extr}).
Making use of the method of Lagrange multipliers and the form of
subdifferentials of the indicator functions $\delta(f,g|\mathcal{D}_f\times
\mathcal{D}_g)$, $\delta(f,g|\mathcal{D}_f\times
\mathcal{D}_p)$ of the sets
$\mathcal{D}_f\times\mathcal{D}_g$, $\mathcal{D}_f\times\mathcal{D}_p$ of spectral densities,
we describe relations that determine least favourable spectral densities in some special classes
of spectral densities (see \cite{Luz_Mokl_book,Moklyachuk2015} for additional details).

\subsection{Least favorable spectral density in classes $\md D_0 \times \md  D_V^U$}\label{set1}

Consider the forecasting problem for the functional $A\vec{\xi}$
 which depends on unobserved values of a sequence $\vec\xi(m)$ with stationary increments based on observations of the sequence $\vec\xi(m)+\vec\eta(m)$ at points $m=-1,-2,\ldots$ under the condition that the sets of admissible spectral densities $\md D_{f0}^k, {\md D_{Vg}^{Uk}},k=1,2,3,4$ are defined as follows:
$$\md D_{f0}^{1} =\bigg\{f(\lambda )\left|\frac{1}{2\pi} \int
_{-\pi}^{\pi}
\frac{|\chi_{\overline{\mu}}^{(d)}(e^{-i\lambda})|^2}{|\beta^{(d)}(i\lambda)|^2}
f(\lambda )d\lambda  =P\right.\bigg\},$$
$$\md D_{f0}^{2} =\bigg\{f(\lambda )\left|\frac{1}{2\pi }
\int _{-\pi }^{\pi}
\frac{|\chi_{\overline{\mu}}^{(d)}(e^{-i\lambda})|^2}{|\beta^{(d)}(i\lambda)|^2}
{\rm{Tr}}\,[ f(\lambda )]d\lambda =p\right.\bigg\},$$
$$\md D_{f0}^{3} =\bigg\{f(\lambda )\left|\frac{1}{2\pi }
\int _{-\pi}^{\pi}
\frac{|\chi_{\overline{\mu}}^{(d)}(e^{-i\lambda})|^2}{|\beta^{(d)}(i\lambda)|^2}
f_{kk} (\lambda )d\lambda =p_{k}, k=\overline{1,T}\right.\bigg\},$$
$$\md D_{f0}^{4} =\bigg\{f(\lambda )\left|\frac{1}{2\pi} \int _{-\pi}^{\pi}
\frac{|\chi_{\overline{\mu}}^{(d)}(e^{-i\lambda})|^2}{|\beta^{(d)}(i\lambda)|^2}
\left\langle B_{1} ,f(\lambda )\right\rangle d\lambda  =p\right.\bigg\},$$

\begin{equation*}
 {\md D_{Vg}^{U1}}=\left\{g(\lambda )\bigg|V(\lambda )\le g(\lambda
)\le U(\lambda ), \frac{1}{2\pi } \int _{-\pi}^{\pi}
g(\lambda )d\lambda=Q\right\},
\end{equation*}
\begin{equation*}
  {\md D_{Vg}^{U}} ^{2}  =\bigg\{g(\lambda )\bigg|{\mathrm{Tr}}\, [V(\lambda
)]\le {\mathrm{Tr}}\,[ g(\lambda )]\le {\mathrm{Tr}}\, [U(\lambda )],
\frac{1}{2\pi } \int _{-\pi}^{\pi}
{\mathrm{Tr}}\,  [g(\lambda)]d\lambda  =q \bigg\},
\end{equation*}
\begin{equation*}
{\md D_{Vg}^{U3}}  =\bigg\{g(\lambda )\bigg|v_{kk} (\lambda )  \le
g_{kk} (\lambda )\le u_{kk} (\lambda ),
\frac{1}{2\pi} \int _{-\pi}^{\pi}
g_{kk} (\lambda
)d\lambda  =q_{k} , k=\overline{1,T}\bigg\},
\end{equation*}
\begin{equation*}
{\md D_{Vg}^{U4}}  =\bigg\{g(\lambda )\bigg|\left\langle B_{2}
,V(\lambda )\right\rangle \le \left\langle B_{2},g(\lambda
)\right\rangle \le \left\langle B_{2} ,U(\lambda)\right\rangle,
\frac{1}{2\pi }
\int _{-\pi}^{\pi}
\left\langle B_{2},g(\lambda)\right\rangle d\lambda  =q\bigg\}.
\end{equation*}

\noindent
Here spectral densities $V( \lambda ),U( \lambda )$ are known and fixed, $p, p_k, q, q_k, k=\overline{1,T}$ are given numbers, $P, B_1, Q, B_2$ are given positive-definite Hermitian matrices.

Define
\[
{\me C^{f0}_{\overline{\mu}}(e^{i\lambda})}
:=
g^0(\lambda)\vec{A}_{\overline{\mu}}(e^{i\lambda}) +
\sum_{k=0}^{\infty}((\me P^0_{\overline{\mu}})^{-1}D^{\overline{\mu}}\me a-(\me P^0_{\overline{\mu}})^{-1}\me T^0_{\overline{\mu}}\me a_{\overline{\mu}})_k e^{ik\lambda},
\]
\[
{\me C}^{g0}_{\overline{\mu}}(e^{i \lambda})
:=
\chi_{\overline{\mu}}^{(d)}(e^{-i\lambda})|\beta^{(d)}(i\lambda)|^{-2}
f^0(\lambda){\vec{A}_{\overline{\mu}}(e^{i\lambda})}
-
{\chi_{\overline{\mu}}^{(d)}(e^{-i\lambda})}
\sum_{k=0}^{\infty}((\me P^0_{\overline{\mu}})^{-1}D^{\overline{\mu}}\me a-(\me P^0_{\overline{\mu}})^{-1}\me T^0_{\overline{\mu}}\me a_{\overline{\mu}})_k e^{ik\lambda}.
\]

From the condition $0\in\partial\Delta_{\mathcal{D}}(f^0,g^0)$
we find the following equations which determine the least favourable spectral densities for these given sets of admissible spectral densities.

For the first set of admissible spectral densities $\md D_{f0}^1 \times\md D_{Vg}^{U1}$   we have equations
\begin{multline} \label{eq_4_1f}
\left(
{\me C^{f0}_{\overline{\mu}}(e^{i\lambda})}
\right)
\left(
{\me C^{f0}_{\overline{\mu}}(e^{i\lambda})}
\right)^{*}=
\left(\frac{|\chi_{\overline{\mu}}^{(d)}(e^{-i\lambda})|^2}{|\beta^{(d)}(i\lambda)|^2} (f^0(\lambda)+|\beta^{(d)}(i\lambda)|^2g^0(\lambda))\right)
\vec{\alpha}_f\cdot \vec{\alpha}_f^{*}
\times
\\
\times
\left(\frac{|\chi_{\overline{\mu}}^{(d)}(e^{-i\lambda})|^2}{|\beta^{(d)}(i\lambda)|^2} (f^0(\lambda)+|\beta^{(d)}(i\lambda)|^2g^0(\lambda))\right),
\end{multline}
\begin{multline} \label{eq_4_1g}
\left(
{\me C}^{g0}_{\overline{\mu}}(e^{i \lambda})
\right)
\left(
{\me C}^{g0}_{\overline{\mu}}(e^{i \lambda})
\right)^{*}=
\\
=
\left(\frac{|\chi_{\overline{\mu}}^{(d)}(e^{-i\lambda})|^2}{|\beta^{(d)}(i\lambda)|^2}(f^0(\lambda)+|\beta^{(d)}(i\lambda)|^2g^0(\lambda))\right)
(\vec{\beta}\cdot \vec{\beta}^{*}+\Gamma _{1} (\lambda )+\Gamma _{2} (\lambda ))
\times
\\
\times
\left(\frac{|\chi_{\overline{\mu}}^{(d)}(e^{-i\lambda})|^2}{|\beta^{(d)}(i\lambda)|^2}(f^0(\lambda)+|\beta^{(d)}(i\lambda)|^2g^0(\lambda))\right),
\end{multline}

\noindent where $\vec{\alpha}_f$ and $ \vec{\beta}$ are vectors of Lagrange multipliers, the matrix $\Gamma _{1} (\lambda )\le 0$ and $\Gamma _{1} (\lambda )=0$ if $g_{0}(\lambda )>V(\lambda ),$ the matrix  $
\Gamma _{2} (\lambda )\ge 0$ and $\Gamma _{2} (\lambda )=0$ if $g_{0}(\lambda )<U(\lambda ).$

For the second set of admissible spectral densities $\md D_{f0}^2 \times\md D_{Vg}^{U2}$  we have equations
\begin{equation} \label{eq_4_2f}
\left(
{\me C^{f0}_{\overline{\mu}}(e^{i\lambda})}
\right)
\left(
{\me C^{f0}_{\overline{\mu}}(e^{i\lambda})}
\right)^{*}=
\alpha_f^{2} \left(\frac{|\chi_{\overline{\mu}}^{(d)}(e^{-i\lambda})|^2}{|\beta^{(d)}(i\lambda)|^2} (f^0(\lambda)+|\beta^{(d)}(i\lambda)|^2g^0(\lambda))\right)^2,
\end{equation}
\begin{equation} \label{eq_4_2g}
\left(
{\me C}^{g0}_{\overline{\mu}}(e^{i \lambda})
\right)
\left(
{\me C}^{g0}_{\overline{\mu}}(e^{i \lambda})
\right)^{*}=
(\beta^{2} +\gamma _{1} (\lambda )+\gamma _{2} (\lambda )) \left(\frac{|\chi_{\overline{\mu}}^{(d)}(e^{-i\lambda})|^2}{|\beta^{(d)}(i\lambda)|^2}(f^0(\lambda)+|\beta^{(d)}(i\lambda)|^2g^0(\lambda))\right)^2,
\end{equation}

\noindent where  $\alpha _{f}^{2}$, $ \beta^{2}$ are Lagrange multipliers,  the function $\gamma _{1} (\lambda )\le 0$ and $\gamma _{1} (\lambda )=0$ if ${\mathrm{Tr}}\,
[g_{0} (\lambda )]> {\mathrm{Tr}}\,  [V(\lambda )],$ the function $\gamma _{2} (\lambda )\ge 0$ and $\gamma _{2} (\lambda )=0$ if $ {\mathrm{Tr}}\,[g_{0}(\lambda )]< {\mathrm{Tr}}\, [ U(\lambda)].$

For the third set of admissible spectral densities $\md D_{f0}^3 \times\md D_{Vg}^{U3}$  we have equations
\begin{multline} \label{eq_4_3f}
\left(
{\me C^{f0}_{\overline{\mu}}(e^{i\lambda})}
\right)
\left(
{\me C^{f0}_{\overline{\mu}}(e^{i\lambda})}
\right)^{*}=
\left(\frac{|\chi_{\overline{\mu}}^{(d)}(e^{-i\lambda})|^2}{|\beta^{(d)}(i\lambda)|^2} (f^0(\lambda)+|\beta^{(d)}(i\lambda)|^2g^0(\lambda))\right)
\left\{\alpha _{fk}^{2} \delta _{kl} \right\}_{k,l=1}^{T}
\times
\\
\times
\left(\frac{|\chi_{\overline{\mu}}^{(d)}(e^{-i\lambda})|^2}{|\beta^{(d)}(i\lambda)|^2} (f^0(\lambda)+|\beta^{(d)}(i\lambda)|^2g^0(\lambda))\right),
\end{multline}
\begin{multline} \label{eq_4_3g}
\left(
{\me C}^{g0}_{\overline{\mu}}(e^{i \lambda})
\right)
\left(
{\me C}^{g0}_{\overline{\mu}}(e^{i \lambda})
\right)^{*}=
\\
=\left(\frac{|\chi_{\overline{\mu}}^{(d)}(e^{-i\lambda})|^2}{|\beta^{(d)}(i\lambda)|^2}(f^0(\lambda)+|\beta^{(d)}(i\lambda)|^2g^0(\lambda))\right)
\left\{(\beta_{k}^{2} +\gamma _{1k} (\lambda )+\gamma _{2k} (\lambda ))\delta _{kl}\right\}_{k,l=1}^{T}
\times
\\
\times
\left(\frac{|\chi_{\overline{\mu}}^{(d)}(e^{-i\lambda})|^2}{|\beta^{(d)}(i\lambda)|^2}(f^0(\lambda)+|\beta^{(d)}(i\lambda)|^2g^0(\lambda))\right),
\end{multline}

\noindent where  $\alpha _{fk}^{2}$,   $\beta_{k}^{2}$ are Lagrange multipliers,
 $\delta _{kl}$ are Kronecker symbols, functions $\gamma _{1k} (\lambda )\le 0$ and $\gamma _{1k} (\lambda )=0$ if $g_{kk}^{0} (\lambda )>v_{kk} (\lambda ),$ functions $\gamma _{2k} (\lambda )\ge 0$ and $\gamma _{2k} (\lambda )=0$ if $g_{kk}^{0} (\lambda )<u_{kk} (\lambda).$

For the fourth set of admissible spectral densities $\md D_{f0}^4 \times\md D_{Vg}^{U4}$  we have equations

\begin{multline} \label{eq_4_4f}
\left(
{\me C^{f0}_{\overline{\mu}}(e^{i\lambda})}
\right)
\left(
{\me C^{f0}_{\overline{\mu}}(e^{i\lambda})}
\right)^{*}=
\alpha_f^{2} \left(\frac{|\chi_{\overline{\mu}}^{(d)}(e^{-i\lambda})|^2}{|\beta^{(d)}(i\lambda)|^2} (f^0(\lambda)+|\beta^{(d)}(i\lambda)|^2g^0(\lambda))\right)
B_{1}^{\top}
\times
\\
\times
\left(\frac{|\chi_{\overline{\mu}}^{(d)}(e^{-i\lambda})|^2}{|\beta^{(d)}(i\lambda)|^2} (f^0(\lambda)+|\beta^{(d)}(i\lambda)|^2g^0(\lambda))\right),
\end{multline}
\begin{multline} \label{eq_4_4g}
\left(
{\me C}^{g0}_{\overline{\mu}}(e^{i \lambda})
\right)
\left(
{\me C}^{g0}_{\overline{\mu}}(e^{i \lambda})
\right)^{*}=
(\beta^{2} +\gamma'_{1}(\lambda )+\gamma'_{2}(\lambda ))
\left(\frac{|\chi_{\overline{\mu}}^{(d)}(e^{-i\lambda})|^2}{|\beta^{(d)}(i\lambda)|^2}(f^0(\lambda)+|\beta^{(d)}(i\lambda)|^2g^0(\lambda))\right)
B_{2}^{\top}
\times
\\
\times
\left(\frac{|\chi_{\overline{\mu}}^{(d)}(e^{-i\lambda})|^2}{|\beta^{(d)}(i\lambda)|^2}(f^0(\lambda)+|\beta^{(d)}(i\lambda)|^2g^0(\lambda))\right),
\end{multline}

\noindent where $ \beta^{2}$, $\alpha _{f}^{2}$,  are Lagrange multipliers, functions $\gamma'_{1}( \lambda )\le 0$ and $\gamma'_{1} ( \lambda )=0$ if $\langle B_{2},g_{0} ( \lambda) \rangle > \langle B_{2},V( \lambda ) \rangle,$ functions $\gamma'_{2}( \lambda )\ge 0$ and $\gamma'_{2} ( \lambda )=0$ if $\langle
B_{2} ,g_{0} ( \lambda) \rangle < \langle B_{2} ,U( \lambda ) \rangle.$

The following theorem  holds true.

\begin{thm}
Let the minimality condition (\ref{umova11_e_st.n_d}) hold true. The least favorable spectral densities $f_{0}(\lambda), $ $g_{0}(\lambda), $ in the classes $\md D_{f0}^k \times\md D_{Vg}^{Uk},k=1,2,3,4$  for the optimal linear extrapolation of the functional  $A\vec{\xi}$ from observations of the sequence $\vec{\xi}(m)+ \vec{\eta}(m)$ at points  $m=-1,-2,\ldots$  are determined by equations
\eqref{eq_4_1f}--\eqref{eq_4_1g},  \eqref{eq_4_2f}--\eqref{eq_4_2g}, \eqref{eq_4_3f}--\eqref{eq_4_3g}, \eqref{eq_4_4f}--\eqref{eq_4_4g},
respectively,
the constrained optimization problem (\ref{minimax1}) and restrictions  on densities from the corresponding classes $\md D_{f0}^k, \md D_{Vg}^{Uk},k=1,2,3,4$.  The minimax-robust spectral characteristic of the optimal estimate of the functional $A\vec{\xi}$ is determined by the formula (\ref{spectr A_e_d}).
\end{thm}

If the spectral densities $f(\lambda)$ and $g(\lambda)$ admit canonical factorizations (\ref{dd}), (\ref{fakt1}) and (\ref{fakt3}), we can derive the  following equation for the least favourable spectral densities.

For the first set of admissible spectral densities $\md D_{f0}^1 \times\md D_{Vg}^{U1}$:
\begin{equation} \label{eq_4_1f_fact}
\left(
{\me r^{0}_{\overline{\mu},f}(e^{i\lambda})}
\right)
\left(
{\me r^{0}_{\overline{\mu},f}(e^{i\lambda})}
\right)^{*}
=(\Theta_{\overline{\mu}}(e^{-i\lambda}))^{\top}
\vec{\alpha}_f\cdot \vec{\alpha}_f^{*}
\overline{\Theta_{\overline{\mu}}(e^{-i\lambda})},
\end{equation}
\begin{equation} \label{eq_4_1g_fact}
\left(
{\me r^{0}_{\overline{\mu},g}(e^{i\lambda})}
\right)
\left(
{\me r^{0}_{\overline{\mu},g}(e^{i\lambda})}
\right)^{*}
=
(\Theta_{\overline{\mu}}(e^{-i\lambda}))^{\top}
(\vec{\beta}\cdot \vec{\beta}^{*}+\Gamma _{1} (\lambda )+\Gamma _{2} (\lambda ))
\overline{\Theta_{\overline{\mu}}(e^{-i\lambda})}
\end{equation}

\noindent where $\vec{\alpha}_f$ and $ \vec{\beta}$ are vectors of Lagrange multipliers, the matrix $\Gamma _{1} (\lambda )\le 0$ and $\Gamma _{1} (\lambda )=0$ if $g_{0}(\lambda )>V(\lambda ),$ the matrix  $
\Gamma _{2} (\lambda )\ge 0$ and $\Gamma _{2} (\lambda )=0$ if $g_{0}(\lambda )<U(\lambda ).$

For the second set of admissible spectral densities $\md D_{f0}^2 \times\md D_{Vg}^{U2}$  we have equations
\begin{equation} \label{eq_4_2f_fact}
\left(
{\me r^{0}_{\overline{\mu},f}(e^{i\lambda})}
\right)
\left(
{\me r^{0}_{\overline{\mu},f}(e^{i\lambda})}
\right)^{*}=
\alpha_f^{2} (\Theta_{\overline{\mu}}(e^{-i\lambda}))^{\top}\overline{\Theta_{\overline{\mu}}(e^{-i\lambda})},
\end{equation}
\begin{equation} \label{eq_4_2g_fact}
\left(
{\me r^{0}_{\overline{\mu},g}(e^{i\lambda})}
\right)
\left(
{\me r^{0}_{\overline{\mu},g}(e^{i\lambda})}
\right)^{*}
=(\beta^{2} +\gamma _{1} (\lambda )+\gamma _{2} (\lambda )) (\Theta_{\overline{\mu}}(e^{-i\lambda}))^{\top}\overline{\Theta_{\overline{\mu}}(e^{-i\lambda})},
\end{equation}

\noindent where  $\alpha _{f}^{2}$, $ \beta^{2}$ are Lagrange multipliers,  the function $\gamma _{1} (\lambda )\le 0$ and $\gamma _{1} (\lambda )=0$ if ${\mathrm{Tr}}\,
[g_{0} (\lambda )]> {\mathrm{Tr}}\,  [V(\lambda )],$ the function $\gamma _{2} (\lambda )\ge 0$ and $\gamma _{2} (\lambda )=0$ if $ {\mathrm{Tr}}\,[g_{0}(\lambda )]< {\mathrm{Tr}}\, [ U(\lambda)].$

For the third set of admissible spectral densities $\md D_{f0}^3 \times\md D_{Vg}^{U3}$  we have equations
\begin{equation} \label{eq_4_3f_fact}
\left(
{\me r^{0}_{\overline{\mu},f}(e^{i\lambda})}
\right)
\left(
{\me r^{0}_{\overline{\mu},f}(e^{i\lambda})}
\right)^{*}
=(\Theta_{\overline{\mu}}(e^{-i\lambda}))^{\top}
\left\{\alpha _{fk}^{2} \delta _{kl} \right\}_{k,l=1}^{T}
\overline{\Theta_{\overline{\mu}}(e^{-i\lambda})},
\end{equation}
\begin{equation} \label{eq_4_3g_fact}
\left(
{\me r^{0}_{\overline{\mu},g}(e^{i\lambda})}
\right)
\left(
{\me r^{0}_{\overline{\mu},g}(e^{i\lambda})}
\right)^{*}
=(\Theta_{\overline{\mu}}(e^{-i\lambda}))^{\top}
\left\{(\beta_{k}^{2} +\gamma _{1k} (\lambda )+\gamma _{2k} (\lambda ))\delta _{kl}\right\}_{k,l=1}^{T}
\overline{\Theta_{\overline{\mu}}(e^{-i\lambda})},
\end{equation}

\noindent where  $\alpha _{fk}^{2}$,   $\beta_{k}^{2}$ are Lagrange multipliers,
 $\delta _{kl}$ are Kronecker symbols, functions $\gamma _{1k} (\lambda )\le 0$ and $\gamma _{1k} (\lambda )=0$ if $g_{kk}^{0} (\lambda )>v_{kk} (\lambda ),$ functions $\gamma _{2k} (\lambda )\ge 0$ and $\gamma _{2k} (\lambda )=0$ if $g_{kk}^{0} (\lambda )<u_{kk} (\lambda).$

For the fourth set of admissible spectral densities $\md D_{f0}^4 \times\md D_{Vg}^{U4}$  we have equations

\begin{equation} \label{eq_4_4f_fact}
\left(
{\me r^{0}_{\overline{\mu},f}(e^{i\lambda})}
\right)
\left(
{\me r^{0}_{\overline{\mu},f}(e^{i\lambda})}
\right)^{*}
=
\alpha_f^{2} (\Theta_{\overline{\mu}}(e^{-i\lambda}))^{\top}
B_{1}
\overline{\Theta_{\overline{\mu}}(e^{-i\lambda})},
\end{equation}
\begin{equation} \label{eq_4_4g_fact}
\left(
{\me r^{0}_{\overline{\mu},g}(e^{i\lambda})}
\right)
\left(
{\me r^{0}_{\overline{\mu},g}(e^{i\lambda})}
\right)^{*}
=
(\beta^{2} +\gamma'_{1}(\lambda )+\gamma'_{2}(\lambda ))
(\Theta_{\overline{\mu}}(e^{-i\lambda}))^{\top}
B_{2}
\overline{\Theta_{\overline{\mu}}(e^{-i\lambda})},
\end{equation}

\noindent where $ \beta^{2}$, $\alpha _{f}^{2}$,  are Lagrange multipliers, functions $\gamma'_{1}( \lambda )\le 0$ and $\gamma'_{1} ( \lambda )=0$ if $\langle B_{2},g_{0} ( \lambda) \rangle > \langle B_{2},V( \lambda ) \rangle,$ functions $\gamma'_{2}( \lambda )\ge 0$ and $\gamma'_{2} ( \lambda )=0$ if $\langle
B_{2} ,g_{0} ( \lambda) \rangle < \langle B_{2} ,U( \lambda ) \rangle.$

The following theorem  holds true.

\begin{thm}
 The least favorable spectral densities $f_{0}(\lambda)$,  $g_{0}(\lambda)$  in the classes $\md D_{f0}^k \times\md D_{Vg}^{Uk}$, $k=1,2,3,4$ for the optimal linear forecasting of the functional  $A\vec{\xi}$ from observations of the sequence $\vec{\xi}(m)+ \vec{\eta}(m)$ at points  $m=-1,-2,\ldots$  are determined
by  canonical factorizations (\ref{dd}), (\ref{fakt1}) and (\ref{fakt3}),
 equations
\eqref{eq_4_1f_fact}--\eqref{eq_4_1g_fact},  \eqref{eq_4_2f_fact}--\eqref{eq_4_2g_fact}, \eqref{eq_4_3f_fact}--\eqref{eq_4_3g_fact}, \eqref{eq_4_4f_fact}--\eqref{eq_4_4g_fact},
respectively,
 constrained optimization problem (\ref{simple_minimax1_e_st.n_d}) and restrictions  on densities from the corresponding classes $\md D_{f0}^k, \md D_{Vg}^{Uk},k=1,2,3,4$.  The minimax-robust spectral characteristic of the optimal estimate of the functional $A\vec{\xi}$ is determined by the formula (\ref{spectr A_e_d}).
\end{thm}

\begin{thm}
 If the spectral density $g(\lambda)$ is known, the least favorable spectral density $f_{0}(\lambda)$ in the classes $\md D_{f0}^k $, $k=1,2,3,4$    for the optimal linear forecasting of the functional  $A\vec{\xi}$ from observations of the sequence $\vec{\xi}(m)+ \vec{\eta}(m)$ at points  $m=-1,-2,\ldots$  is determined
by  canonical factorizations (\ref{dd}) and (\ref{fakt1})  ,
 equations
\eqref{eq_4_1f_fact},  \eqref{eq_4_2f_fact}, \eqref{eq_4_3f_fact}, \eqref{eq_4_4f_fact},
respectively,
 constrained optimization problem (\ref{simple_minimax3_e_st.n_d}) and restrictions  on density from the corresponding classes $\md D_{f0}^k$, $k=1,2,3,4$.  The minimax-robust spectral characteristic of the optimal estimate of the functional $A\vec{\xi}$ is determined by the formula (\ref{spectr A_e_d}).
\end{thm}

\subsection{Least favorable spectral density in classes $\md D_0 \times \md  D_V^U$ for cointegrated vector sequences}\label{set1_coint}

Consider the forecasting problem for the functional $A\vec{\xi}$
 which depends on unobserved values of a sequence $\vec\xi(m)$ with stationary increments based on observations of the vector sequence $\vec\zeta(m)$, cointegrated with $\vec\xi(m)$ in terms of Definition \ref{ozn_coint},  at points $m=-1,-2,\ldots$ for the sets of admissible spectral densities $\md D_{f0}^k, {\md D_{Vp}^{Uk}},k=1,2,3,4$, where the sets $\md D_{f0}^k$ are defined in Subsection \ref{set1}, the sets $\md D_{Vp}^{Uk}$ are the following:

\begin{equation*}
 {\md D_{Vp}^{U1}}=\left\{p(\lambda )\bigg|V(\lambda )\le p(\lambda
)\le U(\lambda ), \frac{1}{2\pi } \int _{-\pi}^{\pi}
\frac{|\chi_{\overline{\mu}}^{(d)}(e^{-i\lambda})|^2}{|\beta^{(d)}(i\lambda)|^2}p(\lambda )d\lambda=Q\right\},
\end{equation*}
\begin{equation*}
  {\md D_{Vp}^{U}} ^{2}  =\bigg\{p(\lambda )\bigg|{\mathrm{Tr}}\, [V(\lambda
)]\le {\mathrm{Tr}}\,[ p(\lambda )]\le {\mathrm{Tr}}\, [U(\lambda )],
\frac{1}{2\pi } \int _{-\pi}^{\pi}\frac{|\chi_{\overline{\mu}}^{(d)}(e^{-i\lambda})|^2}{|\beta^{(d)}(i\lambda)|^2}
{\mathrm{Tr}}\,  [p(\lambda)]d\lambda  =q \bigg\},
\end{equation*}
\begin{equation*}
{\md D_{Vp}^{U3}}  =\bigg\{p(\lambda )\bigg|v_{kk} (\lambda )  \le
p_{kk} (\lambda )\le u_{kk} (\lambda ),
\frac{1}{2\pi} \int _{-\pi}^{\pi}\frac{|\chi_{\overline{\mu}}^{(d)}(e^{-i\lambda})|^2}{|\beta^{(d)}(i\lambda)|^2}
p_{kk} (\lambda
)d\lambda  =q_{k} , k=\overline{1,T}\bigg\},
\end{equation*}
\begin{equation*}
{\md D_{Vp}^{U4}}  =\bigg\{p(\lambda )\bigg|\left\langle B_{2}
,V(\lambda )\right\rangle \le \left\langle B_{2},p(\lambda
)\right\rangle \le \left\langle B_{2} ,U(\lambda)\right\rangle,
\frac{1}{2\pi }
\int _{-\pi}^{\pi}\frac{|\chi_{\overline{\mu}}^{(d)}(e^{-i\lambda})|^2}{|\beta^{(d)}(i\lambda)|^2}
\left\langle B_{2},p(\lambda)\right\rangle d\lambda  =q\bigg\}.
\end{equation*}

\noindent
Here spectral densities $V( \lambda ),U( \lambda )$ are known and fixed, $q, q_k, k=\overline{1,T}$ are given numbers, $Q, B_2$ are given positive-definite Hermitian matrices.

Define
\begin{equation*}
{\me C^{f0}_{\overline{\mu},\alpha}(e^{i\lambda})}
:=
\frac{\overline{\chi_{\overline{\mu}}^{(d)}(e^{-i\lambda})}}{|\beta^{(d)}(i\lambda)|^2}
\vec{A}(e^{i\lambda})(p^0(\lambda)-\alpha^2f^0(\lambda))
+
\sum_{k=0}^{\infty}
 \ld(\ld((\me P_{\overline{\mu}}^{\alpha})^0\rd)^{-1}\ld(D^{\overline{\mu}}\me a
 -(\me T_{\overline{\mu}}^{\alpha})^0 \me a_{\overline{\mu}}\rd)\rd)_ke^{i\lambda k},
 \end{equation*}
\[
{\me C}^{p0}_{\overline{\mu},\alpha}(e^{i \lambda})
:=
\frac{\overline{\chi_{\overline{\mu}}^{(d)}(e^{-i\lambda})}}{|\beta^{(d)}(i\lambda)|^2}
\vec{A}(e^{i\lambda})f^0(\lambda) -
\sum_{k=0}^{\infty}
 \ld(\ld((\me P_{\overline{\mu}}^{\alpha})^0\rd)^{-1}\ld(D^{\overline{\mu}}\me a
 -(\me T_{\overline{\mu}}^{\alpha})^0 \me a_{\overline{\mu}}\rd)\rd)_ke^{i\lambda k}.
\]

The condition $0\in\partial\Delta_{\mathcal{D}}(f^0,p^0)$
implies the following equations which determine the least favourable spectral densities for these given sets of admissible spectral densities.

For the first set of admissible spectral densities $\md D_{f0}^1 \times\md D_{Vp}^{U1}$   we have equations
\begin{multline} \label{eq_4_1f_coint}
\left(
{\me C^{f0}_{\overline{\mu},\alpha}(e^{i\lambda})}
\right)
\left(
{\me C^{f0}_{\overline{\mu},\alpha}(e^{i\lambda})}
\right)^{*}
-\alpha^2\left(
{\me C}^{p0}_{\overline{\mu},\alpha}(e^{i \lambda})
\right)
\left(
{\me C}^{p0}_{\overline{\mu},\alpha}(e^{i \lambda})
\right)^{*}=
\\
=\left(\frac{|\chi_{\overline{\mu}}^{(d)}(e^{-i\lambda})|^2}{|\beta^{(d)}(i\lambda)|^2} p^0(\lambda)\right)
\vec{\alpha}_f\cdot \vec{\alpha}_f^{*}
\left(\frac{|\chi_{\overline{\mu}}^{(d)}(e^{-i\lambda})|^2}{|\beta^{(d)}(i\lambda)|^2} p^0(\lambda)\right),
\end{multline}
\begin{equation}\label{eq_4_1g_coint}
\left(
{\me C}^{p0}_{\overline{\mu},\alpha}(e^{i \lambda})
\right)
\left(
{\me C}^{p0}_{\overline{\mu},\alpha}(e^{i \lambda})
\right)^{*}=
\left(\frac{|\chi_{\overline{\mu}}^{(d)}(e^{-i\lambda})|^2}{|\beta^{(d)}(i\lambda)|^2} p^0(\lambda)\right)
(\vec{\beta}\cdot \vec{\beta}^{*}+\Gamma _{1} (\lambda )+\Gamma _{2} (\lambda ))
\left(\frac{|\chi_{\overline{\mu}}^{(d)}(e^{-i\lambda})|^2}{|\beta^{(d)}(i\lambda)|^2} p^0(\lambda)\right),
\end{equation}

\noindent where $\vec{\alpha}_f$ and $ \vec{\beta}$ are vectors of Lagrange multipliers, the matrix $\Gamma _{1} (\lambda )\le 0$ and $\Gamma _{1} (\lambda )=0$ if $g_{0}(\lambda )>V(\lambda ),$ the matrix  $
\Gamma _{2} (\lambda )\ge 0$ and $\Gamma _{2} (\lambda )=0$ if $p_{0}(\lambda )<U(\lambda ).$

For the second set of admissible spectral densities $\md D_{f0}^2 \times\md D_{Vp}^{U2}$  we have equations
\begin{equation} \label{eq_4_2f_coint}
\left(
{\me C^{f0}_{\overline{\mu},\alpha}(e^{i\lambda})}
\right)
\left(
{\me C^{f0}_{\overline{\mu},\alpha}(e^{i\lambda})}
\right)^{*}
-\alpha^2\left(
{\me C}^{p0}_{\overline{\mu},\alpha}(e^{i \lambda})
\right)
\left(
{\me C}^{p0}_{\overline{\mu},\alpha}(e^{i \lambda})
\right)^{*}
=
\alpha_f^{2} \left(\frac{|\chi_{\overline{\mu}}^{(d)}(e^{-i\lambda})|^2}{|\beta^{(d)}(i\lambda)|^2} p^0(\lambda)\right)^2,
\end{equation}
\begin{equation} \label{eq_4_2g_coint}
\left(
{\me C}^{p0}_{\overline{\mu},\alpha}(e^{i \lambda})
\right)
\left(
{\me C}^{p0}_{\overline{\mu},\alpha}(e^{i \lambda})
\right)^{*}
=(\beta^{2} +\gamma _{1} (\lambda )+\gamma _{2} (\lambda )) \left(\frac{|\chi_{\overline{\mu}}^{(d)}(e^{-i\lambda})|^2}{|\beta^{(d)}(i\lambda)|^2} p^0(\lambda)\right)^2,
\end{equation}

\noindent where  $\alpha _{f}^{2}$, $ \beta^{2}$ are Lagrange multipliers,  the function $\gamma _{1} (\lambda )\le 0$ and $\gamma _{1} (\lambda )=0$ if ${\mathrm{Tr}}\,
[p_{0} (\lambda )]> {\mathrm{Tr}}\,  [V(\lambda )],$ the function $\gamma _{2} (\lambda )\ge 0$ and $\gamma _{2} (\lambda )=0$ if $ {\mathrm{Tr}}\,[p_{0}(\lambda )]< {\mathrm{Tr}}\, [ U(\lambda)].$

For the third set of admissible spectral densities $\md D_{f0}^3 \times\md D_{Vp}^{U3}$  we have equations
\begin{multline} \label{eq_4_3f_coint}
\left(
{\me C^{f0}_{\overline{\mu},\alpha}(e^{i\lambda})}
\right)
\left(
{\me C^{f0}_{\overline{\mu},\alpha}(e^{i\lambda})}
\right)^{*}
-\alpha^2\left(
{\me C}^{p0}_{\overline{\mu},\alpha}(e^{i \lambda})
\right)
\left(
{\me C}^{p0}_{\overline{\mu},\alpha}(e^{i \lambda})
\right)^{*}=
\\
=\left(\frac{|\chi_{\overline{\mu}}^{(d)}(e^{-i\lambda})|^2}{|\beta^{(d)}(i\lambda)|^2} p^0(\lambda)\right)
\left\{\alpha _{fk}^{2} \delta _{kl} \right\}_{k,l=1}^{T}
\left(\frac{|\chi_{\overline{\mu}}^{(d)}(e^{-i\lambda})|^2}{|\beta^{(d)}(i\lambda)|^2} p^0(\lambda)\right),
\end{multline}
\begin{multline} \label{eq_4_3g_coint}
\left(
{\me C}^{p0}_{\overline{\mu},\alpha}(e^{i \lambda})
\right)
\left(
{\me C}^{p0}_{\overline{\mu},\alpha}(e^{i \lambda})
\right)^{*}=
\left(\frac{|\chi_{\overline{\mu}}^{(d)}(e^{-i\lambda})|^2}{|\beta^{(d)}(i\lambda)|^2} p^0(\lambda)\right)
\times
\\
\times
\left\{(\beta_{k}^{2} +\gamma _{1k} (\lambda )+\gamma _{2k} (\lambda ))\delta _{kl}\right\}_{k,l=1}^{T}
\left(\frac{|\chi_{\overline{\mu}}^{(d)}(e^{-i\lambda})|^2}{|\beta^{(d)}(i\lambda)|^2} p^0(\lambda)\right),
\end{multline}

\noindent where  $\alpha _{fk}^{2}$,   $\beta_{k}^{2}$ are Lagrange multipliers,
 $\delta _{kl}$ are Kronecker symbols, functions $\gamma _{1k} (\lambda )\le 0$ and $\gamma _{1k} (\lambda )=0$ if $p_{kk}^{0} (\lambda )>v_{kk} (\lambda ),$ functions $\gamma _{2k} (\lambda )\ge 0$ and $\gamma _{2k} (\lambda )=0$ if $p_{kk}^{0} (\lambda )<u_{kk} (\lambda).$

For the fourth set of admissible spectral densities $\md D_{f0}^4 \times\md D_{Vp}^{U4}$  we have equations
\begin{multline} \label{eq_4_4f_coint}
\left(
{\me C^{f0}_{\overline{\mu},\alpha}(e^{i\lambda})}
\right)
\left(
{\me C^{f0}_{\overline{\mu},\alpha}(e^{i\lambda})}
\right)^{*}
-\alpha^2\left(
{\me C}^{p0}_{\overline{\mu},\alpha}(e^{i \lambda})
\right)
\left(
{\me C}^{p0}_{\overline{\mu},\alpha}(e^{i \lambda})
\right)^{*}=
\\
=
\alpha_f^{2} \left(\frac{|\chi_{\overline{\mu}}^{(d)}(e^{-i\lambda})|^2}{|\beta^{(d)}(i\lambda)|^2} p^0(\lambda)\right)
B_{1}^{\top}
\left(\frac{|\chi_{\overline{\mu}}^{(d)}(e^{-i\lambda})|^2}{|\beta^{(d)}(i\lambda)|^2} p^0(\lambda)\right),
\end{multline}
\begin{multline} \label{eq_4_4g_coint}
\left(
{\me C}^{p0}_{\overline{\mu},\alpha}(e^{i \lambda})
\right)
\left(
{\me C}^{p0}_{\overline{\mu},\alpha}(e^{i \lambda})
\right)^{*}
=
(\beta^{2} +\gamma'_{1}(\lambda )+\gamma'_{2}(\lambda ))
\left(\frac{|\chi_{\overline{\mu}}^{(d)}(e^{-i\lambda})|^2}{|\beta^{(d)}(i\lambda)|^2} p^0(\lambda)\right)
B_{2}^{\top}
\times
\\
\times
\left(\frac{|\chi_{\overline{\mu}}^{(d)}(e^{-i\lambda})|^2}{|\beta^{(d)}(i\lambda)|^2} p^0(\lambda)\right),
\end{multline}

\noindent where $\alpha _{f}^{2}$, $\beta^{2}$,   are Lagrange multipliers, functions $\gamma'_{1}( \lambda )\le 0$ and $\gamma'_{1} ( \lambda )=0$ if $\langle B_{2},p_{0} ( \lambda) \rangle > \langle B_{2},V( \lambda ) \rangle,$ functions $\gamma'_{2}( \lambda )\ge 0$ and $\gamma'_{2} ( \lambda )=0$ if $\langle
B_{2} ,p_{0} ( \lambda) \rangle < \langle B_{2} ,U( \lambda ) \rangle.$

The following theorem  holds true.

\begin{thm}
Let the minimality condition (\ref{umova111_e_st.n_d}) hold true. The least favorable spectral densities $f_{0}(\lambda), $ $p_{0}(\lambda), $ in the classes $\md D_{f0}^k \times\md D_{Vp}^{Uk},k=1,2,3,4$  for the optimal linear forecasting of the functional  $A\vec{\xi}$ from observations of the vector sequence
$\vec{\zeta}(m)$, cointegrated with $\vec\xi(m)$ in terms of Definition \ref{ozn_coint}, at points  $m=-1,-2,\ldots$  are determined by equations
\eqref{eq_4_1f_coint}--\eqref{eq_4_1g_coint},  \eqref{eq_4_2f_coint}--\eqref{eq_4_2g_coint}, \eqref{eq_4_3f_coint}--\eqref{eq_4_3g_coint}, \eqref{eq_4_4f_coint}--\eqref{eq_4_4g_coint},
respectively,
the constrained optimization problem (\ref{minimax1}) with $g(\lambda):=|\beta^{(d)}(i\lambda)|^{-2}(p(\lambda)-\alpha^2f(\lambda))$ and restrictions  on densities from the corresponding classes $\md D_{f0}^k, \md D_{Vp}^{Uk},k=1,2,3,4$.  The minimax-robust spectral characteristic of the optimal estimate of the functional $A\vec{\xi}$ is determined by the formula (\ref{spectr A_co_e_st.n_d}).
\end{thm}

\subsection{Least favorable spectral density in classes $\md D_{\varepsilon}\times \md D_{1\delta}$}\label{set2}

Consider the prediction problem for the functional $A\vec{\xi}$
 which depends on unobserved values of a sequence $\vec\xi(m)$ with stationary increments based on observations of the sequence $\vec\xi(m)+\vec\eta(m)$ at points $m=-1,-2,\ldots$ under the condition that the sets of admissible spectral densities  $ \md D_{f\varepsilon }^{k},  \md D_{g1\delta }^{k}, k=1,2,3,4$ are defined as follows:

\[
\md D_{f\varepsilon }^{1}  =\bigg\{f(\lambda )\bigg|{\mathrm{Tr}}\,
[f(\lambda )]=(1-\varepsilon ) {\mathrm{Tr}}\,  [f_{1} (\lambda
)]+\varepsilon {\mathrm{Tr}}\,  [W(\lambda )],
\frac{1}{2\pi} \int _{-\pi}^{\pi}
\frac{|\chi_{\overline{\mu}}^{(d)}(e^{-i\lambda})|^2}{|\beta^{(d)}(i\lambda)|^2}
{\mathrm{Tr}}\,
\]
\[
\md D_{f\varepsilon }^{2}  =\bigg\{f(\lambda )\bigg|f_{kk} (\lambda)
=(1-\varepsilon )f_{kk}^{1} (\lambda )+\varepsilon w_{kk}(\lambda),
\frac{1}{2\pi} \int _{-\pi}^{\pi}
\frac{|\chi_{\overline{\mu}}^{(d)}(e^{-i\lambda})|^2}{|\beta^{(d)}(i\lambda)|^2}
f_{kk} (\lambda)d\lambda  =p_{k} , k=\overline{1,T}\bigg\};
\]
\begin{equation*}
\md D_{f\varepsilon }^{3} =\bigg\{f(\lambda )\bigg|\left\langle B_{1},f(\lambda )\right\rangle =(1-\varepsilon )\left\langle B_{1},f_{1} (\lambda )\right\rangle+\varepsilon \left\langle B_{1},W(\lambda )\right\rangle,
\frac{1}{2\pi}\int _{-\pi}^{\pi}
\frac{|\chi_{\overline{\mu}}^{(d)}(e^{-i\lambda})|^2}{|\beta^{(d)}(i\lambda)|^2}
\left\langle B_{1} ,f(\lambda )\right\rangle d\lambda =p\bigg\};
\end{equation*}
\begin{equation*}
\md D_{f\varepsilon }^{4}=\bigg\{f(\lambda )\bigg|f(\lambda)=(1-\varepsilon )f_{1} (\lambda )+\varepsilon W(\lambda ),
\frac{1}{2\pi } \int _{-\pi}^{\pi}
\frac{|\chi_{\overline{\mu}}^{(d)}(e^{-i\lambda})|^2}{|\beta^{(d)}(i\lambda)|^2}
f(\lambda )d\lambda=P\bigg\}.
\end{equation*}
\begin{equation*}
\md D_{g1\delta}^{1}=\left\{g(\lambda )\biggl|\frac{1}{2\pi} \int_{-\pi}^{\pi}
\left|{\rm{Tr}}(g(\lambda )-g_{1} (\lambda))\right|d\lambda \le \delta\right\};
\end{equation*}
\begin{equation*}
\md D_{g1\delta}^{2}=\left\{g(\lambda )\biggl|\frac{1}{2\pi } \int_{-\pi}^{\pi}
\left|g_{kk} (\lambda )-g_{kk}^{1} (\lambda)\right|d\lambda  \le \delta_{k}, k=\overline{1,T}\right\};
\end{equation*}
\begin{equation*}
\md D_{g1\delta}^{3}=\left\{g(\lambda )\biggl|\frac{1}{2\pi } \int_{-\pi}^{\pi}
\left|\left\langle B_{2} ,g(\lambda )-g_{1}(\lambda )\right\rangle \right|d\lambda  \le \delta\right\};
\end{equation*}
\begin{equation*}
\md D_{g1\delta}^{4}=\left\{g(\lambda )\biggl|\frac{1}{2\pi} \int_{-\pi}^{\pi}
\left|g_{ij} (\lambda )-g_{ij}^{1} (\lambda)\right|d\lambda  \le \delta_{i}^j, i,j=\overline{1,T}\right\}.
\end{equation*}

\noindent
Here  $f_{1} ( \lambda )$, $g_{1} ( \lambda )$ are fixed spectral densities, $W(\lambda)$ is an unknown spectral density, $p,  p_k, k=\overline{1,T}$, are given numbers, $P$ is a given positive-definite Hermitian matrices,
$\delta,\delta_{k},k=\overline{1,T}$, $\delta_{i}^{j}, i,j=\overline{1,T}$, are given numbers.

From the condition $0\in\partial\Delta_{\mathcal{D}}(f^0,g^0)$
we find the following equations which determine the least favourable spectral densities for these given sets of admissible spectral densities.

For the first set of admissible spectral densities  $\md D_{f\varepsilon}^{1}\times \md D_{g1\delta}^{1}$ we have equations
\begin{equation} \label{eq_5_1f}
\left(
{\me C}^{f0}_{\overline{\mu}}(e^{i \lambda})
\right)
\left(
{\me C}^{f0}_{\overline{\mu}}(e^{i \lambda})
\right)^{*}=
(\alpha_f^{2} +\gamma_1(\lambda ))
\left(\frac{|\chi_{\overline{\mu}}^{(d)}(e^{-i\lambda})|^2}{|\beta^{(d)}(i\lambda)|^2} (f^0(\lambda)+|\beta^{(d)}(i\lambda)|^2g^0(\lambda))\right)
^2,
\end{equation}
\begin{equation} \label{eq_5_1g}
\left(
{\me C}^{g0}_{\overline{\mu}}(e^{i \lambda})
\right)
\left(
{\me C}^{g0}_{\overline{\mu}}(e^{i \lambda})
\right)^{*}=
\beta^{2} \gamma_2( \lambda )\left(\frac{|\chi_{\overline{\mu}}^{(d)}(e^{-i\lambda})|^2}{|\beta^{(d)}(i\lambda)|^2}(f^0(\lambda)+|\beta^{(d)}(i\lambda)|^2g^0(\lambda))\right)^2,
\end{equation}
\begin{equation} \label{eq_5_1c}
\frac{1}{2 \pi} \int_{-\pi}^{ \pi}
\left|{\mathrm{Tr}}\, (g_0( \lambda )-g_{1}(\lambda )) \right|d\lambda =\delta,
\end{equation}

\noindent where $\alpha_f^{2}$, $ \beta^{2}$ are Lagrange multipliers,  the function $\gamma_1(\lambda )\le 0$ and $\gamma_1(\lambda )=0$ if ${\mathrm{Tr}}\,[f_{0} (\lambda )]>(1-\varepsilon ) {\mathrm{Tr}}\, [f_{1} (\lambda )]$, the function $\left| \gamma_2( \lambda ) \right| \le 1$ and
\[\gamma_2( \lambda )={ \mathrm{sign}}\; ({\mathrm{Tr}}\, (g_{0} ( \lambda )-g_{1} ( \lambda ))): \; {\mathrm{Tr}}\, (g_{0} ( \lambda )-g_{1} ( \lambda )) \ne 0.\]

For the second set of admissible spectral densities $\md D_{f\varepsilon}^{2}\times \md D_{g1\delta}^{2}$ we have equation
\begin{multline}   \label{eq_5_2f}
\left(
{\me C}^{f0}_{\overline{\mu}}(e^{i \lambda})
\right)
\left(
{\me C}^{f0}_{\overline{\mu}}(e^{i \lambda})
\right)^{*}=
\\
=
\left(\frac{|\chi_{\overline{\mu}}^{(d)}(e^{-i\lambda})|^2}{|\beta^{(d)}(i\lambda)|^2} (f^0(\lambda)+|\beta^{(d)}(i\lambda)|^2g^0(\lambda))\right)
\left\{(\alpha_{fk}^{2} +\gamma_{k}^1 (\lambda ))\delta _{kl} \right\}_{k,l=1}^{T}
\times
\\
\times
\left(\frac{|\chi_{\overline{\mu}}^{(d)}(e^{-i\lambda})|^2}{|\beta^{(d)}(i\lambda)|^2} (f^0(\lambda)+|\beta^{(d)}(i\lambda)|^2g^0(\lambda))\right),
\end{multline}
\begin{multline}   \label{eq_5_2g}
\left(
{\me C}^{g0}_{\overline{\mu}}(e^{i \lambda})
\right)
\left(
{\me C}^{g0}_{\overline{\mu}}(e^{i \lambda})
\right)^{*}=
\\
=
\left(\frac{|\chi_{\overline{\mu}}^{(d)}(e^{-i\lambda})|^2}{|\beta^{(d)}(i\lambda)|^2}(f^0(\lambda)+|\beta^{(d)}(i\lambda)|^2g^0(\lambda))\right)
\left \{ \beta_{k}^{2} \gamma^2_{k} ( \lambda ) \delta_{kl} \right \}_{k,l=1}^{T}
\times
\\
\times
\left(\frac{|\chi_{\overline{\mu}}^{(d)}(e^{-i\lambda})|^2}{|\beta^{(d)}(i\lambda)|^2}(f^0(\lambda)+|\beta^{(d)}(i\lambda)|^2g^0(\lambda))\right),
\end{multline}
\begin{equation} \label{eq_5_2c}
\frac{1}{2 \pi} \int_{- \pi}^{ \pi}  \left|g^0_{kk} ( \lambda)-g_{kk}^{1} ( \lambda ) \right| d\lambda =\delta_{k},
\end{equation}

\noindent where $\alpha_{fk}^{2}$, $\beta_{k}^{2}$ are Lagrange multipliers,  functions $\gamma_{k}^1(\lambda )\le 0$ and $\gamma_{k}^1 (\lambda )=0$ if $f_{kk}^{0}(\lambda )>(1-\varepsilon )f_{kk}^{1} (\lambda )$, functions $\left| \gamma^2_{k} ( \lambda ) \right| \le 1$ and
\[\gamma_{k}^2( \lambda )={ \mathrm{sign}}\;(g_{kk}^{0}( \lambda)-g_{kk}^{1} ( \lambda )): \; g_{kk}^{0} ( \lambda )-g_{kk}^{1}(\lambda ) \ne 0, \; k= \overline{1,T}.\]

For the third set of admissible spectral densities $\md D_{f\varepsilon}^{3}\times \md D_{g1\delta}^{3}$ we have equation
\begin{multline}   \label{eq_5_3f}
\left(
{\me C}^{f0}_{\overline{\mu}}(e^{i \lambda})
\right)
\left(
{\me C}^{f0}_{\overline{\mu}}(e^{i \lambda})
\right)^{*}=
\\
=
(\alpha_f^{2} +\gamma_1'(\lambda ))
\left(\frac{|\chi_{\overline{\mu}}^{(d)}(e^{-i\lambda})|^2}{|\beta^{(d)}(i\lambda)|^2} (f^0(\lambda)+|\beta^{(d)}(i\lambda)|^2g^0(\lambda))
\right)
 B_{1}^{\top}
 \times
\\
\times
\left(\frac{|\chi_{\overline{\mu}}^{(d)}(e^{-i\lambda})|^2}{|\beta^{(d)}(i\lambda)|^2} (f^0(\lambda)+|\beta^{(d)}(i\lambda)|^2g^0(\lambda))\right),
\end{multline}
\begin{multline}   \label{eq_5_3g}
\left(
{\me C}^{g0}_{\overline{\mu}}(e^{i \lambda})
\right)
\left(
{\me C}^{g0}_{\overline{\mu}}(e^{i \lambda})
\right)^{*}=
\\
=
\beta^{2} \gamma_2'( \lambda )
\left(\frac{|\chi_{\overline{\mu}}^{(d)}(e^{-i\lambda})|^2}{|\beta^{(d)}(i\lambda)|^2}(f^0(\lambda)+|\beta^{(d)}(i\lambda)|^2g^0(\lambda))\right)
B_{2}^{ \top}
\times
\\
\times
\left(\frac{|\chi_{\overline{\mu}}^{(d)}(e^{-i\lambda})|^2}{|\beta^{(d)}(i\lambda)|^2}(f^0(\lambda)+|\beta^{(d)}(i\lambda)|^2g^0(\lambda))\right),
\end{multline}
\begin{equation} \label{eq_5_3c}
\frac{1}{2 \pi} \int_{- \pi}^{ \pi}  \left| \left \langle B_{2}, g_0( \lambda )-g_{1} ( \lambda ) \right \rangle \right|d\lambda
= \delta,
\end{equation}
where $\alpha_f^{2}$,  $\beta^{2}$ are Lagrange multipliers,  function $\gamma_1' ( \lambda )\le 0$ and $\gamma_1' ( \lambda )=0$ if $\langle B_{1} ,f_{0} ( \lambda ) \rangle>(1- \varepsilon ) \langle B_{1} ,f_{1} ( \lambda ) \rangle$, function
$\left| \gamma_2' ( \lambda ) \right| \le 1$ and
\[\gamma_2' ( \lambda )={ \mathrm{sign}}\; \left \langle B_{2},g_{0} ( \lambda )-g_{1} ( \lambda ) \right \rangle : \; \left \langle B_{2},g_{0} ( \lambda )-g_{1} ( \lambda ) \right \rangle \ne 0.\]

For the fourth set of admissible spectral densities $\md D_{f\varepsilon}^{4}\times \md D_{g1\delta}^{4}$ we have equation
\begin{multline}  \label{eq_5_4f}
\left(
{\me C}^{f0}_{\overline{\mu}}(e^{i \lambda})
\right)
\left(
{\me C}^{f0}_{\overline{\mu}}(e^{i \lambda})
\right)^{*}=
\\
=
\left(\frac{|\chi_{\overline{\mu}}^{(d)}(e^{-i\lambda})|^2}{|\beta^{(d)}(i\lambda)|^2} (f^0(\lambda)+|\beta^{(d)}(i\lambda)|^2g^0(\lambda))\right)
(\vec{\alpha}_f\cdot \vec{\alpha}_f^{*}+\Gamma(\lambda))
\times
\\
\times
\left(\frac{|\chi_{\overline{\mu}}^{(d)}(e^{-i\lambda})|^2}{|\beta^{(d)}(i\lambda)|^2} (f^0(\lambda)+|\beta^{(d)}(i\lambda)|^2g^0(\lambda))\right),
\end{multline}
\begin{multline}  \label{eq_5_4g}
\left(
{\me C}^{g0}_{\overline{\mu}}(e^{i \lambda})
\right)
\left(
{\me C}^{g0}_{\overline{\mu}}(e^{i \lambda})
\right)^{*}=
\\
=
\left(\frac{|\chi_{\overline{\mu}}^{(d)}(e^{-i\lambda})|^2}{|\beta^{(d)}(i\lambda)|^2}(f^0(\lambda)+|\beta^{(d)}(i\lambda)|^2g^0(\lambda))\right)
\left \{ \beta_{ij}( \lambda ) \gamma_{ij} ( \lambda ) \right \}_{i,j=1}^{T}
\times
\\
\times
\left(\frac{|\chi_{\overline{\mu}}^{(d)}(e^{-i\lambda})|^2}{|\beta^{(d)}(i\lambda)|^2}(f^0(\lambda)+|\beta^{(d)}(i\lambda)|^2g^0(\lambda))\right),
\end{multline}
\begin{equation} \label{eq_5_4c}
\frac{1}{2 \pi} \int_{- \pi}^{ \pi} \left|g^0_{ij}(\lambda)-g_{ij}^{1}( \lambda ) \right|d\lambda = \delta_{i}^{j},
\end{equation}
where $\vec{\alpha}_f$, $ \beta_{ij}$ are Lagrange multipliers,  function $\Gamma(\lambda )\le 0$ and $\Gamma(\lambda )=0$ if $f_{0}(\lambda )>(1-\varepsilon )f_{1} (\lambda )$, functions $\left| \gamma_{ij} ( \lambda ) \right| \le 1$ and
\[
\gamma_{ij} ( \lambda )= \frac{g_{ij}^{0} ( \lambda )-g_{ij}^{1} (\lambda )}{ \left|g_{ij}^{0} ( \lambda )-g_{ij}^{1}(\lambda) \right|}: \; g_{ij}^{0} ( \lambda )-g_{ij}^{1} ( \lambda ) \ne 0, \; i,j= \overline{1,T}.
\]

The following theorem  holds true.

\begin{thm}
Let the minimality condition (\ref{umova11_e_st.n_d}) hold true. The least favorable spectral densities $f_{0}(\lambda)$, $g_{0}(\lambda)$ in classes
$\md D_{f\varepsilon}^{k}\times \md D_{g1\delta}^{k},k=1,2,3,4$ for the optimal linear extrapolation of the functional  $A\vec{\xi}$ from observations of the vector sequence $\vec{\xi}(m)+ \vec{\eta}(m)$ at points  $m=-1,-2,\ldots$ are determined by equations
\eqref{eq_5_1f} -- \eqref{eq_5_1c},  \eqref{eq_5_2f} -- \eqref{eq_5_2c}, \eqref{eq_5_3f} -- \eqref{eq_5_3c}, \eqref{eq_5_4f} -- \eqref{eq_5_4c},
respectively,
the constrained optimization problem (\ref{minimax1}) and restrictions  on densities from the corresponding classes
$ \md D_{f\varepsilon}^{k}, \md D_{g1\delta}^{k},k=1,2,3,4$.  The minimax-robust spectral characteristic of the optimal estimate of the functional $A\vec{\xi}$ is determined by the formula (\ref{spectr A_e_d}).
\end{thm}

Let the spectral densities $f(\lambda)$ and $g(\lambda)$ admit canonical factorizations (\ref{dd}), (\ref{fakt1}) and (\ref{fakt3}). Then we  derive the  following equation for the least favourable spectral densities.

For the first set of admissible spectral densities  $\md D_{f\varepsilon}^{1}\times \md D_{g1\delta}^{1}$ we have equations
\begin{equation} \label{eq_5_1f_fact}
\left(
{\me r^{0}_{\overline{\mu},f}(e^{i\lambda})}
\right)
\left(
{\me r^{0}_{\overline{\mu},f}(e^{i\lambda})}
\right)^{*}
=
(\alpha_f^{2} +\gamma_1(\lambda ))
(\Theta_{\overline{\mu}}(e^{-i\lambda}))^{\top}\overline{\Theta_{\overline{\mu}}(e^{-i\lambda})},
\end{equation}
\begin{equation} \label{eq_5_1g_fact}
\left(
{\me r^{0}_{\overline{\mu},g}(e^{i\lambda})}
\right)
\left(
{\me r^{0}_{\overline{\mu},g}(e^{i\lambda})}
\right)^{*}=
\beta^{2} \gamma_2( \lambda )(\Theta_{\overline{\mu}}(e^{-i\lambda}))^{\top}\overline{\Theta_{\overline{\mu}}(e^{-i\lambda})},
\end{equation}
\begin{equation} \label{eq_5_1c_fact}
\frac{1}{2 \pi} \int_{-\pi}^{ \pi}
\left|{\mathrm{Tr}}\, (g_0( \lambda )-g_{1}(\lambda )) \right|d\lambda =\delta,
\end{equation}

\noindent where $\alpha_f^{2}$, $ \beta^{2}$ are Lagrange multipliers,  the function $\gamma_1(\lambda )\le 0$ and $\gamma_1(\lambda )=0$ if ${\mathrm{Tr}}\,[f_{0} (\lambda )]>(1-\varepsilon ) {\mathrm{Tr}}\, [f_{1} (\lambda )]$, the function $\left| \gamma_2( \lambda ) \right| \le 1$ and
\[\gamma_2( \lambda )={ \mathrm{sign}}\; ({\mathrm{Tr}}\, (g_{0} ( \lambda )-g_{1} ( \lambda ))): \; {\mathrm{Tr}}\, (g_{0} ( \lambda )-g_{1} ( \lambda )) \ne 0.\]

For the second set of admissible spectral densities $\md D_{f\varepsilon}^{2}\times \md D_{g1\delta}^{2}$ we have equation
\begin{equation}   \label{eq_5_2f_fact}
\left(
{\me r^{0}_{\overline{\mu},f}(e^{i\lambda})}
\right)
\left(
{\me r^{0}_{\overline{\mu},f}(e^{i\lambda})}
\right)^{*}
=
(\Theta_{\overline{\mu}}(e^{-i\lambda}))^{\top}
\left\{(\alpha_{fk}^{2} +\gamma_{k}^1 (\lambda ))\delta _{kl} \right\}_{k,l=1}^{T}
\overline{\Theta_{\overline{\mu}}(e^{-i\lambda})},
\end{equation}
\begin{equation}   \label{eq_5_2g_fact}
\left(
{\me r^{0}_{\overline{\mu},g}(e^{i\lambda})}
\right)
\left(
{\me r^{0}_{\overline{\mu},g}(e^{i\lambda})}
\right)^{*}
=
(\Theta_{\overline{\mu}}(e^{-i\lambda}))^{\top}
\left \{ \beta_{k}^{2} \gamma^2_{k} ( \lambda ) \delta_{kl} \right \}_{k,l=1}^{T}
\overline{\Theta_{\overline{\mu}}(e^{-i\lambda})},
\end{equation}
\begin{equation} \label{eq_5_2c_fact}
\frac{1}{2 \pi} \int_{- \pi}^{ \pi}  \left|g^0_{kk} ( \lambda)-g_{kk}^{1} ( \lambda ) \right| d\lambda =\delta_{k},
\end{equation}

\noindent where $\alpha_{fk}^{2}$, $\beta_{k}^{2}$ are Lagrange multipliers,  functions $\gamma_{k}^1(\lambda )\le 0$ and $\gamma_{k}^1 (\lambda )=0$ if $f_{kk}^{0}(\lambda )>(1-\varepsilon )f_{kk}^{1} (\lambda )$, functions $\left| \gamma^2_{k} ( \lambda ) \right| \le 1$ and
\[\gamma_{k}^2( \lambda )={ \mathrm{sign}}\;(g_{kk}^{0}( \lambda)-g_{kk}^{1} ( \lambda )): \; g_{kk}^{0} ( \lambda )-g_{kk}^{1}(\lambda ) \ne 0, \; k= \overline{1,T}.\]

For the third set of admissible spectral densities $\md D_{f\varepsilon}^{3}\times \md D_{g1\delta}^{3}$ we have equation
\begin{equation}   \label{eq_5_3f_fact}
\left(
{\me r^{0}_{\overline{\mu},f}(e^{i\lambda})}
\right)
\left(
{\me r^{0}_{\overline{\mu},f}(e^{i\lambda})}
\right)^{*}
=
(\alpha_f^{2} +\gamma_1'(\lambda ))
(\Theta_{\overline{\mu}}(e^{-i\lambda}))^{\top}
 B_{1}
 \overline{\Theta_{\overline{\mu}}(e^{-i\lambda})},
\end{equation}
\begin{equation}   \label{eq_5_3g_fact}
\left(
{\me r^{0}_{\overline{\mu},g}(e^{i\lambda})}
\right)
\left(
{\me r^{0}_{\overline{\mu},g}(e^{i\lambda})}
\right)^{*}
=
\beta^{2} \gamma_2'( \lambda )
(\Theta_{\overline{\mu}}(e^{-i\lambda}))^{\top}
B_{2}
\overline{\Theta_{\overline{\mu}}(e^{-i\lambda})},
\end{equation}
\begin{equation} \label{eq_5_3c_fact}
\frac{1}{2 \pi} \int_{- \pi}^{ \pi}  \left| \left \langle B_{2}, g_0( \lambda )-g_{1} ( \lambda ) \right \rangle \right|d\lambda
= \delta,
\end{equation}
where $\alpha_f^{2}$,  $\beta^{2}$ are Lagrange multipliers,  function $\gamma_1' ( \lambda )\le 0$ and $\gamma_1' ( \lambda )=0$ if $\langle B_{1} ,f_{0} ( \lambda ) \rangle>(1- \varepsilon ) \langle B_{1} ,f_{1} ( \lambda ) \rangle$, function
$\left| \gamma_2' ( \lambda ) \right| \le 1$ and
\[\gamma_2' ( \lambda )={ \mathrm{sign}}\; \left \langle B_{2},g_{0} ( \lambda )-g_{1} ( \lambda ) \right \rangle : \; \left \langle B_{2},g_{0} ( \lambda )-g_{1} ( \lambda ) \right \rangle \ne 0.\]

For the fourth set of admissible spectral densities $\md D_{f\varepsilon}^{4}\times \md D_{g1\delta}^{4}$ we have equation
\begin{equation}  \label{eq_5_4f_fact}
\left(
{\me r^{0}_{\overline{\mu},f}(e^{i\lambda})}
\right)
\left(
{\me r^{0}_{\overline{\mu},f}(e^{i\lambda})}
\right)^{*}
=
(\Theta_{\overline{\mu}}(e^{-i\lambda}))^{\top}
(\vec{\alpha}_f\cdot \vec{\alpha}_f^{*}+\Gamma(\lambda))
\overline{\Theta_{\overline{\mu}}(e^{-i\lambda})},
\end{equation}
\begin{equation}  \label{eq_5_4g_fact}
\left(
{\me r^{0}_{\overline{\mu},g}(e^{i\lambda})}
\right)
\left(
{\me r^{0}_{\overline{\mu},g}(e^{i\lambda})}
\right)^{*}
=
(\Theta_{\overline{\mu}}(e^{-i\lambda}))^{\top}
\left \{ \beta_{ij}( \lambda ) \gamma_{ij} ( \lambda ) \right \}_{i,j=1}^{T}
\overline{\Theta_{\overline{\mu}}(e^{-i\lambda})},
\end{equation}
\begin{equation} \label{eq_5_4c_fact}
\frac{1}{2 \pi} \int_{- \pi}^{ \pi} \left|g^0_{ij}(\lambda)-g_{ij}^{1}( \lambda ) \right|d\lambda = \delta_{i}^{j},
\end{equation}
where $\vec{\alpha}_f$, $ \beta_{ij}$ are Lagrange multipliers,  function $\Gamma(\lambda )\le 0$ and $\Gamma(\lambda )=0$ if $f_{0}(\lambda )>(1-\varepsilon )f_{1} (\lambda )$, functions $\left| \gamma_{ij} ( \lambda ) \right| \le 1$ and
\[
\gamma_{ij} ( \lambda )= \frac{g_{ij}^{0} ( \lambda )-g_{ij}^{1} (\lambda )}{ \left|g_{ij}^{0} ( \lambda )-g_{ij}^{1}(\lambda) \right|}: \; g_{ij}^{0} ( \lambda )-g_{ij}^{1} ( \lambda ) \ne 0, \; i,j= \overline{1,T}.
\]

The following theorem  holds true.

\begin{thm}
The least favorable spectral densities $f_{0}(\lambda)$,  $g_{0}(\lambda)$  in the classes
$\md D_{f\varepsilon}^{k}\times \md D_{g1\delta}^{k},k=1,2,3,4$ for the optimal linear forecasting of the functional  $A\vec{\xi}$ from observations of the vector sequence $\vec{\xi}(m)+ \vec{\eta}(m)$ at points  $m=-1,-2,\ldots$ by  canonical factorizations (\ref{dd}), (\ref{fakt1}) and (\ref{fakt3}),
 equations
\eqref{eq_5_1f_fact} -- \eqref{eq_5_1c_fact},  \eqref{eq_5_2f_fact} -- \eqref{eq_5_2c_fact}, \eqref{eq_5_3f_fact} -- \eqref{eq_5_3c_fact}, \eqref{eq_5_4f_fact} -- \eqref{eq_5_4c_fact},
respectively,
the constrained optimization problem (\ref{simple_minimax1_e_st.n_d}) and restrictions  on densities from the corresponding classes
$ \md D_{f\varepsilon}^{k}, \md D_{g1\delta}^{k},k=1,2,3,4$.  The minimax-robust spectral characteristic of the optimal estimate of the functional $A\vec{\xi}$ is determined by the formula (\ref{spectr A_e_d}).
\end{thm}

\begin{thm}
 If the spectral density $g(\lambda)$ is known, the least favorable spectral density $f_{0}(\lambda)$ in the classes $\md D_{f\varepsilon}^{k}$, $k=1,2,3,4$    for the optimal linear forecasting of the functional  $A\vec{\xi}$ from observations of the sequence $\vec{\xi}(m)+ \vec{\eta}(m)$ at points  $m=-1,-2,\ldots$  is determined
by  canonical factorizations (\ref{dd}) and (\ref{fakt1})  ,
 equations
\eqref{eq_5_1f_fact},  \eqref{eq_5_2f_fact}, \eqref{eq_5_3f_fact}, \eqref{eq_5_4f_fact},
respectively,
 constrained optimization problem (\ref{simple_minimax3_e_st.n_d}) and restrictions  on density from the corresponding classes $\md D_{f\varepsilon}^{k}$, $k=1,2,3,4$.  The minimax-robust spectral characteristic of the optimal estimate of the functional $A\vec{\xi}$ is determined by the formula (\ref{spectr A_e_d}).
\end{thm}

\subsection{Least favorable spectral density for cointegrated vector sequences in classes $\md D_{\varepsilon}\times \md D_{1\delta}$}\label{set2_coint}

Consider the minimax forecasting problem for the functional $A\vec{\xi}$
 which depends on unobserved values of a vector sequence $\vec\xi(m)$ with GM increments based on observations of the Vector sequence $\vec\zeta(m)$, cointegrated with $\vec\xi(m)$ in terms of Definition \ref{ozn_coint}, at points $m=-1,-2,\ldots$ for  the sets of admissible spectral densities  $ \md D_{f\varepsilon}^{k}$, $k=1,2,3,4$, defined in Subsection \ref{set2} and  $\md D_{p1\delta}^{k}$, $k=1,2,3,4$, defined as follows:
\begin{equation*}
\md D_{p1\delta}^{1}=\left\{g(\lambda )\biggl|\frac{1}{2\pi} \int_{-\pi}^{\pi}\frac{|\chi_{\overline{\mu}}^{(d)}(e^{-i\lambda})|^2}{|\beta^{(d)}(i\lambda)|^2}
\left|{\rm{Tr}}(p(\lambda )-p_{1} (\lambda))\right|d\lambda \le \delta\right\};
\end{equation*}
\begin{equation*}
\md D_{p1\delta}^{2}=\left\{p(\lambda )\biggl|\frac{1}{2\pi } \int_{-\pi}^{\pi}\frac{|\chi_{\overline{\mu}}^{(d)}(e^{-i\lambda})|^2}{|\beta^{(d)}(i\lambda)|^2}
\left|p_{kk} (\lambda )-p_{kk}^{1} (\lambda)\right|d\lambda  \le \delta_{k}, k=\overline{1,T}\right\};
\end{equation*}
\begin{equation*}
\md D_{p1\delta}^{3}=\left\{p(\lambda )\biggl|\frac{1}{2\pi } \int_{-\pi}^{\pi}\frac{|\chi_{\overline{\mu}}^{(d)}(e^{-i\lambda})|^2}{|\beta^{(d)}(i\lambda)|^2}
\left|\left\langle B_{2} ,p(\lambda )-p_{1}(\lambda )\right\rangle \right|d\lambda  \le \delta\right\};
\end{equation*}
\begin{equation*}
\md D_{p1\delta}^{4}=\left\{p(\lambda )\biggl|\frac{1}{2\pi} \int_{-\pi}^{\pi}\frac{|\chi_{\overline{\mu}}^{(d)}(e^{-i\lambda})|^2}{|\beta^{(d)}(i\lambda)|^2}
\left|p_{ij} (\lambda )-p_{ij}^{1} (\lambda)\right|d\lambda  \le \delta_{i}^j, i,j=\overline{1,T}\right\}.
\end{equation*}

\noindent
Here   $p_{1} ( \lambda )$ is a fixed spectral density,
$\delta,\delta_{k},k=\overline{1,T}$, $\delta_{i}^{j}, i,j=\overline{1,T}$, are given numbers.

The condition $0\in\partial\Delta_{\mathcal{D}}(f^0,g^0)$
implies the following equations which determine the least favourable spectral densities for these given sets of admissible spectral densities.

For the first set of admissible spectral densities  $\md D_{f\varepsilon}^{1}\times \md D_{p1\delta}^{1}$ we have equations
\begin{multline} \label{eq_5_1f_coint}
\left({\me C^{f0}_{\overline{\mu},\alpha}(e^{i\lambda})}
\right)
\left(
{\me C^{f0}_{\overline{\mu},\alpha}(e^{i\lambda})}
\right)^{*}
-\alpha^2\left(
{\me C}^{p0}_{\overline{\mu},\alpha}(e^{i \lambda})
\right)
\left(
{\me C}^{p0}_{\overline{\mu},\alpha}(e^{i \lambda})
\right)^{*}=
(\alpha_f^{2} +\gamma_1(\lambda ))
\left(\frac{|\chi_{\overline{\mu}}^{(d)}(e^{-i\lambda})|^2}{|\beta^{(d)}(i\lambda)|^2} p^0(\lambda)\right)^2,
\end{multline}
\begin{equation} \label{eq_5_1g_coint}
\left(
{\me C}^{p0}_{\overline{\mu},\alpha}(e^{i \lambda})
\right)
\left(
{\me C}^{p0}_{\overline{\mu},\alpha}(e^{i \lambda})
\right)^{*}
=
\beta^{2} \gamma_2( \lambda )\left(\frac{|\chi_{\overline{\mu}}^{(d)}(e^{-i\lambda})|^2}{|\beta^{(d)}(i\lambda)|^2} p^0(\lambda)\right)^2,
\end{equation}
\begin{equation} \label{eq_5_1c_coint}
\frac{1}{2 \pi} \int_{-\pi}^{ \pi}\frac{|\chi_{\overline{\mu}}^{(d)}(e^{-i\lambda})|^2}{|\beta^{(d)}(i\lambda)|^2}
\left|{\mathrm{Tr}}\, (p_0( \lambda )-p_{1}(\lambda )) \right|d\lambda =\delta,
\end{equation}

\noindent where $\alpha_f^{2}$, $ \beta^{2}$ are Lagrange multipliers,  the function $\gamma_1(\lambda )\le 0$ and $\gamma_1(\lambda )=0$ if ${\mathrm{Tr}}\,[f_{0} (\lambda )]>(1-\varepsilon ) {\mathrm{Tr}}\, [f_{1} (\lambda )]$, the function $\left| \gamma_2( \lambda ) \right| \le 1$ and
\[\gamma_2( \lambda )={ \mathrm{sign}}\; ({\mathrm{Tr}}\, (p_{0} ( \lambda )-p_{1} ( \lambda ))): \; {\mathrm{Tr}}\, (p_{0} ( \lambda )-p_{1} ( \lambda )) \ne 0.\]

For the second set of admissible spectral densities $\md D_{f\varepsilon}^{2}\times \md D_{p1\delta}^{2}$ we have equation
\begin{multline}   \label{eq_5_2f_coint}
\left({\me C^{f0}_{\overline{\mu},\alpha}(e^{i\lambda})}
\right)
\left(
{\me C^{f0}_{\overline{\mu},\alpha}(e^{i\lambda})}
\right)^{*}
-\alpha^2\left(
{\me C}^{p0}_{\overline{\mu},\alpha}(e^{i \lambda})
\right)
\left(
{\me C}^{p0}_{\overline{\mu},\alpha}(e^{i \lambda})
\right)^{*}=
\\
=
\left(\frac{|\chi_{\overline{\mu}}^{(d)}(e^{-i\lambda})|^2}{|\beta^{(d)}(i\lambda)|^2} p^0(\lambda)\right)
\left\{(\alpha_{fk}^{2} +\gamma_{k}^1 (\lambda ))\delta _{kl} \right\}_{k,l=1}^{T}
\left(\frac{|\chi_{\overline{\mu}}^{(d)}(e^{-i\lambda})|^2}{|\beta^{(d)}(i\lambda)|^2} p^0(\lambda)\right),
\end{multline}
\begin{equation}   \label{eq_5_2g_coint}
\left(
{\me C}^{p0}_{\overline{\mu},\alpha}(e^{i \lambda})
\right)
\left(
{\me C}^{p0}_{\overline{\mu},\alpha}(e^{i \lambda})
\right)^{*}=
\left(\frac{|\chi_{\overline{\mu}}^{(d)}(e^{-i\lambda})|^2}{|\beta^{(d)}(i\lambda)|^2} p^0(\lambda)\right)
\left \{ \beta_{k}^{2} \gamma^2_{k} ( \lambda ) \delta_{kl} \right \}_{k,l=1}^{T}
\left(\frac{|\chi_{\overline{\mu}}^{(d)}(e^{-i\lambda})|^2}{|\beta^{(d)}(i\lambda)|^2} p^0(\lambda)\right),
\end{equation}
\begin{equation} \label{eq_5_2c_coint}
\frac{1}{2 \pi} \int_{- \pi}^{ \pi} \frac{|\chi_{\overline{\mu}}^{(d)}(e^{-i\lambda})|^2}{|\beta^{(d)}(i\lambda)|^2} \left|p^0_{kk} ( \lambda)-p_{kk}^{1} ( \lambda ) \right| d\lambda =\delta_{k},
\end{equation}

\noindent where $\alpha _{fk}^{2}$, $\beta_{k}^{2}$ are Lagrange multipliers,  functions $\gamma_{k}^1(\lambda )\le 0$ and $\gamma_{k}^1 (\lambda )=0$ if $f_{kk}^{0}(\lambda )>(1-\varepsilon )f_{kk}^{1} (\lambda )$, functions $\left| \gamma^2_{k} ( \lambda ) \right| \le 1$ and
\[\gamma_{k}^2( \lambda )={ \mathrm{sign}}\;(p_{kk}^{0}( \lambda)-p_{kk}^{1} ( \lambda )): \; p_{kk}^{0} ( \lambda )-p_{kk}^{1}(\lambda ) \ne 0, \; k= \overline{1,T}.\]

For the third set of admissible spectral densities $\md D_{f\varepsilon}^{3}\times \md D_{p1\delta}^{3}$ we have equation
\begin{multline}   \label{eq_5_3f_coint}
\left({\me C^{f0}_{\overline{\mu},\alpha}(e^{i\lambda})}
\right)
\left(
{\me C^{f0}_{\overline{\mu},\alpha}(e^{i\lambda})}
\right)^{*}
-\alpha_f^2\left(
{\me C}^{p0}_{\overline{\mu},\alpha}(e^{i \lambda})
\right)
\left(
{\me C}^{p0}_{\overline{\mu},\alpha}(e^{i \lambda})
\right)^{*}=
\\
=
(\alpha^{2} +\gamma_1'(\lambda ))
\left(\frac{|\chi_{\overline{\mu}}^{(d)}(e^{-i\lambda})|^2}{|\beta^{(d)}(i\lambda)|^2} p^0(\lambda)\right)
 B_{1}^{\top}
\left(\frac{|\chi_{\overline{\mu}}^{(d)}(e^{-i\lambda})|^2}{|\beta^{(d)}(i\lambda)|^2} p^0(\lambda)\right),
\end{multline}
\begin{equation}   \label{eq_5_3g_coint}
\left(
{\me C}^{p0}_{\overline{\mu},\alpha}(e^{i \lambda})
\right)
\left(
{\me C}^{p0}_{\overline{\mu},\alpha}(e^{i \lambda})
\right)^{*}=
\beta^{2} \gamma_2'( \lambda )
\left(\frac{|\chi_{\overline{\mu}}^{(d)}(e^{-i\lambda})|^2}{|\beta^{(d)}(i\lambda)|^2} p^0(\lambda)\right)
B_{2}^{ \top}
\left(\frac{|\chi_{\overline{\mu}}^{(d)}(e^{-i\lambda})|^2}{|\beta^{(d)}(i\lambda)|^2} p^0(\lambda)\right),
\end{equation}
\begin{equation} \label{eq_5_3c_coint}
\frac{1}{2 \pi} \int_{- \pi}^{ \pi} \frac{|\chi_{\overline{\mu}}^{(d)}(e^{-i\lambda})|^2}{|\beta^{(d)}(i\lambda)|^2} \left| \left \langle B_{2}, p_0( \lambda )-p_{1} ( \lambda ) \right \rangle \right|d\lambda
= \delta,
\end{equation}
where $\alpha_f^{2}$,  $\beta^{2}$ are Lagrange multipliers,  function $\gamma_1' ( \lambda )\le 0$ and $\gamma_1' ( \lambda )=0$ if $\langle B_{1} ,f_{0} ( \lambda ) \rangle>(1- \varepsilon ) \langle B_{1} ,f_{1} ( \lambda ) \rangle$, function
$\left| \gamma_2' ( \lambda ) \right| \le 1$ and
\[\gamma_2' ( \lambda )={ \mathrm{sign}}\; \left \langle B_{2},p_{0} ( \lambda )-p_{1} ( \lambda ) \right \rangle : \; \left \langle B_{2},p_{0} ( \lambda )-p_{1} ( \lambda ) \right \rangle \ne 0.\]

For the fourth set of admissible spectral densities $\md D_{f\varepsilon}^{4}\times \md D_{p1\delta}^{4}$ we have equation
\begin{multline} \label{eq_5_4f_coint}
\left({\me C^{f0}_{\overline{\mu},\alpha}(e^{i\lambda})}
\right)
\left(
{\me C^{f0}_{\overline{\mu},\alpha}(e^{i\lambda})}
\right)^{*}
-\alpha^2\left(
{\me C}^{p0}_{\overline{\mu},\alpha}(e^{i \lambda})
\right)
\left(
{\me C}^{p0}_{\overline{\mu},\alpha}(e^{i \lambda})
\right)^{*}=
\\
=
\left(\frac{|\chi_{\overline{\mu}}^{(d)}(e^{-i\lambda})|^2}{|\beta^{(d)}(i\lambda)|^2} p^0(\lambda)\right)
(\vec{\alpha}_f\cdot \vec{\alpha}_f^{*}+\Gamma(\lambda))
\left(\frac{|\chi_{\overline{\mu}}^{(d)}(e^{-i\lambda})|^2}{|\beta^{(d)}(i\lambda)|^2} p^0(\lambda)\right),
\end{multline}
\begin{equation}  \label{eq_5_4g_coint}
\left(
{\me C}^{p0}_{\overline{\mu},\alpha}(e^{i \lambda})
\right)
\left(
{\me C}^{p0}_{\overline{\mu},\alpha}(e^{i \lambda})
\right)^{*}=
\left(\frac{|\chi_{\overline{\mu}}^{(d)}(e^{-i\lambda})|^2}{|\beta^{(d)}(i\lambda)|^2} p^0(\lambda)\right)
\left \{ \beta_{ij}( \lambda ) \gamma_{ij} ( \lambda ) \right \}_{i,j=1}^{T}
\left(\frac{|\chi_{\overline{\mu}}^{(d)}(e^{-i\lambda})|^2}{|\beta^{(d)}(i\lambda)|^2} p^0(\lambda)\right),
\end{equation}
\begin{equation} \label{eq_5_4c_coint}
\frac{1}{2 \pi} \int_{- \pi}^{ \pi}\frac{|\chi_{\overline{\mu}}^{(d)}(e^{-i\lambda})|^2}{|\beta^{(d)}(i\lambda)|^2}\left|p^0_{ij}(\lambda)-p_{ij}^{1}( \lambda ) \right|d\lambda = \delta_{i}^{j},
\end{equation}
where $\vec{\alpha}_f$, $ \beta_{ij}$ are Lagrange multipliers,  function $\Gamma(\lambda )\le 0$ and $\Gamma(\lambda )=0$ if $f_{0}(\lambda )>(1-\varepsilon )f_{1} (\lambda )$, functions $\left| \gamma_{ij} ( \lambda ) \right| \le 1$ and
\[
\gamma_{ij} ( \lambda )= \frac{p_{ij}^{0} ( \lambda )-p_{ij}^{1} (\lambda )}{ \left|p_{ij}^{0} ( \lambda )-p_{ij}^{1}(\lambda) \right|}: \; p_{ij}^{0} ( \lambda )-p_{ij}^{1} ( \lambda ) \ne 0, \; i,j= \overline{1,T}.
\]

The following theorem  holds true.

\begin{thm}
Let the minimality condition (\ref{umova111_e_st.n_d}) hold true. The least favorable spectral densities $f_{0}(\lambda)$, $p_{0}(\lambda)$ in classes
$\md D_{f\varepsilon}^{k}\times \md D_{p1\delta}^{k},k=1,2,3,4$ for the optimal linear filtering of the functional  $A\vec{\xi}$ from observations of the vector sequence $\vec{\zeta}(m)$ cointegrated with $\vec\xi(m)$ in terms of Definition \ref{ozn_coint}, at points  $m=-1,-2,\ldots$ are determined by equations
\eqref{eq_5_1f_coint} -- \eqref{eq_5_1c_coint},  \eqref{eq_5_2f_coint}--\eqref{eq_5_2c_coint}, \eqref{eq_5_3f_coint} -- \eqref{eq_5_3c_coint}, \eqref{eq_5_4f_coint} -- \eqref{eq_5_4c_coint},
respectively,
the constrained optimization problem (\ref{minimax1}) with $g(\lambda):=|\beta^{(d)}(i\lambda)|^{-2}(p(\lambda)-\alpha^2f(\lambda))$ and restrictions  on densities from the corresponding classes
$ \md D_{f\varepsilon}^{k}, \md D_{p1\delta}^{k},k=1,2,3,4$.  The minimax-robust spectral characteristic of the optimal estimate of the functional $A\vec{\xi}$ is determined by the formula (\ref{spectr A_co_e_st.n_d}).
\end{thm}

\begin{section}{Conclusions}

In this article, we dealt with stochastic sequences with periodically stationary GM increments introduced in \cite{Luz_Mokl_extra_GMI}.
We give a definition of one class of vector seasonally cointegrated sequences related to stationary GM increment.
These non-stationary stochastic sequences combine  periodic structure of covariation functions of sequences as well as integrating one.

We derived solutions of the  forecasting problem for the linear functionals constructed from the  unobserved values of a sequence with periodically stationary GM increments.
Estimates are based on observations of the sequence with a periodically stationary noise.
We obtained the estimates   by representing the sequence under investigation as a vector-valued sequence with stationary GM increments. Based on the  solutions for these type of sequences, we solved the corresponding problem for the defined class of seasonally cointegrated vector sequences.
The problem is investigated in the case of spectral certainty, where spectral densities of sequences are exactly known.
In this case  we propose an approach based on the Hilbert space projection method.
We derive formulas for calculating the spectral characteristics and the mean-square errors of the optimal estimates of the functionals.
In the case of spectral uncertainty where the spectral densities are not exactly known while, instead, some sets of admissible spectral densities are specified, the minimax-robust method is applied.
We propose a representation of the mean square error in the form of a linear
functional in $L_1$ with respect to spectral densities, which allows
us to solve the corresponding constrained optimization problem and
describe the minimax (robust) estimates of the functionals. Formulas
that determine the least favorable spectral densities and minimax
(robust) spectral characteristic of the optimal linear estimates of
the functionals are derived  for a wide list of specific classes
of admissible spectral densities.

These least favourable spectral density matrices are solutions of the optimization problem $\Delta_{\mathcal{D}}(f,g)=\widetilde{\Delta}(f,g)+ \delta(f,g|\mathcal{D}_f\times
\mathcal{D}_g)\to\inf,$
 where $\delta(f,g|\mathcal{D}_f\times
\mathcal{D}_g)$ is the indicator function of the set
$\mathcal{D}=\mathcal{D}_f\times\mathcal{D}_g$.
 Solution $(f^0,g^0)$ to this unconstrained optimisation problem is characterized by the condition $0\in
\partial\Delta_{\mathcal{D}}(f^0,g^0)$, where
$\partial\Delta_{\mathcal{D}}(f^0,g^0)$ is the subdifferential of the functional $\Delta_{\mathcal{D}}(f,g)$ at point $(f^0,g^0)\in \mathcal{D}=\mathcal{D}_f\times\mathcal{D}_g$.
 This condition makes it possible to find the least favourable spectral densities in some special classes of spectral densities.
These are: classes $D_0$ of densities with the moment restrictions,
 classes $D_{1\delta}$ which describe the ``$\delta$-neighborhood''\, models in the space $L_1$ of a fixed bounded spectral density,
 classes $D_{\varepsilon}$ which describe the  ``$\varepsilon$-contaminated''\, models of a fixed bounded spectral density,
 classes $ D_V^U$ which describe the ``strip'' models of spectral densities.

\end{section}

\end{document}